\documentclass[11pt, reqno, a4paper]{amsart}
\usepackage{amsmath,amssymb,amsthm}
\usepackage{color}
\usepackage{comment}
\usepackage[utf8]{inputenc}
\usepackage[dvipdfmx]{graphicx}
\usepackage{enumerate}	
\theoremstyle{definition}
 \newtheorem{dfn}{Definition}
 \newtheorem{remark}[dfn]{Remark}

\theoremstyle{plain}
 \newtheorem{thm}[dfn]{Theorem}
 \newtheorem{prop}[dfn]{Proposition}
 \newtheorem{lem}[dfn]{Lemma}
 \newtheorem{cor}[dfn]{Corollary}
 \newtheorem{assump}[dfn]{Assumption}

\numberwithin{equation}{section}
\newcommand{\bn}{{\bold n}}

\newcommand{\bk}{{\bold k}}
\newcommand{\bm}{{\bold m}}
\newcommand{\bu}{{\bold u}}

\newcommand{\bv}{{\bold v}}
\newcommand{\bw}{{\bold w}}

\newcommand{\bff}{{\bold f}}
\newcommand{\bA}{{\bold A}}

\newcommand{\bG}{{\bold G}}
\newcommand{\bH}{{\bold H}}
\newcommand{\bI}{{\bold I}}

\newcommand{\bM}{{\bold M}}
\newcommand{\bN}{{\bold N}}

\newcommand{\bU}{{\bold U}}

\newcommand{\DV}{{\rm Div}\,}
\newcommand{\dv}{{\rm div}}

\newcommand{\BD}{{\Bbb D}}
\newcommand{\BR}{{\Bbb R}}
\newcommand{\BC}{{\Bbb C}}

\newcommand{\BN}{{\Bbb N}}

\newcommand{\BI}{{\Bbb I}}

\newcommand{\BZ}{{\Bbb Z}}

\newcommand{\CA}{{\mathcal A}}
\newcommand{\CB}{{\mathcal B}}
\newcommand{\CC}{{\mathcal C}}
\newcommand{\CD}{{\mathcal D}}

\newcommand{\CF}{{\mathcal F}}

\newcommand{\CJ}{{\mathcal J}}

\newcommand{\CL}{{\mathcal L}}
\newcommand{\CM}{{\mathcal M}}

\newcommand{\CR}{{\mathcal R}}
\newcommand{\CS}{{\mathcal S}}
\newcommand{\CT}{{\mathcal T}}
\newcommand{\CH}{{\mathcal H}}
\newcommand{\CO}{{\mathcal O}}
\newcommand{\CP}{{\mathcal P}}
\newcommand{\CQ}{{\mathcal Q}}
\newcommand{\CU}{{\mathcal U}}
\newcommand{\CV}{{\mathcal V}}
\newcommand{\CW}{{\mathcal W}}

\newcommand{\CZ}{{\mathcal Z}}

\newcommand{\bg}{{\bold g}}
\newcommand{\bh}{{\bold h}}
\newcommand{\pd}{\partial}

\newcommand{\supp}{{\rm supp\,}}
\newcommand{\HS}{\BR^N_+}

\renewcommand{\d}{{\rm d}}

\begin{document}
\title[]{On the $L_1$-maximal regularity in the study of free boundary problem 
for the compressible fluid flows}

\author[]{Yuko Enomoto}
\address{(Y. Enomoto)  Department of Mathematical Sciences, 
Shibaura Institute of Technology \\
Fukasaku 307, Minuma-ku, Saitama-shi 337-8570, Japan }		
\email{e-yuko@shibaura-it.ac.jp} 
\author[]{Yoshihiro Shibata}
\address{(Y. Shibata) 
{ Professor Emeritus of Waseda University \\
3-4-5 Ohkubo Shinjuku-ku, Tokyo, 169-8555, Japan. \\
Adjunct faculty member in the Department of Mechanical Engineering 
and	Materials Science, University of Pittsburgh, USA
}}
\email{yshibata325@gmail.com}

\thanks{{MSC Classification of 2020: Primary: 76N10; Secondary: 35Q30. }
\\
\quad \enskip
The first author is supported by JST SPRING, Grant Number JPMJSP2128.  
The second author is  partially supported by JSPS KAKENHI 
Grant Number 23K22405}


\maketitle

\begin{abstract}
In this paper, we consider the Stokes equations with non-homogeneous
free boundary conditions, which is obtained by the linearization procedure
of the free boundary problem of the Navier-Stokes equations describing
the viscous compressible fluid flows.  We prove the $L_1$ maximal regularity
of solutions to this Stokes equations. This is an extension result of 
$L_p$-$L_q$ maximal regularity result obtained by D.~G\"otz and Y.~Shibata 
\cite{GS14} to the $L_1$ in time maximal regularity case. 
\end{abstract}
\section{Introduction} 
\subsection{Problem and Result}
Let $\Omega$ be a domain in the $N$ dimensional Euclidean space $\BR^N$, 
whose boundary $\pd\Omega$ is a $C^3$ compact hypersurface. Let $\bn$
denote the unit outer normal to $\pd\Omega$. 
In this paper, we consider the Stokes equations with free boundary conditions, which 
read as
\begin{equation}\label{Eq:Linear}\left\{\begin{aligned}
\pd_t\rho+\eta_0\dv\bu  = F&&\quad&\text{in $\Omega\times(0, \infty)$}, \\
\eta_0\pd_t\bu - \DV(\alpha\BD(\bu) + (\beta-\alpha)\dv\bu\BI
-P'(\eta_0)\rho\BI)  = \bG&&\quad&\text{in $\Omega\times(0, \infty)$}, \\
(\alpha\BD(\bu) + (\beta-\alpha)\dv\bu\BI-P'(\eta_0)\rho\BI)\bn
= \bH&&\quad&\text{on $\pd\Omega\times(0, \infty)$}, \\
(\rho, \bu)|_{t=0}=(\rho_0, \bu_0)&&\quad&\text{in $\Omega$}.
\end{aligned}\right.\end{equation}
Here, $\rho$ and $\bu$ denote unknown density field and velocity field, while 
$F$, $\bG$, $\bH$, $\rho_0$ and $\bu_0$ are given right members and initial data.
The $\alpha$ and $\beta$ denote viscosity constants such that $\alpha > 0$ and 
$\alpha+\beta>0$.  The $P(\rho)$ denote a pressure which is a $C^\infty$ function of
$\rho \in (0, \infty)$ such that $P'(\rho) > 0$.  Let $\rho_*$ be a positive constant
describing the mass density of the reference body and $\eta_0(x) = \rho_* + \tilde\eta_0(x)$,
where $\tilde\eta_0(x)$ belongs to some Besov space. Throughout the paper, we assume that 
there exist two positive constants $\rho_1 < \rho_2$ such that 
\begin{equation}
\rho_1 < \rho_* < \rho_2, \quad \rho_1 < \eta_0(x) < \rho_2, \quad
\rho_1 < P'(\rho_*) < \rho_2, \quad \rho_1 < P'(\eta_0(x)) < \rho_2.
\label{assump:0}
\end{equation}
Moreover, $\BD(\bu) = \nabla\bu + \nabla\bu^\top$, where $A^\top$ denotes the
transposed $A$ and $\nabla \bu$ the gradient of $N$ vector of function
$\bu=(u_1, \ldots, u_N)$, and $\dv \bu = \sum_{j=1}^N \pd_ju_j$, where $\pd_j=\pd/\pd x_j$,
and $\BI$ denotes the $N\times N$ identity matrix. 
For any $N\times N$ matrix of functions $A=(A_{jk})$, $\DV A$ denotes the 
$N$ vector of functions whose $i$-th component is $\sum_{j=1}^N\pd_jA_{ij}$.

The equations in \eqref{Eq:Linear} are obtained by linearization at $\rho=\eta_0$ and $\bu=0$
of the Navier-Stokes equations describing the viscous compressible fluid motion without surface
tension in a
time dependent domain $\Omega_t$ with free boundary conditions on its boundary $\pd\Omega_t$,
which read as
\begin{equation}\label{Eq:Nonlinear}\left\{\begin{aligned}
\pd_t\rho+\dv(\rho\bu)  = 0&&\quad&\text{in $\Omega_t$}, \\
\rho(\pd_t\bu + \bu\cdot\nabla\bu) - \DV(\alpha\BD(\bu) + (\beta-\alpha)\dv\bu\BI
-P'(\rho)\rho\BI)  = \bG&&\quad&\text{in $\Omega_t$}, \\
(\alpha\BD(\bu) + (\beta-\alpha)\dv\bu\BI-P'(\rho)\rho\BI)\bn_t
= P(\rho_*)\bn_t&&\quad&\text{on $\pd\Omega_t$}, \\
V_{\pd\Omega_t} = \bn_t\cdot\bu&&\quad&\text{on $\pd\Omega_t$}, \\
(\rho, \bu)|_{t=0} = (\rho_0, \bu_0)&&\quad&\text{in $\Omega$}.
\end{aligned}\right.\end{equation}
for $t \in (0, T)$.  Here, $\bn_t$ denotes the unit outer normal to $\pd\Omega_t$ and 
$V_{\pd\Omega_t}$ denotes the evolution speed of $\pd\Omega_t$ in the normal direction.
Thus, $V_{\pd\Omega_t}=\bn_t\cdot\bu$ denotes the non-slip condition, which is a kinematic
boundary condition. \par
The purpose of this paper is to prove the $L_1$ in time maximal regularity of solutions
$\rho$ and $\bu$ of equations \eqref{Eq:Linear}.  The local well-posedness in $\Omega$
and the global well-posedness for small initial data in a bounded domain for
\eqref{Eq:Nonlinear} will be proved
as a direct application of $L_1$ maximal regularity of this paper, and the global
well-posedness for small initial data in exterior domains for \eqref{Eq:Nonlinear}
will be proved some combination of the $L_1$ maximal regularity and decay properties of 
solutions to \eqref{Eq:Linear} with $\eta_0=\rho_*$.  The topics for the nonlinear equations
\eqref{Eq:Nonlinear} will be treated in a forthcoming paper. 
\par
To state the main result, we introduce the opdrator $\Lambda^{1/2}_\gamma$ which is
defined by 
$$\Lambda^{1/2}_\gamma f = \CL^{-1}[\lambda^{1/2}\CL[f](\lambda)](t)
= \frac{1}{2\pi}\int_{\BR}e^{\lambda t}\lambda^{1/2}\CL[f](\lambda)\, d \tau
\quad (\lambda = \gamma + i\tau \in \BC).$$
 Here $\CL$ and $\CL^{-1}$ denote the 
Laplace transform and the inverse Laplace transform defined by
$$\CL[f](\lambda) = \int_{\BR}e^{-\lambda t}f(t)\,d t, \quad
\CL^{-1}[f](t) = \frac{1}{2\pi} \int_{\BR} e^{\lambda t}f(\tau)\,d\tau
\quad(\lambda = \gamma + i\tau \in \BC).
$$
\par
In this paper, we shall prove the following theorem.
\begin{thm}\label{thm:main.1}
Assume that the following conditions \thetag1 or \thetag2 holds. 
\begin{itemize}
\item[\thetag1] If $\eta_0(x) = \rho_*$, then $1 < q < \infty$ and $-1+1/q < s < 1/q$. 
\item[\thetag2] If $\tilde\eta_0(x)\not\equiv0$, then 
$N-1 < q < 2N$ and $-1+N/q \leq  s < 1/q$. 
\end{itemize}
Then, there exists a large constant $\gamma_0 > 0$ such that for any initial data 
$(\rho_0, \bu_0) \in \CH^s_{q,1}(\Omega):=B^{s}_{q,1}(\Omega)^N\times B^{s+1}_{q,1}(\Omega)^N$, 
and right members $F$, $\bG$, and $\bH$ satisfying the 
conditions:
\begin{align*}
e^{-\gamma t}F &\in L_1(\BR, B^{s+1}_{q,1}(\Omega)), \quad 
e^{-\gamma t}\bG \in L_1(\BR, B^s_{q,1}(\Omega)^N), \\
e^{-\gamma t}\bH &\in W^{1/2}_1(\BR, B^s_{q,1}(\Omega)^N) \cap
L_1(\BR, B^{s+1}_{q,1}(\Omega)^N),
\end{align*}
for some $\gamma \geq \gamma_0$, 
then problem \eqref{Eq:Linear} admits unique solutions $\rho$ and $\bu$ with
\begin{align*}
e^{-\gamma t}\rho &\in W^1_1((0, \infty), B^{s+1}_{q,1}(\Omega)), \\
e^{-\gamma t}\bu & \in L_1((0, \infty), B^s_{q,1}(\Omega)^N) \cap W^1_1((0, \infty), 
B^{s+2}_{q,1}(\Omega)^N)
\end{align*}
possessing the estimate:
\begin{align*}
&\|e^{-\gamma t} \rho\|_{W^1_1((0, \infty), B^{s+1}_{q,1}(\Omega))}
+ \|e^{-\gamma t}\pd_t\bu\|_{L_1((0, \infty), B^s_{q,1}(\Omega))}
+ \|e^{-\gamma t}\Lambda^{1/2}_\gamma \bu\|_{L_1((0, \infty), B^{s+1}_{q,1}(\Omega))}\\
&\qquad + \|e^{-\gamma t}\bu\|_{L_1((0, \infty), B^{s+2}_{q,1}(\Omega))}\\
&\leq C(\|\rho_0\|_{B^{s+1}_{q,1}(\Omega)} + \|\bu_0\|_{B^s_{q,1}(\Omega)}
+ \|e^{-\gamma t}F\|_{L_1(\BR, B^{s+1}_{q,1}(\Omega))}
+ \|e^{-\gamma t}\bG\|_{L_1(\BR, B^s_{q,1}(\Omega))} \\
&\qquad + 
\|e^{-\gamma t}\Lambda_\gamma^{1/2}\bH\|_{L_1(\BR, B^s_{q,1}(\Omega))}
+ \|e^{-\gamma t}\bH\|_{L_1(\BR, W^{s+1}_{q,1}(\Omega))}).
\end{align*}
Here, $C$ depends on $\gamma_0$ but is independent of $\gamma$ whenever 
$\gamma \geq \gamma_0$.
\end{thm}
Theorem \ref{thm:main.1} follows from the spectral properties of solutions to the 
generalized resolvent problem: 
\begin{equation}\label{Eq:Resolvent}\left\{\begin{aligned}
\lambda\rho+\eta_0\dv\bv  = f&&\quad&\text{in $\Omega$}, \\
\eta_0\lambda\bv - \DV(\alpha\BD(\bv) + (\beta-\alpha)\dv\bv\BI
-P'(\eta_0)\rho\BI)  = \bg&&\quad&\text{in $\Omega$}, \\
(\alpha\BD(\bv) + (\beta-\alpha)\dv\bv\BI-P'(\eta_0)\rho\BI)\bn
= \bh&&\quad&\text{on $\pd\Omega$}
\end{aligned}\right.\end{equation}
where the spectral parameter $\lambda$ runs through the parabolic sector 
$\Sigma_{\epsilon, \lambda_0}$ for $0 < \epsilon < \pi/2$ and large $\lambda_0>0$,
where 
\begin{equation*}
\Sigma_\epsilon = \{\lambda \in \BC\setminus\{0\} \mid |\arg\lambda| \leq \pi-\epsilon\},
\quad 
\Sigma_{\epsilon, \lambda_0} = \{\lambda \in \Sigma_\epsilon \mid |\lambda| \geq \lambda_0\}.
\end{equation*}
And the spectral properties of solutions to equations \eqref{Eq:Resolvent} will be drived
as perturbation for large $\lambda$ 
from the spectral properties of solutions to Lam\'e equations, which read as
\begin{equation}\label{Eq:Lame}\left\{\begin{aligned}
\eta_0\lambda\bv - \DV(\alpha\BD(\bv) + (\beta-\alpha)\dv\bv\BI)
  = \bg&&\quad&\text{in $\Omega$}, \\
(\alpha\BD(\bv) + (\beta-\alpha)\dv\bv\BI)\bn
= \bh&&\quad&\text{on $\pd\Omega$}.
\end{aligned}\right.\end{equation}
In fact, setting $\rho = \lambda^{-1}(f-\eta_0\dv\bu)$ in the first equation of
\eqref{Eq:Resolvent} and inserting this formula 
into the second equations in \eqref{Eq:Resolvent} imply that 
\begin{equation*}
\left\{\begin{aligned}
\eta_0\lambda\bv - \DV(\alpha\BD(\bv) + (\beta-\alpha)\dv\bv\BI 
+ \lambda^{-1}P'(\eta_0)\eta_0\dv\bu\BI)
 & = \bg-\lambda^{-1} \nabla (P'(\eta_0)f)&\quad&\text{in $\Omega$}, \\
(\alpha\BD(\bv) + (\beta-\alpha)\dv\bv\BI+\lambda^{-1}P'(\eta_0)\dv\bu\BI)\bn
& = \bh-\lambda^{-1}P'(\eta_0)f&\quad&\text{on $\pd\Omega$}.
\end{aligned}\right.\end{equation*}
From this observation, we see that \eqref{Eq:Resolvent} can be regarded as a perturbation 
from \eqref{Eq:Lame} for large $\lambda$.  Thus, the main part of this paper is devoted 
to analysis of the spectral properties of solutions to \eqref{Eq:Lame}.
\subsection{Short History} Mathematical studies on the compressible Navier-Stokes
equations started the uniqueness results by \cite{Grafi, Itaya71, Nash62, Serrin, VH}.
After these works, concerning the study of the local and global well-posedness 
for the Cauchy problem and the initial boundary value problem with non-slip conditions
(Dirichlet conditions), many  studies have been done. The local well-posedness
has been studied by Tani \cite{Tani77} in the H\"older space and by Solonnikov 
\cite{Sol80} in the Sobolev-Slobodetskii spaces. Danchin \cite{D10} improved 
Solonnikov's result in the critical space in the Cauchy problem case. 
Matsumura and Nishida \cite{MN80, MN83} made a breakthrough
in proving the global well-posedness for small initial data using the energy method. 
This result was extended to the optimal regularity of initial data in the $L_2$ space 
 by Kawashita \cite{K}.  Kobayashi and Shibata \cite{KS1} improved the decay properties
of solutions in the exterior domains combining the energy method and $L_p$-$L_q$ decay 
properties of solutions to the linearized equations, where the condition:
$1 < p \leq 2 \leq q \leq \infty$ is assumed.  In the  no restrictions of
exponents case, 
so called the diffusion wave properties has been studied by Hopf and Zumbrun \cite{HZ}
and Liu and Wang \cite{LW}.  Kobayashi and Shibata \cite{KS1, KS2}  improved their results.
In the half space case, the decay properties were studied by Kagei and Kobayashi \cite{KK02, 
KK05}. 
The global well-posedness results were extensively studied  in the energy spaces of exterior
domains by \cite{SE18, ST04, VZ86, WT} 
and in the critical space of the whole space by \cite{AP07, CMZ, DX, CD10, H11, H, O}. 
Valli \cite{Vali83} , Kagei and Tsuda \cite{Kagei-Tsuda},
and Tsuda \cite{Tsuda16}
studied time periodic solutions in the bounded domains and in $\BR^N$, respectively. 
The analytic semigroup approach  using  Lagrange 
transformation was started by Str\"ohmer \cite{Str:89, Str:90} and the $L_p$-$L_q$ 
maximal regularity of Stokes semigroup  in general domains was proved by 
Enomoto and Shibata \cite{ES13} and the global well-posedness in the maximal $L_p$-$L_q$
regularity class was proved by Shibata \cite{S22} in exterior domains. 
More references can be found in a survey paper \cite{SE18}, and also references of the papers
mentioned above.  \par
Concerning the viscous compressible fluid motion without surface tension, 
the local well-posedness of the free boundary problem without surface tension
for compressible viscous fluid flow in the multi-dimensional case was first proved by
Secchi and Valli \cite{SV83} in the $L_2$ framework and by Tani \cite{T81}
in the H\"older spaces,
respectively. The gaseous stars case was established by Secchi \cite{PS90.1, 
PS90.2}.  Enomoto, von Below and Shibata \cite{EBS14} proved the local
well-posedness in the maximal $L_p$-$L_q$ regularity class, where they used the 
$L_p$-$L_q$ maximal regularity theory for the linearized equations established by
G\"otz abd Shibata \cite{GS14}.  The global well-posedness has been established
by Shibata \cite{YS16} in the bounded domain case.  More recently, Shibata and 
Zhang \cite{SZ23} prove the global well-posedness for small data in 
exterior domains combining the $L_p$-$L_q$ maximal regularity theorem with 
$L_p$-$L_q$ decay properties of solutions to the linearized equations 
obtained in Shibata and Zhang \cite{SZ22}. 

Concerning the viscous compressible fluid motion with surface tension, 
the local well-posedness has been studied by Denisova and Solonnikov 
\cite{DS03.1, DS03.2} in the H\"older space framework, 
and the global well-posedness was proved by
Solonnikov and Tani \cite{ST91}, Zadrzynska and Zajaczkowski \cite{ZZ99},
and Zajaczkowski \cite{Z94} in the bounded domains and $L_2$ framework,
where  the reference domain is close to a ball, respectively, the initial density is close
to a positive constant and the initial velocity is small. Concerning the 
$L_p$-$L_q$ approach to the surface tension problem has been done by
Zhang \cite{XZ20}\par

On the other hand, $L_1$ in time maximal regularity theorem was first investigated 
by Danchin and Mucha \cite{DM09} for the viscous incompressible fluid 
with non-slip boundary conditions in exterior domains. 
 The global well-posedness for the free boundary problem of the 
viscous incompressible fluid in the half-space
 has been studied  first by Danchin et al \cite{DHMT}
using the extension of Da Prato and Grisvard \cite{DG} $L_1$ in time maximal regularity
theorem for abstract evolution equations under some regularity assumption 
of  initial data.  The restriction in \cite{DHMT} was removed by 
Ogawa and Shimizu \cite{OS} and Shibata and Watanabe \cite{SW}.
The argument based on Littlewood-Payley decomposition in the half space
plays an essential role in \cite{OS}, while  the real interpolation method 
due to  Shibata \cite{S23} was used in \cite{SW}. 
  The method due to Shibata \cite{S23} is based on the spectral analysis
of linearized equations and so $L_1$ in time maximal regularity can be
obtained for initial boundary value problems of parabolic and hyperbolic-
parabolic system of equations appearing in the mathematical fluid mechanics,
like Navier-Stokes equations both in the viscous incompressible fluid  case and
in the viscous compressible fluid  case. 
In fact, the $L_1$ in time maximal regularity was obtained for the 
compressible Navier-Stokes equations with non-slip boundary conditions by
Kuo \cite{Kuo23} and Kuo and Shibata \cite{KS23}. 
In this paper, we consider the $L_1$ maximal regularity in the 
free boundary condition case. 
This paper  is a continuation of the study of free boundary problem for the compressible 
viscous fluid flows in $L_p$-$L_q$ maximal regularity framework due to \cite{GS14, EBS14}
\par
{\bf Why is the $L_1$ maximal regularity important ?}   
Assuming that the velocity field $\bu$ is expected to be found as an element of the function space 
$L_p((0, T), W^2_q(\Omega)) \cap W^1_p((0,  T), L_q(\Omega))$, the trace theorem
tells us that $\bu|_{t=0} \in B^{2(1-1/p)}_{q,p}(\Omega)$.  Thus, $L_1$ maximal regularity 
gives the best order regularity  class  of the initial data.  Moreover, without surface tension case,
we usually use the Lagrange transform, because of the lack of regularity of functions representing
the free surface. Thus, $L_1$ maximal regularity is the best space to treat the Lagrange
transform. 

\section{Preparations for latter sections}

 \subsection{Symbols used throughout the paper}
Let us explain the symbols used in this paper. Let $\BR$, $\BN$, and $\BC$ be the set of all real, natural, 
complex numbers, respectively, while let $\BZ$ be the set of all integers. 
Set $\BN_0 := \BN \cup \{0\}$. For multi-index $\kappa=(\kappa_1, \ldots, \kappa_N)
\in \BN_0^N$ 
and $x=(x_1, \ldots, x_N) \in \BR^N$, $\pd_x^\kappa = \pd^\kappa=
\pd^{|\kappa|}/\pd^{\kappa_1}x_1\cdots \pd^{\kappa_N}x_N$ stands for standard 
partial derivatives with respect to $x$ 
of order $\kappa$.  For the dual variable $\xi=(\xi_1, \ldots, \xi_N)\in \BR^N$, 
$D^\kappa_\xi = \pd^{|\kappa|}/\pd^{\kappa_1}\xi_1\cdots \pd^{\kappa_N}\xi_N$. 
For differentiations, we use symbols
$\nabla f= \{\pd^\kappa f \mid |\kappa|=1\}, \bar\nabla f = \{\pd^\kappa f \mid
|\kappa| \leq 1\}$, $\nabla^2f = \{\pd^\kappa f \mid |\kappa|=2\}$, 
$\bar\nabla^2f = \{\pd^\kappa f \mid |\kappa| \leq 2\}$. \par 
Let $\HS$ and $\BR^N_0$ denote the half space and its boundary defined by
$$\HS = \{x=(x_1, \ldots, x_N) \in \BR^N \mid x_N > 0\},\quad
\BR^N_0 = \{x=(x_1, \ldots, x_N) \in \BR^N \mid x_N=0\},$$
and $\bn_0=(0, \ldots, 0, -1)$ denotes the unit outer normal to $\pd\HS$.
For $N\in \BN$ and a Banach space $X$, let $\CS(\BR^N;X)$ be the Schwartz class of 
$X$-valued rapidly decreasing functions on $\BR^N$. Let
 $\mathcal{S'}(\BR^N;X)$ denote 
the space of $X$-valued tempered distributions, which means the set of all 
continuous linear mappings from $\CS(\BR^N)$ to  $X$.  The Fourier transform 
$\CF [f]$  and the inverse Fourier transform
  $\CF^{-1}_\xi[g]$ are defined by 
\begin{equation*}
\CF [f](\xi):=\int_{\BR^{N}}f(x) e^{-ix\cdot\xi }\,dx, \qquad 
\CF_{\xi}^{-1} [g] (x) :=\frac{1}{(2\pi)^N}\int_{\BR^{N}}g(\xi) e^{ix\cdot\xi }\,d\xi,
\end{equation*}
respectively. In addition,  the partial Fourier transform
$\CF'[f(\,\cdot\,, x_N)] = \hat f(\xi', x_N)$ 
and the partial inverse Fourier transform $\CF^{-1}_{\xi'}$ are defined by
\begin{align*}
\CF'[f(\,\cdot\,, x_N)] (\xi') &:= \hat f(\xi', x_N) 
= \int_{\BR^{N-1}} f(x', x_N) e^{-ix'\cdot\xi'} \,dx',\\
\CF^{-1}_{\xi'}[g (\,\cdot\,, x_N)](x') &:= \frac{1}{(2\pi)^{N-1}}\int_{\BR^{N-1}}
g(\xi', x_N) e^{ix'\cdot\xi'}\,d\xi',
\end{align*}
respectively, 
where  $x'=(x_1,\cdots,x_{N-1}) \in \BR^{N-1}$ and 
$\xi' = (\xi_1, \cdots, \xi_{N-1})\in\BR^{N-1}$. 
The Laplace transform $\CL[f](\lambda)$ and inverse
Laplace transform $\CL^{-1}[g](t)$ are defined by 
$$\CL[f](\lambda) = \int_{\BR} e^{-\lambda t}f(t)\,dt, \quad
\CL^{-1}[g](t) = \frac{1}{2\pi i}\int_{\BR} e^{\lambda t} g(\lambda)\,d\tau
\quad (\lambda = \gamma + i\tau),
$$
respectively. 
\par
For a  domain $D$ and a Banach space $X$,  $L_p(D, X)$, $W^m_p(D, X)$ $(m \in \BN$) 
and $W^s_{p}(D, X)$ ($s>0$, $s\not\in \BN$)
stand for standard $X$ valued Lebesgue spaces, 
 Sobolev spaces and Sobolev-Slobodetskii spaces,
while $\|\cdot\|_{L_p(D,X)}$, $\|\cdot\|_{W^m_p(D,X)}$ and $\|\cdot\|_{W^s_p(D, X)}$
denote their norms. In particular, we write
\begin{align*}\|e^{-\gamma t}f\|_{L_p((0, T), X)} 
&= \Bigr(\int_0^T (e^{-\gamma t}\|f(t)\|_X)^p \,dt\Bigr)^{1/p}, \\
\|e^{-\gamma t}\Lambda^{1/2}_\gamma f\|_{L_p(\BR, X)} &=
\Bigr(\int_{\BR} (e^{-\gamma t}\|\Lambda^{1/2}_\gamma f(t)
\|_X)^p \,dt\Bigr)^{1/p}.
\end{align*}
We set 
$$W^{1/2}_p(\BR, X) = \{f \mid \|e^{-\gamma t}f \|_{W^{1/2}_p(\BR, X)}
= \|e^{-\gamma t}f\|_{L_p(\BR, X)} + \|e^{-\gamma t}\Lambda^{1/2}_\gamma
f\|_{L_p(\BR, X)} < \infty\}. 
$$
When $X = \BR^N$, we omit $X=\BR^N$, namely, we write
$L_p(D)$, $W^m_p(D)$, $W^s_p(D)$, $\|\cdot\|_{L_p(D)}$, $\|\cdot\|_{W^m_p(D)}$ 
and $\|\cdot\|_{W^s_q(D)}$. 
In particular, $W^0_p(D) = L_p(D)$ for the notational simplicity. \par
For any domain $D$, the functional space $\CH^\nu_{q,1}(D)$ 
of data $(\bg, \bh)$ for spectral problems
is defined by 
\begin{equation*}
\CH^\nu_{q,1}(D) = B^\nu_{q,1}(D)^N\times B^{\nu+1}_{q,1}(D)^N, 
\quad \|(\bg, \bh)\|_{\CH^\nu_{q,1}(D)} = \|\bg\|_{B^\nu_{q,1}(D)} + \|\bh\|_{B^{\nu+1}_{q,1}(D)}
\enskip\text{for $(\bg, \bh) \in \CH^\nu_{q,1}(D)$}.
\end{equation*}
For $(\bg, \bh) \in \CH^\nu_{q,1}(D)$, $H_i=\{H_{ij} \mid j=1, \ldots, N\}$ ($i=1,2$) 
are respective corresponding variables to 
$\bg=(g_1, \ldots, g_N)$ and $\lambda^{1/2}\bh=(\lambda^{1/2}h_1, \ldots, \lambda^{1/2}h_N)$.
And, $\{H_{3ij} \mid i,j=1, \ldots, N\}$ are corresponding variables to $\pd_i h_j$, and 
set $H_3=\{H_{3ij} \mid i, j=1, \ldots, N\}$.   Let
$m(N) = 2N+N^2$. Set 
$H=(H_1, H_2, H_3)$, which is an $m(N)$ vector of functions.\par 
For a domain $D$ in $\BR^N$ and 
$N \ge 2$, we set $(\bff, \bg)_D = \int_D \bff (x) \cdot \bg (x) \,dx$ 
for $N$-vector functions $\bff$ and $\bg$ on $D$, 
where we will write $(\bff, \bg) = (\bff, \bg)_{D}$ for short if there is no confusion.
\par
For a Banach space $X$, $\|\cdot\|_X$ denotes its
norm.  For Banach spaces $X$ and $Y$, $\CL(X, Y)$ denotes the set of all bounded
linear operators from $X$ into $Y$, and we write $\CL(X)=\CL(X,X)$.  
Let ${\rm Hol}\, (U, X)$ denote the set of 
all $X$ valued holomorphic functions defined on  a domain $U \subset \BC$.
The $n$ product space of $X$ is written as $X^n
=\{x =(x_1, \ldots, x_n) \mid x_i \in X (i=1, \ldots, n)\}$, 
while its norm is denoted by
$\|x\|_X= \sum_{i=1}^n\|x_i\|_X$. 
$X \hookrightarrow Y$ means that $X$ is continuously imbedded into $Y$, 
that is $X \subset Y$ and $\|x\|_Y \leq C\|x\|_X$ with some constant $C$. 
\par
For any interpolation couple $(X, Y)$ of Banach spaces $X$ and $Y$, 
the maps $(X, Y) \to (X, Y)_{\theta, p}$ and $(X, Y) \to (X, Y)_{[\theta]}$ 
denote the real interpolation functor for each $\theta \in (0, 1)$ and 
$p \in [1, \infty]$ and 
the complex interpolation functor for each $\theta \in (0,1)$, respectively.
By $C > 0$ we will often denote a generic constant 
that does not depend on the quantities at stake. 
 And, by $C_{a, b, c, \cdots}$ we denote 
generic constants depending on the quantities $a$, $b$, $c, \cdots$.  $C$ and 
$C_{a, b, c, \cdots}$ may change from line to line. 
\subsection{Definition of Besov spaces and some properties}

To define  Besov pace $B^s_{q,1}$,  we  introduce Littlewood-Paley decomposition.
Let $\phi \in \CS (\BR^N)$ with $\supp \phi = \{ \xi \in \BR^N \mid 1 
\slash 2 \le \lvert \xi \rvert \le 2 \}$ such that
$\sum_{k \in \BZ} \phi (2^{- k} \xi) = 1$ for all $\xi \in \BR^N \setminus \{0\}$. 
Then, define
\begin{equation*}
\phi_k := \CF^{- 1}_\xi [\phi (2^{- k} \xi)] \quad (k \in \BZ), \qquad 
\psi = 1 - \sum_{k \in \BN} \phi_k.
\end{equation*}
For $1 \le p, q \le \infty$ and $s \in \BR$ we denote
\begin{equation*}
\lVert f \rVert_{B^s_{p, q} (\BR^N)} := 
\left\{\begin{aligned}
& \lVert \psi * f \rVert_{L_p (\BR^N)} + \bigg(\sum_{k \in \BN} 
\Big(2^{s k} \lVert \phi_k * f \rVert_{L_p (\BR^N)} \Big)^q \bigg)^{1 \slash q}
& \quad & \text{if $1\le q < \infty$}, \\
& \lVert \psi * f \rVert_{L_p (\BR^N)} + \sup_{k \in \BN} 
\Big(2^{s k} \lVert \phi_k * f \rVert_{L_p (\BR^N)} \Big)
& \quad & \text{if $q=\infty$}.
\end{aligned}\right.
\end{equation*}
Here, $f \, * \, g$ means the convolution between $f$ and $g$.
The inhomogeneous Besov spaces $B^s_{p, q} (\BR^N)$ is 
defined as the sets of all 
$f \in \CS' (\BR^N)$ such that $\lVert f \rVert_{B^s_{p, q} (\BR^N)} < \infty$. 
In particular, we define $B^s_{q, \infty-}(\BR^N)$ by
$$B^s_{p, \infty-}(\BR^N) = \{f \in B^s_{p,\infty}(\BR^N) \mid 
\lim_{j\to\infty} \, 2^{js}\|\phi_k*f\|_{L_p(\BR^N)} = 0\}.
$$
In this paper, we use the following conventions: $r < \infty- < \infty$ for $r \in \BR$ and 
$1/\infty-=1/\infty=0$.  In particular, the H\"older conjugate of $1$ and $\infty-$
are  read as $1' = \infty-$ and $ (\infty-)'=1$.
\par
For any domain $D$ in $\BR^N$,  $B^s_{q,1}(D)$ is defined by the restriction of 
$B^s_{q,1}(\BR^N)$, that is 
\begin{gather*}B^s_{q,1}(D) = \{f \in \CD'(D) \mid \text{there exists a $g\in B^s_{q,1}(\BR^N)$
such that $g|_{D} = f$}\}, \\
\|f\|_{B^s_{q,1}(D)} = \inf \{\|g\|_{B^s_{q,1}(\BR^N)} \mid g \in B^s_{q,1}(\BR^N), \enskip g|_{D} = f\}.
\end{gather*}
Here, $\CD'(D)$ denotes the set of all distributions on $D$ and $g|_{D}$ denotes the 
restriction of $g$ to $D$. Notice that $W^s_p(D) = B^s_{p,p}(D)$,
which are called a Sobolev-Slobodeckij spaces.  In particular, 
if $s$ is a non-negative integer, $W^s_q$ is a usual Sobolev space. \par 
It is well-known that if $D$ satisfies the cone property, then 
$B^s_{p,q}(D)$ may be \textit{characterized} by means of real interpolation.
In fact, for $-\infty<s_0<s_1<\infty$, $1<p<\infty$, $1\le q\le\infty$ and $0<\theta<1$, 
it follows that
\begin{equation*}
B^{\theta s_0+(1-\theta)s_1}_{p,q}(D)=\left(W^{s_0}_p(D), 
W^{s_1}_p(D)\right)_{\theta,q}
\end{equation*}
cf. \cite[Theorem 8]{M74},  \cite[Theorem 2.4.2]{Tbook78}. 
If $D$ is a uniform $C^1$ domain, then $D$ satisfies the cone property.

\subsection{The  estimate for the product of 2 functions and 
some composite  functions using Besov norms} 
The following lemma is concerned the estimate of product of two functions
using Besov norms. 
\begin{lem}\label{lem:prod} 
Let $D$ be a uniform $C^3$ domain whose boundary is compact hypersurface. 
Let  $N-1 < q  < 2N$ and $-1+N/q \leq s < 1/q$.  Then, 
for any $u \in B^s_{q,1}(D)$ and 
$v \in B^{N/q}_{q,\infty}(D) \cap L_\infty(D)$, there holds
\begin{equation}\label{prod.1}
\|uv\|_{B^s_{q,1}(D)} \leq C_{D,s,q,r}\|u\|_{B^s_{q,1}(D)}
\|v\|_{B^{N/q}_{q, \infty}(D) \cap L_\infty(D)}.
\end{equation}
\end{lem}
\begin{remark} \thetag1 The condition that $N-1 < q$ is required to obtain 
$-1+N/q < 1/q$. \\
\thetag2 We know that 
$\|v\|_{B^{N/q}_{q, \infty}(D) \cap L_\infty(D)} \leq C\|v\|_{B^{N/q}_{q,1}(D)}$. 
\end{remark}
\begin{proof} By the Abidi-Paicu estimate \cite{AP07} and the Haspot estimate
\cite{H11}, when $2< q < \infty$ and 
$-N/q < s < N/q$ or when $1 \leq q < 2$ and $-N/q' < s < N/q$,  the estimate \eqref{prod.1}
holds. When $2 \leq q < \infty$, the condition: $-N/q <-1+N/q$ implies that 
$q < 2N$.  When $1 \leq q < 2$, the condition: $-N/q' < -1 + N/q$ implies that 
$N \geq 1$.  Thus, we have Lemma \ref{lem:prod}.
\end{proof}
We now use the following lemma for the Besov norm estimate of composite functions cf. 
\cite[Proposition 2.4]{H11} and \cite[Theorem 2.87]{BCD}.
\begin{lem}\label{lem:Hasp} Let $1 < q < \infty$. 
Let $I$ be an open interval of $\BR$.  Let $\omega>0$ and let $\tilde\omega$ be the smallest
integer such that $\tilde\omega \geq \omega$. Let $F:I \to \BR$ satisfy $F(0) = 0$ and 
$F' \in BC^{\tilde\omega}(I, \BR)$ that is $F' \in C^{\tilde\omega}(I, \BR)$ and 
${\displaystyle \|F'\|_{BC^{\tilde\omega}(I,\BR)} :=\sum_{\ell=0}^{\tilde\omega} 
\sup_{t \in I} \| \partial^{\ell}_t F' \|_{L_{\infty}(\BR)} < \infty}$. 
Assume that $v\in B^\omega_{q,1}$
has valued in $J \subset I$.  Then, 
$F(v) \in B^\omega_{q,1}$ and there exists a constant $C$ depending only on 
$\nu$, $I$, $J$, and $N$, such that 
$$\|F(v)\|_{B^\omega_{q,1}} \leq C(1 + \|v\|_{L_\infty})^{\tilde\omega}
\|F'\|_{BC^{\tilde\omega}(I,\BR)}
\|v\|_{B^\omega_{q,1}}.$$
\end{lem}
\subsection{Fourier multiplier theorems in $\BR^N$} 
To estimate solutions of the equations in $\BR^N$, we use the following 
Fourier multiplier theorem of Mihlin - H\"ormander type \cite{Mih, Hor}, 
which is stated as follows: 
Let $m(\xi)$ be a $C^\infty(\BR^N)$ function such that for any multi-index 
$\kappa \in \BN_0^{N}$ there exists a constant $C_\kappa$ such that 
\begin{equation*}
\label{sumbol.1}|D_\xi^\kappa m(\xi)| \leq C_\kappa|\xi|^{-|\kappa|}.
\end{equation*}
We call $m$ a multiplier symbol of order $0$.  Set 
$[m] = \max_{|\kappa| \leq N}C_{\kappa}$.  For any multiplier symbol of order $0$,
we define an operator $T_m$ by 
$$T_mf = \CF^{-1}[m \CF[f]].$$  
We call $T_m$ the Fourier multiplier with symbol $m$. Then, for any $1 < p < \infty$,
there exists a constant $C_p$ depending on $p$ such that there holds
\begin{equation*}\label{fmt.1}
\|T_mf\|_{L_p(\BR^N)} \leq C_p[m]\|f\|_{L_p(\BR^N)}.
\end{equation*}
We extend this result to the Besov space case as follows: 
\begin{prop}\label{FMT} 
Let $1<q<\infty$, $1 \leq r \leq \infty$ and $s\in\BR$. 
Let $m(\xi)$ be a multiplier symbol of order $0$ and let $T_m$ be the Fourier multiplier 
with symbol $m$.  Then, there exists a constant $C_{q,r}$ depending on $q$ and $r$ 
such that for any $f \in B^s_{q,r}(\BR^n)$, there holds
$$\|T_mf\|_{B^s_{q,r}(\BR^N)} \leq C_{q,r}[m]\|f\|_{B^s_{q,r}(\BR^N)}.$$
\end{prop}
\begin{proof}
First of all, we recall the definition of Besov spaces in Subsection 2.2. 
Since $\phi_k * (T_mf ) 
= \CF^{-1}[m \CF[\phi_k*f]]$, by the standard Fourier multiplier
theorem of Mihlin-H\"ormander type, we have
$$\|\phi_k *(T_mf) \|_{L_q(\BR^N)} \leq C[m]\|\phi_k*f\|_{L_q(\BR^N)}.$$
Similarly, we have
$$\|\psi *( T_mf)\|_{L_q(\BR^N)} \leq C[m]\|\psi*f\|_{L_q(\BR^N)}.$$
Thus, by the definition of the Besov norm, we have
$$\|T_mf\|_{B^s_{q,r}(\BR^N)} \leq C[m]\|f\|_{B^s_{q,r}(\BR^N)}.$$
This completes the proof of Proposition \ref{FMT}. 
\end{proof}
\subsection{Symbol classes and estimates of the integral operators in $\HS$}
To state main results of this subsection, we introduce a symbol class.  In the following, 
let $m(\lambda, \xi')$ be a function defined on $\Sigma_{\epsilon, \lambda_0}\times
(\BR^{N-1}\setminus\{0\})$ such that for each $\xi' \in \BR^{N-1}\setminus\{0\}$ 
$m(\lambda, \xi')$ is  holomorphic with respect to $\lambda \in \Sigma_{\epsilon, \lambda_0}$ 
and for each $\lambda \in \Sigma_{\epsilon, \lambda_0}$ $m(\lambda, \xi')$ is a
$C^\infty$ function with respect to $\xi' \in
\BR^{N-1}\setminus\{0\}$.
Let $\ell \in \BZ$. We say that  $m(\lambda, \xi')$ is an order $\ell$  symbol 
if for any $\kappa' \in \BN_0^{N-1}$
and $\lambda \in \Sigma_{\epsilon}$ there exists a constant $C_{\kappa'}$ being
depending on $\kappa'$, $\epsilon$, $\lambda_0$ and $\ell$ such that 
\begin{equation*}
|D_{\xi}^{\kappa'} m(\lambda, \xi')| \leq C_{\kappa'}(|\lambda|^{1/2}+|\xi'|)^{\ell-|\kappa'|}.
\end{equation*}
Let $\bM_\ell$ be the set of all order $\ell$ symbols.
Let 
$$\|m\| = \max_{|\kappa'| \leq N} C_{\kappa'}.$$
Let $$A=\sqrt{(\alpha+\beta)^{-1}\lambda + |\xi'|^2}, 
\quad B=\sqrt{\alpha^{-1}\lambda +|\xi'|^2}.$$
Here, $A$ and $B$ are characteristic roots of the Lam\'e equations 
given in \eqref{L.4.1} in Sect. \ref{sect.4.1} below.  The solution
formulas of equations \eqref{L.4.1} will be given in \eqref{L.0}. 
There exist two constants $d_1 < d_2$ depending on $0 < \epsilon <\pi/2$
such that 
\begin{equation}\label{top.1}
d_1(|\lambda|^{1/2} + |\xi'|) \leq {\rm Re}\, E \leq |E| \leq d_2(|\lambda|^{1/2}+|\xi'|)
\end{equation}
for any $(\lambda, \xi') \in \Sigma_\epsilon\times(\BR^{N-1}\setminus\{0\})$,
where $E \in \{A, B\}$. 
We can show the following two propositions using the same argument as in the proof of
Lemma 4.4 in Enomoto and Shibata \cite{ES13}.
\begin{prop} \label{prop:2}
Let $1 < q < \infty$, $0< \epsilon <\pi/2$, $\lambda_0>0$, and 
$\lambda\in\Sigma_{\epsilon, \lambda_0}$. Let  $m_0(\lambda, \xi') \in \bM_0$. Set 
\begin{equation}\label{Stokes-kernel}
M(x_N) = \frac{e^{-Bx_N}-e^{-Ax_N}}{B-A}.
\end{equation}
Define the integral operators $L_i$ ($i=1, 2, 3, 4$) by the following formulas:
\allowdisplaybreaks
\begin{align*}
L_1(\lambda)f &= \int^\infty_0\CF^{-1}_{\xi'}\left[m_0(\lambda, \xi')
Be^{-B(x_N+y_N)}\CF'[f](\xi', y_N)\right](x')\,dy_N, \\
L_2(\lambda)f &= \int^\infty_0\CF^{-1}_{\xi'}\left[m_0(\lambda, \xi')
B^2M(x_N + y_N)\CF'[f](\xi', y_N)\right](x')\,dy_N, \\
L_3(\lambda)f &= \int^\infty_0\CF^{-1}_{\xi'}\left[m_0(\lambda, \xi')
B^2\pd_\lambda(Be^{-B(x_N+y_N)})\CF'[f](\xi', y_N)\right](x')\,dy_N, \\
L_4(\lambda)f &= \int^\infty_0\CF^{-1}_{\xi'}\left[m_0(\lambda, \xi')
B^2\pd_\lambda(B^2M(x_N + y_N))\CF'[f](\xi', y_N)\right](x')\,dy_N,
\end{align*}
respectively. 
Then for every $f\in L_q(\HS)$, it holds
\begin{equation*}
\|L_i(\lambda)f\|_{L_q(\HS)} \leq C_q\|m_0\|\|f\|_{L_q(\HS)}\quad (i=1,2,3,4).
\end{equation*}
\end{prop}
\subsection{Estimates of operator valued holomorphic functions with respect to Besov norms}\label{subsec:2.5}
We consider two operator valued holomorphic 
functions $Q_i(\lambda)$ ($i=1,2$) defined on $\Sigma_\epsilon$ acting on
$f \in C^\infty_0(\HS)$. We denote the dual operator of $Q_i(\lambda)$ by $Q_i(\lambda)^*$
which satisfies the equality:
$$(Q_i(\lambda)f, \varphi)_{\HS} = (f, Q_i(\lambda)^*\varphi)_{\HS} \quad(i=1,2)
$$
for any $f$ and $\varphi \in C^\infty_0(\HS)$.  Here, $(f, g) = \int_{\HS} f(x)g(x)\,dx$.  
Let $Q_i(\lambda)$ satisfy the following assumptions.
\begin{assump}\label{assump.1} Let $1 < q < \infty$ and $q' = q/(q-1)$.  
 For any $f \in C^\infty_0(\HS)$ and 
$\lambda\in \Sigma_{\epsilon, \lambda_0}$, the following estimates hold:
\allowdisplaybreaks
\begin{align}
\|Q_1(\lambda)f\|_{W^i_q(\HS)}&\leq C\|f\|_{W^i_q(\HS)}, \label{as.1}\\
\|Q_1(\lambda)f\|_{L_q(\HS)}&\leq C|\lambda|^{-1/2}\|f\|_{W^1_q(\HS)}, \label{as.3}\\
\|Q_1(\lambda)^*f\|_{W^i_{q'}(\HS)}&\leq C\|f\|_{W^i_{q'}(\HS)}, \label{as.5} \\
\|Q_1(\lambda)^*f\|_{L_{q'}(\HS)}&\leq C|\lambda|^{-1/2}\|f\|_{W^1_{q'}(\HS)}, \label{as.6} \\
\|Q_2(\lambda)f\|_{W^i_q(\HS)}&\leq C|\lambda|^{-1}\|f\|_{W^i_q(\HS)}, \label{as.7} \\
\|Q_2(\lambda)f\|_{W^1_q(\HS)}&\leq C|\lambda|^{-1/2}\|f\|_{L_q(\HS)}, \label{as.9} \\
\|Q_2(\lambda)^*f\|_{W^i_{q'}(\HS)}&\leq C|\lambda|^{-1}
\|f\|_{W^i_{q'}(\HS)}, \label{as.11}\\
\|Q_2(\lambda)^*f\|_{W^1_{q'}(\HS)}&\leq C|\lambda|^{-1/2}\|f\|_{L_{q'}(\HS)} \label{as.12} 
\end{align}
for $i=0,1$, where $W^0_r(\HS) = L_r(\HS)$. 
\end{assump}
The following theorem has been proved in  \cite{S23} and \cite{SW}.
\begin{thm}\label{thm:5.2} Let $1 < q < \infty$ and $-1+1/q < s < 1/q$. 
Let $\sigma > 0$ be a number such that $-1+1/q < s-\sigma < s+\sigma < 1/q$ and 
let $\nu \in \{ s-\sigma, s, s+\sigma\}$. 
Let $Q_i(\lambda)$ $(i=1,2)$ be operator valued holomorphic 
functions defined on $\Sigma_{\epsilon, \lambda_0}$ acting on $C^\infty_0(\HS)$ functions.
 Then, for any $\lambda \in \Sigma_{\epsilon, \lambda_0}$ and $f \in C^\infty_0(\HS)$, 
the following assertions hold.\\
\thetag1 If $Q_1(\lambda)$ satisfies \eqref{as.1}, and \eqref{as.5}, 
then there holds
$$\|Q_1(\lambda)f \|_{B^\nu_{q,1}(\HS)}\leq C\|f\|_{B^\nu_{q,1}(\HS)}.$$
If $Q_1(\lambda)$ satisfies \eqref{as.3} and \eqref{as.6} in addition, then there holds
$$\|Q_1(\lambda)f\|_{B^s_{q,1}(\HS)}  
\leq C|\lambda|^{-\frac{\sigma}{2}}\|f\|_{B^{s+\sigma}_{q.r}(\HS)}.
$$
\thetag2 If $Q_2(\lambda)$ satisfies \eqref{as.7} and \eqref{as.11}, then there holds
$$\|Q_2(\lambda)f\|_{B^\nu_{q,1}(\HS)} \leq C|\lambda|^{-1}\|f\|_{B^\nu_{q,1}(\HS)}.$$
If $Q_2(\lambda)$ satisfies \eqref{as.9} and \eqref{as.12} in addition, then there holds
$$\|Q_2(\lambda)f\|_{B^s_{q,1}(\HS)}  \leq C|\lambda|^{-(1-\frac{\sigma}{2})}
\|f\|_{B^{s-\sigma}_{q,1}(\HS)}.$$
\end{thm}
\section{$L_1$ integrability of Laplace inverse transform.}

In this section, we consider the $L_1$ integrability of solutions to equations
\eqref{Eq:Linear}. The equations \eqref{Eq:Linear} is treated as a perturbation 
of Lam\'e equations with free boundary conditions.  In particular, 
one of main  issues is to prove the $L_1$ integrability of 
solutions in time, which is represented by the Laplace inverse transform
of solutions to the corresponding generalized resolvent problem. 
Thus, in this section, we introduce a solution operator $\CL_\Omega(\lambda)$
to the Lam\'e equations with free boundary conditions \eqref{Eq:Lame}.
We will construct a solution operator $\CL_\Omega(\lambda)$ which has
two properties.  One is stated in the following definition. 
\begin{dfn}
Let $\lambda_{lame} >0$ and $0 < \epsilon < \pi/2$.
Let $1 < q < \infty$ and $-1+1/q < s < 1/q$.  Let $\sigma > 0$ 
be a small number such that $-1+1/q < s-\sigma < s+\sigma < 1/q$.
Let $\nu \in \{s-\sigma, s, s+\sigma\}$.
Let $\CL_\Omega(\lambda) \in {\rm Hol}\,(\Sigma_{\epsilon, \lambda_{lame}}, 
\CL(B^s_{q,1}(\Omega)^{m(N)}, B^{s+2}_{q,1}(\Omega)^N))$.  We say that
$\CL_\Omega$ has  $(s, \sigma, q)$  properties in $\Omega$ if for any
$\lambda \in \Sigma_{\epsilon, \lambda_{lame}}$ there hold
\begin{equation}\label{L1prop}\begin{aligned}
\|(\lambda, \lambda^{1/2}\bar\nabla,
\bar\nabla^2)\pd_\lambda^\ell \CL_\Omega(\lambda)H\|_{B^\nu_{q,1}(\Omega)}
& \leq C|\lambda|^{-\ell}\|H\|_{B^\nu_{q,1}(\Omega)} \quad(\ell=0,1), \\
\|(\lambda^{1/2}\bar\nabla, \bar\nabla^2)\CL_\Omega(\lambda)H
\|_{B^s_{q,1}(\Omega)}& \leq C|\lambda|^{-\frac{\sigma}{2}}\|H\|_{B^{s+\sigma}_{q,1}
(\Omega)}, \\
\|(1, \lambda^{-1/2}\bar\nabla)\CL_\Omega(\lambda)H
\|_{B^s_{q,1}(\Omega)} & \leq C|\lambda|^{-(1-\frac{\sigma}{2})}\|H\|_{B^{s-\sigma}_{q,1}
(\Omega)}, \\
\|(\lambda, \lambda^{1/2}\bar\nabla, \bar\nabla^2)
\pd_\lambda\CL_\Omega(\lambda)H
\|_{B^s_{q,1}(\Omega)}& \leq C|\lambda|^{-(1-\frac{\sigma}{2})}
\|H\|_{B^{s-\sigma}_{q,1}(\Omega)}
\end{aligned}\end{equation}
provided that $H\in B^{s+\sigma}_{q,1}(\Omega)^{m(N)}$. 
\end{dfn}
\begin{remark} Since $s-\sigma < s < s+\sigma$, that 
$H\in B^{s+\sigma}_{q,1}(\Omega)$ implies that 
$H\in B^{\nu}_{q,1}(\Omega)$ for $\nu=s$ and $\nu=s-\sigma$.
\end{remark}
For equations \eqref{Eq:Lame}, we shall prove the existence of 
operators having the $L_1$ spectrum properties as follows.
\begin{thm}\label{thm:Lame}Let $0<\epsilon <\pi/2$. 
Assume that the following \thetag1 or \thetag2 holds:
\begin{itemize}
\item[\thetag1]~ If $\eta_0=\rho_*$, then
$1 < q < \infty$ and $-1+1/q < s < 1/q$. Let $\sigma>0$ be a number such that 
$-1+1/q < s-\sigma < s+\sigma < 1/q$. 
\item[\thetag2]~ If $\tilde\eta_0\not\equiv0$
and $\tilde\eta_0 \in B^{s+1}_{q,1}(\Omega)$.  then, $N-1 < q < 2N$,
$-1+N/q \leq s < 1/q$.  Let $\sigma>0$ be a small number such that
$s+\sigma < 1/q$ and $\sigma < 2N/q-1$. 
\end{itemize}
Then, there exists a positive 
constant $\lambda_{lame}$ and an operator $\CL_\Omega(\lambda)$ with
$$
\CL_\Omega(\lambda) \in {\rm Hol}\,(\Sigma_{\epsilon, \lambda_{lame}}, 
\CL(B^s_{q,1}(\Omega)^{m(N)}, B^{s+2}_{q,1}(\Omega)^N))
$$
having $(s, \sigma, q)$ properties such that 
$\bu = \CL_\Omega(\lambda)\CO_\lambda(\bg, \bh)$ is a unique solution
of equations \eqref{Eq:Lame} for $\lambda \in \Sigma_{\epsilon, \lambda_{lame}}$
and $(\bg, \bh) \in \CH^s_{q,1}(\Omega)$.  Here and in the following, 
$\CO_\lambda$ is an operation defined by 
$$\CO_\lambda(\bg, \bh) =(\bg, \lambda^{1/2}\bh, \nabla\bh).$$
\end{thm}
Let $L_\Omega(t)$ be the Laplace inverse transform 
of $\CL_\Omega$, which is defined by
\begin{equation*}
L_\Omega(t)H = \CL^{-1}[\CL_\Omega(\lambda)H](t)
= \frac{1}{2\pi i}\int_{\BR} e^{(\gamma +i\tau)t}\CL_\Omega(\gamma+i\tau)H\,d\tau.
\end{equation*}
Then, we have the following proposition about $L_\Omega(t)$.
\begin{prop}\label{prop:L1} Let $0< \epsilon < \pi/2$ and 
$\lambda_{lame} > 0$. Let $q$, $s$ and $\sigma$ be numbers given in
Theorem \ref{thm:Lame} and let $\CL_\Omega(\lambda)
\in {\rm Hol}\,(\Sigma_{\epsilon, \lambda_{lame}}, \CL(B^s_{q,1}(\Omega)^{m(N)},
B^{s+2}_{q,1}(\Omega)^N))$ be the operator 
having the $(s,\sigma,q)$ properties.  Then, $L_\Omega(t)H$ and 
$\Lambda_\gamma^{1/2}L_\Omega(t)H$ vanish for $t < 0$.  
Moreover, $e^{-\gamma t}L_\Omega(t)H \in L_1(\BR, B^{s+2}_{q,1}(\Omega)^N)$ and 
$e^{-\gamma t}\Lambda^{1/2}_\gamma L_\Omega(t)H \in L_1(\BR, B^{s+1}_{q,1}(\Omega)^N)$,  
and there holds  
\begin{equation}\label{Lame.L1}
\int_{\BR} e^{-\gamma t}\|L_\Omega(t)H\|_{B^{s+2}_{q,1}(\Omega)}\,dt
+ \int_{\BR} e^{-\gamma t}\|\Lambda^{1/2} L_\Omega(t)H\|_{B^{s+1}_{q,1}(\Omega)}\,dt
\leq C\|H\|_{B^s_{q,1}(\Omega)}
\end{equation}
for any $H \in B^s_{q,1}(\Omega)^{m(N)}$. \par
If for  any $G$ with $e^{-\gamma t}G \in L_1(\BR, B^s_q(\Omega))$, 
we define  $L_\Omega(t)G(t)$ by 
$$L_\Omega(t)G(t) = \CL^{-1}[\CL_\Omega(\lambda)\CL[G](\lambda)](t), $$
then,
there holds
\begin{equation}\label{Lame.L2}
\int_{\BR}e^{-\gamma t}(\|L_\Omega(t)G(t)\|_{B^{s+2}_{q,1}(\Omega)} + 
\|\Lambda^{1/2}L_\Omega(t)G(t)\|_{B^{s+1}_{q,1}(\Omega)})\,dt \leq 
C\int_{\BR}e^{-\gamma t}\|G(t)\|_{B^s_{q,1}(\Omega)}\, dt.
\end{equation}
\end{prop}
\begin{proof}
Since $C^\infty_0(\Omega)$ is dense in $B^{s+\sigma}_{q,1}(\Omega)$ and 
$B^s_{q,1}(\Omega)$, we may assume that $H \in C^\infty_0(\Omega)^{m(N)}$ 
below. 
First, we shall show that 
\begin{equation}\label{L1.1}
L_\Omega(t)H = 0 \quad\text{for $t < 0$}, \quad 
\Lambda^{1/2}_\gamma L_\Omega(t)H = 0\quad\text{for $t < 0$}.
\end{equation}
To prove \eqref{L1.1}, we  represent $L_\Omega(t)$ 
by using the contour integral in the complex plane
$\BC$. Let $C_R$ be a contour defined by 
\begin{align*}
C_R = \{\lambda \in \BC \mid \lambda =Re^{i\theta}, 
\, -\frac{\pi}{2} \leq \theta \leq \frac{\pi}{2}\}.
\end{align*} 
Let $\gamma > \lambda_{lame}$. 
By the Cauchy theorem in theory of one complex variable, we have
\begin{equation}\label{L1.2}
0 = \int_{-R}^{R} e^{(\gamma + i\tau)t}
\CL_\Omega(\gamma+i\tau )H\,d\tau
+ \int_{C_R+\gamma} e^{\lambda t}\CL_\Omega(\lambda)H\,d\lambda.
\end{equation}
Using \eqref{L1prop}, we know that 
$$\|\CL_\Omega(\lambda)H\|_{B^s_{q,1}(\Omega)}
\leq C|\lambda|^{-1}\|H\|_{B^s_{q,1}(\Omega)}.
$$
Thus, for $t < 0$ we have
\begin{align*}
\Bigl\|\int_{C_R+\gamma} e^{\lambda t}\CL_\Omega(\lambda)H
\,d\lambda\Bigr\|_{B^s_{q,1}(\Omega)} 
\leq Ce^{\gamma t}\int^{\pi/2}_0 e^{-|t|R\cos\theta}\,d\theta
\|H\|_{B^{s}_{q,1}(\Omega)}.
\end{align*}
Since $|e^{-|t|R\cos\theta}| \leq 1$, by Lebesgue's dominated convergence theorem, 
we have
$$\lim_{R\to\infty}\int^{\pi/2}_0 e^{-|t|R\cos\theta}\,d\theta
= \int^{\pi/2}_0\lim_{R\to\infty} e^{-|t|R\cos\theta}\,d\theta = 0.
$$
Therefore, letting $R\to\infty$ in \eqref{L1.2}, we have
$$0 = \frac{1}{2\pi i}\int_{\BR}e^{(\gamma + i\tau)t}\CL_\Omega(\gamma+i\tau )H
\,d\tau
= L_\Omega(t)H
$$
which proves the first part of \eqref{L1.1}. \par 
To prove the second part of  \eqref{L1.1}, we use  the similar argument.
Notice that
$$\Lambda_\gamma^{1/2} f = \frac{1}{2\pi}\int_{\BR} (\gamma+i\tau)^{1/2}e^{(\gamma + i\tau)t}
\CL[f](\gamma+i\tau)\, d\tau = \frac{1}{2\pi i}\int_{{\rm Re}\, \lambda =  \gamma}
\lambda^{1/2}e^{\lambda t}\CL[f](\lambda)\,d\lambda.$$
By the Cauchy theorem in theory of one complex variable, we have
\begin{equation}\label{L1.2*}
0 = \int_{-R}^{R} e^{(\gamma + i\tau)t} (\gamma+i\tau)^{1/2}
\CL_\Omega(\gamma+i\tau )H\,d\tau
+ \int_{C_R+\gamma} e^{\lambda t}\lambda^{1/2}\CL_\Omega(\lambda)H\,d\lambda.
\end{equation}
Using \eqref{L1prop}, we know that 
$$\|\lambda^{1/2}\CL_\Omega(\lambda)H\|_{B^s_{q,1}(\Omega)}
\leq C|\lambda|^{-1/2}\|H\|_{B^s_{q,1}(\Omega)}.
$$
Thus, for $t < 0$ we have
\begin{align*}
\Bigl\|\int_{C_R+\gamma} e^{\lambda t}\lambda^{1/2}\CL_\Omega(\lambda)H
\,d\lambda\Bigr\|_{B^s_{q,1}(\Omega)} 
\leq Ce^{\gamma t}\int^{\pi/2}_0 (R+\gamma)^{1/2}e^{-|t|R\cos\theta}\,d\theta
 = (*)
\end{align*}
Using the change of variable $\theta = \pi/2-\tau$ and the inequality: 
$\sin\tau \geq (2/\pi)\tau$ for $\tau \in (0, \pi/2)$,  we have 
\begin{align*}
(*) 
= C e^{\gamma t} \int^{\pi/2}_0 (R+\gamma)^{1/2}e^{-(2|t|R/\pi)\tau}\,d\tau 
\leq C e^{\gamma t}(R+\gamma)^{1/2}((|t|R/\pi)^{-1}
\to 0 \quad( R\to\infty).
\end{align*}
Therefore,  we have
$$0 = \frac{1}{2\pi i}\int_{\BR}(\gamma+i\tau)^{1/2}
e^{(\gamma + i\tau)t}\CL_\Omega(\gamma+i\tau )H \,d\tau
= \Lambda^{1/2}_\gamma L_\Omega(t)H
$$
which proves the second part of \eqref{L1.1}.

We next consider the case where $t>0$.  Let $\Gamma_\pm$ be contours defined by
$$\Gamma_\pm = \{\lambda = re^{\pm i(\pi-\epsilon)} \mid r \in (0, \infty)\}.$$
We shall show that for $t>0$ $L_\Omega(t)$ is represented by 
\begin{equation}\label{repr.L1}
L_\Omega(t)H = \frac{1}{2\pi i}\int_{\Gamma_+ \cup \Gamma_- + \gamma}
e^{\lambda t}\CL_\Omega(\lambda)H\, d\lambda.
\end{equation}
In fact, for $R > 0$, define $\tilde C_{R\pm}$ and $\Gamma_{R \pm}$ by
$$\tilde C_{R\pm} = \{\lambda = Re^{i\theta} \mid \pi/2 < \pm \theta < \pi-\epsilon\},
\quad
\Gamma_{R\pm} = \{\lambda \in \Gamma_\pm \mid |\lambda| < R\}.$$
By the Cauchy theorem in theory of one complex variable, we have
\begin{equation}\label{repr.L2}
\frac{1}{2\pi} \int^R_{-R} e^{i\tau t}\CL_\Omega(\gamma + i\tau)H\, d\tau
+ \frac{1}{2\pi i}\Bigl\{\int_{\tilde C_{R+}+\gamma} -\int_{\Gamma _{R+}+\gamma} 
- \int_{\Gamma_{R-}+\gamma }
+ \int_{\tilde C_{R-}+\gamma }\Bigr\}e^{\lambda t}\CL_\Omega(\lambda)H\,d\lambda
=0.
\end{equation}
Using \eqref{L1prop} and the change of variable: $\theta=\tau+\pi/2$, 
for $R > \gamma$ we have
\begin{align*}
\Bigl\|\int_{\tilde C_{R+}+\gamma} e^{\lambda t}\CL_\Omega(\lambda)H
\,\d\lambda\Bigr\|_{B^s_{q,1}(\Omega)}
&\leq \frac{e^{\gamma t}}{2\pi}
\int^{\pi-\epsilon}_{\pi/2}e^{R\cos\theta t}\frac{R}{|Re^{i\theta}+\gamma|}
\,d\theta \|H\|_{B^s_{q,1}(\Omega)}\\
&\quad \leq e^{\gamma t}\frac{R}{R-\gamma} 
\int^{\pi/2-\epsilon}_0 e^{-R\sin\theta t}\,d\theta
\|H\|_{B^s_{q,1}(\Omega)}\\
&\leq \frac{e^{\gamma t}}{2\pi}\frac{R}{R-\gamma} 
\int^{\pi/2}_0 e^{-\frac{2Rt}{\pi}\tau}\,d\tau
\|H\|_{B^s_{q,1}(\Omega)} \\
&\leq \frac{e^{\gamma t}}{2\pi}
\frac{R}{R-\gamma} \frac{\pi}{2Rt}\|H\|_{B^s_{q,1}(\Omega)}.
\end{align*}
Thus, for $t > 0$ 
$$\lim_{R\to\infty} \int_{\tilde C_{R+}+\gamma} e^{\lambda t}\CL_\Omega(\lambda)H = 0.$$
Analogously, we have
$$\lim_{R\to\infty} \int_{\tilde C_{R-}+\gamma} e^{\lambda t}\CL_\Omega(\lambda)H = 0.$$
Combining these facts and letting $R\to\infty$ in \eqref{repr.L2} yield \eqref{repr.L1}.
\par
Using the representation formula \eqref{repr.L1}, we shall show that 
\begin{align}\label{L1.3}
\|\bar\nabla^2L_\Omega(t)H\|_{B^s_{q,1}(\Omega)} &\leq Ce^{\gamma t}
t^{-(1-\frac{\sigma}{2})}
\|H\|_{B^{s+\sigma}_{q,1}(\Omega)}, \\
\label{L1.4}
\|\bar\nabla^2L_\Omega(t)H\|_{B^s_{q,1}(\Omega)} &\leq Ce^{\gamma t}
t^{-(1+\frac{\sigma}{2})}
\|H\|_{B^{s-\sigma}_{q,1}(\Omega)}.
\end{align}
In fact, noticing that $|e^{\lambda t}| 
= e^{-tr\cos\epsilon}$ and $|\gamma + re^{\pm i (\pi-\epsilon)}|
\geq (1-\cos\epsilon)^{1/2}r$ 
for $\lambda \in \Gamma_+\cup \Gamma_-+\gamma$ in \eqref{repr.L1}
and using \eqref{L1prop}, we have 
\begin{align*}
&\|\bar\nabla^2L_\Omega(t)H\|_{B^s_{q,1}(\Omega)} \\
&\quad \leq Ce^{\gamma t}\int^\infty_0 e^{-tr\cos\epsilon}
((1-\cos\epsilon)^{1/2}r)^{-\sigma/2}\, dr \\
&\quad 
\leq Ce^{\gamma t}t^{-1+\frac{\sigma}{2}}
\int^\infty_0 e^{-\tau\cos\epsilon}
((1-\cos\epsilon)^{1/2}\tau)^{-\frac{\sigma}{2}}\, d\tau
\|H\|_{B^{s+\sigma}_{q,1}(\Omega)}.
\end{align*}
From this \eqref{L1.3} follows. \par 
Using integration by parts concerning $\lambda$ in \eqref{repr.L1}, we have 
$$\bar\nabla^2L_\Omega(t)H = \frac{-1}{2\pi it}
\int_{\Gamma_+\cup \Gamma_-+\gamma}e^{\lambda t} 
\bar\nabla^2 \pd_\lambda \CL_\Omega(\lambda)H \,d\lambda.
$$
Using \eqref{L1prop}, 
we have 
\begin{align*}
\|\bar\nabla^2L_\Omega(t)H\|_{B^s_{q,1}(\Omega)} 
&\leq \Bigl\|\frac{1}{2\pi i t}\int_{\Gamma_+\cup\Gamma_-+\gamma} e^{\lambda t} 
\bar\nabla^2 \pd_\lambda \CL_\Omega(\lambda)H\,d\lambda\Bigr\|_{B^s_{q,1}} \\
&\leq Ce^{\gamma t} t^{-1}
\int^\infty_0 e^{-tr\cos\epsilon}
((1-\cos\epsilon)^{1/2}r)^{-(1-\frac{\sigma}{2})}\, dr
\|H\|_{B^{s-\sigma}_{q,1}(\Omega)} \\
&=Ce^{\gamma t}
t^{-1-\frac{\sigma}{2}}
\int^\infty_0 e^{-\tau\cos\epsilon}
((1-\cos\epsilon)^{1/2}\tau)^{-(1-\frac{\sigma}{2})}\, d\tau
\|H\|_{B^{s-\sigma}_{q,1}(\Omega)}.
\end{align*}
Therefore, we have \eqref{L1.4}. \par
Now, we shall estimate $\Lambda_\gamma^{1/2}L\Omega(t)H$. 
We shall show that  
\begin{align}
\label{L1.3*}
\|\bar\nabla \Lambda^{1/2} L_\Omega(t)H
\|_{B^s_{q,1}(\Omega)} &\leq Ce^{\gamma t}
t^{-(1-\frac{\sigma}{2})}
\|H\|_{B^{s+\sigma}_{q,1}(\Omega)}, \\
\label{L1.4*}
\|\bar\nabla \Lambda^{1/2}
L_\Omega(t)H\|_{B^s_{q,1}(\Omega)} &\leq Ce^{\gamma t}
t^{-(1+\frac{\sigma}{2})}
\|H\|_{B^{s-\sigma}_{q,1}(\Omega)}.
\end{align}
In fact, for $\gamma > \gamma_b$ and $t>0$ from \eqref{repr.L1} we have 
\begin{equation}\label{repr.L1*}
\Lambda_\gamma^{1/2}L(t)H 
= \frac{1}{2\pi i}\int_{\Gamma_+ \cup \Gamma_- + \gamma}
\lambda^{1/2} e^{\lambda t}\CL(\lambda) H \, d\lambda.
\end{equation}
Using the assumption \eqref{L1prop} and the same argument as in
the proof of \eqref{L1.3}, we have \eqref{L1.3*}. \par 
Using integration by parts in \eqref{repr.L1*},  we have
\begin{align*}
\bar\nabla\Lambda^{1/2}L_\Omega(t)H &= \frac{-1}{2\pi it}
\int_{\Gamma_+\cup \Gamma_-+\gamma}e^{\lambda t} 
\bar\nabla \pd_\lambda(\lambda^{1/2}  \CL_\Omega(\lambda)H) \,d\lambda\\
&= \frac{-1}{2\pi it}
\int_{\Gamma_+\cup \Gamma_-+\gamma}e^{\lambda t} 
 ((1/2)\bar\nabla\lambda^{-1/2}\CL_\Omega(\lambda)H
+\bar\nabla^{1/2}\lambda^{1/2}\pd_\lambda(\CL_\Omega(\lambda)H) ) 
\,d\lambda.
\end{align*}
Using \eqref{L1prop} and the same argument as in the proof 
 of \eqref{L1.4}, we have \eqref{L1.4*}. 
\par

Now, we shall prove \eqref{Lame.L1} 
by using \eqref{L1.3},  \eqref{L1.4}, \eqref{L1.3*}, and  \eqref{L1.4*}.
We write
\begin{align*}
\int^\infty_0 e^{-\gamma t}\|(\bar\nabla^2, \bar\nabla\Lambda^{1/2}) L_\Omega(t)H
\|_{B^{s}_{q,1}(\Omega)}\, dt 
&= \sum_{j \in \BZ} \int^{2^{(j+1)}}_{2^j}e^{-\gamma t}
\|(\bar\nabla^2, \bar\nabla\Lambda^{1/2}) L_\Omega(t)H\|_{B^{s}_{q,1}(\Omega)}\, dt\\
&\leq
\sum_{j \in \BZ} \, 2^j
\sup_{t \in (2^j, 2^{j+1})}(e^{-\gamma t}\|(\bar\nabla^2, \bar\nabla\Lambda^{1/2}) 
L_\Omega(t)H \|_{B^{s}_{q,1}(\Omega)}).
\end{align*}
Setting $a_j = \sup_{t \in (2^j, 2^{j+1})}e^{-\gamma t}
\|(\bar\nabla^2, \bar\nabla\Lambda^{1/2}) L_\Omega(t)H\|_{B^{s}_{q,1}(\Omega)}$, 
we have
$$
\int^\infty_0 e^{-\gamma t}\|(\bar\nabla^2, \bar\nabla\Lambda^{1/2}) 
L_\Omega(t)H\|_{B^{s}_{q,1}(\Omega)}\, dt
\leq 2 \|(2^ja_j)\|_{\ell_1} = 2\|(a_j)\|_{\ell^1_1}. 
$$
Here and in the following, 
$\ell^s_q$ denotes the set of all sequences $(2^{js}a_j)_{j \in \BZ}$ such that 
\begin{alignat*}2
\|(a_j)\|_{\ell^s_q} &= \Bigl\{\sum_{j \in \BZ} 
(2^{js}|a_j|)^q \Bigr\}^{1/q} < \infty &\quad &\text{for $1 \leq q < \infty$},  \\
\|(a_j)\|_{\ell^s_\infty} &= \sup_{j \in \BZ} 
2^{js}|a_j| < \infty &\quad &\text{for  $q = \infty$}.
\end{alignat*}
By \eqref{L1.3} and \eqref{L1.4},  we have
$$\sup_{j \in \BZ} 2^{j(1-\frac{\sigma}{2})}a_j \leq C\|H\|_{B^{s+\sigma}_{q,1}(\Omega)}, 
\quad\sup_{j \in \BZ} 2^{j(1+\frac{\sigma}{2})}a_j \leq C\|H\|_{B^{s-\sigma}_{q,1}(\Omega)}
\quad (H \in B^s_{q,1}(\Omega)).
$$
Namely, we have
$$\|(a_j)_{j \in \BZ}\|_{\ell_\infty^{1-\frac{\sigma}{2}}} 
\leq C\|H\|_{B^{s+\sigma}_{q,1}(\Omega)}, \quad 
\|(a_j)_{j \in \BZ}\|_{\ell_\infty^{1+ \frac{\sigma}{2}}} 
\leq C\|H\|_{B^{s-\sigma}_{q,1}(\Omega)} \quad (H \in B^s_{q,1}(\Omega)).
$$
According to \cite[5.6.1.Theorem]{BL},  
we know that $\ell^1_1 = (\ell^{1-\frac{\sigma}{2}}_\infty,
\ell^{1+\frac{\sigma}{2}}_\infty)_{1/2, 1}$, 
where
$(\cdot, \cdot)_{\theta, q}$ denotes the real interpolation functor, 
and therefore we have 
$$
\int^\infty_0 e^{-\gamma t}\|(\bar\nabla^2, \bar\nabla\Lambda^{1/2}) L_\Omega(t)H
\|_{B^{s}_{q,1}(\Omega)}\, dt
\leq C
\|H\|_{(B^{s+\sigma}_{q,1}(\Omega), B^{s-\sigma}_{q,1}(\Omega))_{1/2, 1}} 
=C
\|H\|_{B^{s}_{q,1}(\Omega)}
$$
for any $H \in B^s_{q,1}(\Omega)^{m(N)}$. 
This  completes the proof of \eqref{Lame.L1}.\par
To prove \eqref{Lame.L2}, we write 
\begin{align*}
(1, \Lambda^{1/2}) 
L(t)G(t) &= \frac{1}{2\pi i}\int_{\BR}e^{\lambda t}(1, \lambda^{1/2})
\CL_\Omega(\lambda)\CL[G](\lambda)\, d\lambda \\
&= \frac{1}{2\pi }\int_{\BR}e^{\lambda t}(1, \lambda^{1/2})\CL_\Omega(\lambda)
\CL[G](\lambda)\Bigl(\int_{\BR} e^{-\lambda s}G(s)\, ds\Bigr)\,d\tau\\ 
& = \int_{\BR}\Bigl(\frac{1}{2\pi} \int_{\BR} e^{\lambda(t-s)}
(1, \lambda^{1./2})\CL_\Omega(\lambda)G(s)\,
d\lambda\Bigr)\, ds \\
&= \int_{\BR}(1, \Lambda^{1/2}) L_\Omega(t-s)G(s)\,ds \\
& = \int^t_{-\infty}(1, \Lambda^{1/2}) L_\Omega(t-s)G(s)\,ds.
\end{align*}
Here, $\lambda = \gamma + i\tau$  and we have used the fact that
$L_\Omega(t)=0$ and $\Lambda^{1/2}L_\Omega(t) = 0$ for $t < 0$.
By Fubini's theorem, we have \eqref{Lame.L2}.  
This completes the proof of Proposition \ref{prop:L1}.
\end{proof}
To treat the perturbation term, we introduce one more definition.
\begin{dfn}\label{dfn.2} 
Let $\lambda_0>0$ and $0 < \epsilon < \pi/2$. 
Let $X$ be a Banach space and $\CM(\lambda) \in {\rm Hol}\,(\Sigma_{\epsilon,
\lambda_0}, \CL(X))$.  We say that $\CM(\lambda)$ has a generalized
resolvent properties if there hold
\begin{equation*}
\|\pd_\lambda^\ell \CM(\lambda)f\|_X \leq C|\lambda|^{-\ell-1}
\|f\|_X
\quad\text{for $f \in X$ and $\ell=0,1$}.
\end{equation*}
\end{dfn}
Let $M(t)$ be the Laplace inverse transform of $\CM(\lambda)$ defined by
$$M(t)f = \CL^{-1}[\CM(\lambda)f] = \frac{1}{2\pi i} \int_{\BR}e^{(\gamma  + i\tau)t}
\CM(\gamma + i\tau)f\,d\tau.
$$
Then, we have the following proposition about the $L_1$ integrability of $M(t)$.
\begin{prop} \label{prop:L2} 
Let $\lambda_0>0$ and $0 < \epsilon < \pi/2$. 
Let $X$ be a Banach space and $\CM(\lambda) \in {\rm Hol}\,(\Sigma_{\epsilon,
\lambda_0}, \CL(X))$. If $\CM(\lambda)$ has  generalized resolvent properties.
Then, for $f \in X$, $M(t)=0$ for $t < 0$,  and for any $\gamma > \lambda_0$ it holds that 
$$\int_{\BR}e^{-\gamma t}\|M(t)f\|_X \, dt \leq C\|f\|_X$$
with some constant $C$ depending on $\lambda_0$. \par 
Moreover, if  we define  $M(t)g(t)$ by
$$M(t)g(t) = \CL^{-1}[\CM(\lambda)\CL[g](\lambda)]$$
for $g(t)$ with $e^{-\gamma t}g(t) \in L_1(\BR, X)$, then there holds
$$\int_{\BR}e^{-\gamma t}\|M(t)g(t)\|_X \, dt \leq C\int_{\BR} e^{-\gamma t}\|g(t)\|_X\,dt.
$$
with some constant $C$ depending on $\lambda_0$. 
\end{prop}
\begin{proof} For $\lambda \in \Sigma_{\epsilon, \lambda_0}$, we have
\begin{align*}
\|\CM(\lambda)f\|_X \leq C|\lambda|^{-1}\|f\|_X
&\leq C\lambda_0^{-(1-\frac{\sigma}{2})}|\lambda|^{-\frac{\sigma}{2}}\|f\|_X,
\\
\|\pd_\lambda \CM(\lambda)f\|_X \leq C|\lambda|^{-2}\|f\|_X
&\leq C\lambda_0^{-(1+\frac{\sigma}{2})}|\lambda|^{-(1-\frac{\sigma}{2})}\|f\|_X
\end{align*}
for any $\lambda \in \Sigma_{\epsilon, \lambda_0}$.  Thus, employing the 
same argument as in the proof of Proposition \ref{prop:L1}, we can prove
Proposition \ref{prop:L2}.  This completes the proof. 
\end{proof}
It is stated in the last section \ref{sec:6} to prove Theorem \ref{thm:main.1} 
using Propositions \ref{prop:L1} and \ref{prop:L2}. 

\section{The spectral analysis of the Lam\'e equations}
In view of Proposition \ref{prop:L1}, to prove the $L_1$ integrability of 
solutions to the evolution equations \eqref{Eq:Linear}, we start with proving
Theorem \ref{thm:Lame}.  Our proof is divided into several steps.
\subsection{Estimate of solutions in $\BR^N$} \label{sect.4.1}
In this subsection, we consider the Lam\'e equations in $\BR^N$: 
\begin{equation}\label{L.4.1}
\lambda \bu - \DV(\alpha\BD(\bu) + (\beta -\alpha)
\dv\bu \BI ) = \bg \quad\text{in $\BR^N$}
\end{equation}
for $\lambda \in \Sigma_\epsilon$ with $0 < \epsilon < \pi/2$. 
Notice that $\DV(\alpha\BD(\bu) + (\beta -\alpha)
\dv\bu \BI ) = \alpha\Delta\bu+ \beta\nabla\dv\bu$. 
We shall prove the following theorem. 
\begin{thm}\label{thm:L.1} Let $1 < q < \infty$,  $-1+1/q < s < 1/q$ and 
$\lambda_0>0$.
Let  $\sigma>0$ be a small number such that 
$-1+1/q < s-\sigma < s+\sigma<1/q$. 
Then, there exists an operator $\CS(\lambda) \in 
{\rm Hol}\,( \Sigma_\epsilon, \CL(B^s_{q,1}(\BR^N)^N, B^{s+2}_{q,1}(\BR^N)^N))$ 
having $(s, \sigma, q)$ properties such that 
for any $\bg \in B^s_{q,1}(\BR^N)^N$, 
$\bu = \CS(\lambda)\bg$ is a unique solution of equations \eqref{L.4.1}. 
\end{thm}
\begin{proof} Let $\nu \in \{s-\sigma, s, s+\sigma\}$. 
Applying the divergence to equations \eqref{L.4.1} gives 
$$\lambda \dv \bu - (\alpha + \beta)\Delta \dv\bu = \dv\bg \quad\text{in $\BR^N$}.
$$
Using the Fourier transform $\CF$ and its inverse transform $\CF^{-1}$, 
we have
$$\dv \,\bu = \CF^{-1}\Bigl[\frac{i\xi\cdot
\CF[\bg](\xi)}{\lambda+(\alpha+\beta)|\xi|^2}\Bigr].$$
Inserting this formula into \eqref{L.4.1} gives that 
\begin{equation}\label{L.0}
\begin{aligned}
\bu & = \CF^{-1}\Bigl[\frac{\CF[\bg](\xi)}{\lambda + \alpha|\xi|^2}\Bigr]
+ \beta\CF^{-1}\Bigl[\frac{i\xi i\xi\cdot\CF[\bg](\xi)}{(\lambda + \alpha|\xi|^2)
(\lambda + (\alpha+\beta)|\xi|^2)}\Bigr].
\end{aligned}\end{equation}
As we know well, there exist positive constants $c_1$ and $c_2$ depending on $\alpha$,
$\beta$ and $\epsilon$ such that for any $\lambda \in \Sigma_\epsilon$ there hold:
\begin{equation}\label{L.2}\begin{aligned}
c_1(|\lambda|^{1/2}+|\xi|)  \leq &|\lambda + \alpha|\xi|^2| \leq c_2(|\lambda|^{1/2}+|\xi|), \\
c_1(|\lambda|^{1/2}+|\xi|)  \leq |\lambda &+ (\alpha+\beta)|\xi|^2| \leq c_2(|\lambda|^{1/2}+|\xi|).
\end{aligned}\end{equation}
Moreover, 
\begin{equation*}
\pd_\lambda \bu 
 = \CF^{-1}\Bigl[\pd_\lambda\frac{\CF[\bg](\xi)}{\lambda + \alpha|\xi|^2}\Bigr]
+ \beta\CF^{-1}\Bigl[\pd_\lambda \frac{i\xi i\xi\cdot\CF[\bg](\xi)}{(\lambda + \alpha|\xi|^2)
(\lambda + (\alpha+\beta)|\xi|^2)}\Bigr]
\end{equation*}
and we see that 
\begin{equation*}
\begin{aligned}
c_1(|\lambda|^{1/2}+|\xi|)^{-4}  \leq |\pd_\lambda(\lambda + \alpha|\xi|^2)^{-1}| 
\leq c_2(|\lambda|^{1/2}+|\xi|)^{-4}, \\
c_1(|\lambda|^{1/2}+|\xi|)^{-6}  \leq |\pd_\lambda((\lambda+\alpha|\xi|^2)
(\lambda + (\alpha+\beta)|\xi|^2))^{-1}| \leq c_2(|\lambda|^{1/2}+|\xi|)^{-6}.
\end{aligned}\end{equation*}
Here, if necessary, we choose different $c_1$ and $c_2$ from \eqref{L.2}. 
Thus, applying the Fourier multiplier theorem of 
Mikhlin-H\"ormander type, we have 
$$
\|(\lambda, \lambda^{1/2}\bar\nabla, \bar\nabla^2)\pd_\lambda^\ell
\bu\|_{B^\nu_{q,1}(\BR^N)} 
\leq C|\lambda|^{-\ell}\|\bg\|_{B^\nu_{q,1}(\BR^N)} \qquad (\ell=0,1).
$$
For $\bg \in B^{s+\sigma}_{q,1}(\BR^N)$, we write
\begin{align*}
&\lambda^{1/2}\lambda^{\sigma/2}\bar\nabla \bu \\
&\quad
= \CF^{-1}\Bigl[\frac{\lambda^{\frac12+\frac{\sigma}{2}}(1, i\xi)(1+|\xi|^2)^{\sigma/2}\CF[\bg](\xi)}
{(\lambda + \alpha|\xi|^2)(1+|\xi|^2)^{\sigma/2}}\Bigr] 
+ \beta\CF^{-1}\Bigl[
\frac{\lambda^{\frac12+\frac{\sigma}{2}}(1, i\xi)i\xi i\xi\cdot
((1+|\xi|^2)^{\sigma/2}\CF[\bg](\xi))}{(\lambda + \alpha|\xi|^2)
(\lambda + (\alpha+\beta)|\xi|^2)(1+|\xi|^2)^{\sigma/2}}\Bigr], \\
&\lambda^{1/2}\lambda^{1-\sigma/2}\bar\nabla \pd_\lambda \bu
= \CF^{-1}\Bigl[\lambda^{\frac32-\frac{\sigma}{2}}(1, i\xi)
\{(1+|\xi|^2)^{-\sigma/2}\CF[\bg](\xi)\}
\Bigl(\pd_\lambda\frac{1}
{(\lambda + \alpha|\xi|^2)}\Bigr)(1+|\xi|^2)^{\sigma/2}\Bigr] \\
&\quad
+ \beta\CF^{-1}\Bigl[\lambda^{\frac32-\frac{\sigma}{2}}(1, i\xi)i\xi i\xi\cdot
\{(1+|\xi|^2)^{-\sigma/2}\CF[\bg](\xi)\}
\Bigl(\pd_\lambda\frac{1}{(\lambda + \alpha|\xi|^2))
(\lambda + (\alpha+\beta)|\xi|^2)}\Bigr)(1+|\xi|^2)^{\sigma/2}\Bigr]. 
\end{align*}
Applying the Fourier multiplier theorem of Mikhilin-H\"ormander type, we have
\begin{align*}
\|\lambda^{1/2}\bar\nabla \bu\|_{B^s_{q,1}(\BR^N)}
&\leq C(1+\lambda_0^{-1/2})|\lambda|^{-\frac{\sigma}{2}}
\|\bg\|_{B^{s+\sigma}_{q,1}(\BR^N)}, \\
\|\lambda^{1/2}\bar\nabla \pd_\lambda\bu\|_{B^s_{q,1}(\BR^N)}
&\leq C(1+\lambda_0^{-1/2})|\lambda|^{-(1-\frac{\sigma}{2})}
\|\bg\|_{B^{s-\sigma}_{q,1}(\BR^N)}.
\end{align*}
Analogously, we have
\begin{align*}
\|\bar\nabla^2 \bu\|_{B^s_{q,1}(\BR^N)}
&\leq C(1+\lambda_0^{-1/2}+\lambda_0^{-1})|\lambda|^{-\frac{\sigma}{2}}
\|\bg\|_{B^{s+\sigma}_{q,1}(\BR^N)}, \\
\|\bar\nabla^2 \pd_\lambda\bu\|_{B^s_{q,1}(\BR^N)}
&\leq C(1+\lambda_0^{-1/2}+\lambda_0^{-1})|\lambda|^{-(1-\frac{\sigma}{2})}
\|\bg\|_{B^{s-\sigma}_{q,1}(\BR^N)},
\\
\|(1, \lambda^{-1/2}\bar\nabla)\bu\|_{B^s_{q,1}(\BR^N)}
&\leq C(1+\lambda_0^{-1/2})
|\lambda|^{-(1-\frac{\sigma}{2})}\|\bg\|_{B^{s-\sigma}_{q,1}(\BR^N)}.
\end{align*}
Define $\CS(\lambda)$ by $\CS(\lambda)\bg = \bu$, and then  we see that 
$\CS(\lambda)$ has $(s, \sigma, q)$ properties and $\bu = \CS(\lambda)\bg$
is a solution of equations \eqref{L.4.1}. The uniqueness follows from the existence 
of solutions to the dual problem. This completes the proof of Theorem \ref{thm:L.1}.
\end{proof}
\par
\subsection{Solution formulas of Lam\'e equations in the half-space}
In this subsection, we find solution formulas of  the Lam\'e equations in the half-space 
with  free boundary conditions which reads as 
\begin{equation}\label{free.1}\begin{aligned}
\lambda\bu - \DV(\alpha\BD(\bu) + (\beta-\alpha)\dv\bu\BI) = 0&&\quad&
\text{in $\HS$}, \\
(\alpha\BD(\bu) + (\beta-\alpha)\dv\bu\BI)\bn_0 = \bh&&\quad&
\text{on $\pd\BR^N_+$}.
\end{aligned}\end{equation} 
Since $\DV \BD(\bu) = \Delta \bu + \nabla\dv\bu$, \eqref{free.1} is rewritten 
as follows: 
\begin{equation*}
\begin{aligned}
&\lambda\bu - \alpha\Delta \bu - \beta\nabla\dv\bu= 0 &\quad &\text{in $\HS$}, 
\\
&-\alpha(\pd_Nu_j + \pd_ju_N)|_{x_N=0} = h_j|_{x_N=0}, & & \\
&-(2\alpha \pd_Nu_N + (\beta-\alpha)\dv\bu)|_{x_N=0}=h_N|_{x_N=0} & &
\end{aligned}\end{equation*}
for $j=1, \ldots, N-1$. 
In this subsection, indices $j$ and $J$ run from $1$ through $N-1$ and 
from 1 through $N$, respectively.  
We apply the partial Fourier transform $\CF'$ with respect to 
$x' = (x_1, \ldots, x_{N-1})$ variables and  then
equations \eqref{free.1} are transformed to the following system of 
ordinary differential equations:
\begin{equation}\label{fp1}\begin{aligned}
&(\lambda +\alpha |i\xi'|^2)\CF'[u_j]-\alpha \pd_N^2\CF'[u_j] 
-\beta i\xi_j(i\xi'\cdot \CF'[\bu'] + \pd_N\CF'[u_N]) =0 \qquad(x_N>0), \\
&(\lambda +\alpha |i\xi'|^2)\CF'[u_N]-\alpha \pd_N^2\CF'[u_N] 
-\beta \pd_N(i\xi'\cdot \CF'[\bu'] + \pd_N\CF'[u_N]) =0 \quad(x_N>0), \\
&(\pd_N\CF'[u_j] + i\xi_j \CF'[u_N])|_{x_N=0} = -\alpha^{-1}\CF'[h_j](\xi', 0), \\
&(\alpha (\pd_N\CF'[u_N] - i\xi'\cdot\CF'[\bu']) 
+ \beta(i\xi'\cdot\CF'[\bu'] +\pd_N\CF'[u_N]))|_{x_N=0}= -\CF'[h_N](\xi', 0).
\end{aligned}\end{equation}
To obtain solutions formula, we define $\CF[u_J]$  by 
\begin{equation*}
\CF[u_J](\xi', x_N) = m_J e^{-Bx_N} + n_J (e^{-Ax_N} - e^{-Bx_N}).
\end{equation*}
From the first two equations in \eqref{fp1}, we have 
\begin{equation}\label{fp2}\begin{aligned}
\alpha(B^2-A^2)n_j - \beta i\xi_j(i\xi'\cdot\bn'-An_N) = 0,&&\quad&
\beta i\xi_j(i\xi'\cdot(\bm'-\bn') -B(m_N-n_N))=0,\\
\alpha(B^2-A^2)n_N+\beta A(i\xi'\cdot\bn'-An_N) = 0,&&\quad&
\beta B(i\xi'\cdot(\bm'-\bn') - B(m_N-n_N))=0. 
\end{aligned}\end{equation}
From the boundary conditions it follows that
\begin{align}\label{fp3}
& B(m_j-n_j)+An_j - i\xi_j m_N = \alpha^{-1}\hat h_j, \\
&\alpha(B(m_N-n_N)+An_N +i\xi'\cdot\bm') + \beta(B(m_N-n_N) +An_N-i\xi'\cdot\bm') 
= \hat h_N.
\label{fp3*}
\end{align}
Setting  $\ell =\alpha^{-1}\beta (i\xi'\cdot\bn'-An_N)$,  from \eqref{fp2} we have
\begin{align}
&n_j = \frac{ i\xi_j}{B^2-A^2}\ell, \quad n_N = -\frac{A}{B^2-A^2}\ell, \label{fp4}\\
&i\xi'\cdot\bm'-i\xi'\cdot\bn' -B(m_N-  n_N) = 0. \label{fp5}
\end{align}
Since $-B(m_N-n_N) +i\xi'\cdot\bm' = i\xi'\cdot\bn'$ as follows from \eqref{fp5}, 
the last boundary condition in \eqref{fp3*} becomes
\begin{equation}\label{fp6}
Bm_N + (A-B)n_N + i\xi'\cdot\bm' -\ell = \alpha^{-1}\hat h_N.
\end{equation}
From \eqref{fp4} we obtain 
\begin{equation}\label{f1.26}
i\xi'\cdot\bn' = \frac{-|\xi'|^2}{B^2-A^2}\ell.
\end{equation}
Substituting \eqref{fp4} and \eqref{f1.26} into \eqref{fp5}, we get 
\begin{align}
m_N & = B^{-1}\Bigl(i\xi'\cdot\bm' +\frac{|\xi'|^2-AB}{B^2-A^2}\ell\Bigr).
\label{f1.27}
\end{align}
From \eqref{fp3}, \eqref{f1.26} and \eqref{f1.27}, we have  
\begin{equation}\label{f1.28}
\alpha^{-1}Bi\xi'\cdot\hat \bh' 
= (B^2+|\xi'|^2)i\xi'\cdot\bm' + |\xi'|^2\frac{B^2+|\xi'|^2-2AB}{B^2-A^2}\ell.
\end{equation}
Substituting \eqref{fp4} and \eqref{f1.27} into \eqref{fp6}, we obtain 
\begin{align}
\alpha^{-1}\hat h_N 
& = 2i\xi'\cdot\bm' + \frac{|\xi'|^2-B^2}{B^2-A^2}\ell. \label{f1.29}
\end{align}
Let $L$ be a $2\times 2$ matrix defined by
$$L = \left(\begin{matrix} B^2+|\xi'|^2 & |\xi'|^2\dfrac{B^2+|\xi'|^2-2AB}{B^2-A^2} \vspace{1pc} \\
2 & \dfrac{|\xi'|^2-B^2}{B^2-A^2} \end{matrix}\right).
$$
From \eqref{f1.28} and \eqref{f1.29}, it follows that 
\begin{equation*}
L\left(\begin{matrix} i\xi'\cdot\bm' \\ \ell \end{matrix}\right) = \left(\begin{matrix}
\alpha^{-1}Bi\xi'\cdot\hat \bh' \\ \alpha^{-1}\hat h_N\end{matrix}\right).
\end{equation*}
We observe that 
\begin{align*}
\det L & = \frac{4AB|\xi'|^2-(B^2+|\xi'|^2)^2}{B^2-A^2}.
\end{align*}
Thus, setting 
\begin{equation*}
\CL =4AB|\xi'|^2- (B^2+|\xi'|^2)^2,
\end{equation*}
we have
$$L^{-1} = \frac{B^2-A^2}{\CL}\left(\begin{matrix} \dfrac{|\xi'|^2-B^2}{B^2-A^2} 
& -\dfrac{|\xi'|^2(B^2+|\xi'|^2-2AB)}{B^2-A^2} \vspace{1pc} \\ 
-2 &( B^2+|\xi'|^2) \end{matrix}\right).
$$
Thus, we have
\begin{align*}
i\xi'\cdot\bm' & = \frac{1}{\alpha\CL}((|\xi'|^2-B^2)Bi\xi'\cdot\hat\bh'
-|\xi'|^2(B^2+|\xi'|^2-2AB)\hat h_N), \\
\ell & = \frac{B^2-A^2}{\alpha \CL}(-2Bi\xi'\cdot\hat\bh'
 + (B^2+|\xi'|^2)\hat h_N).
\end{align*}
From \eqref{fp3}, \eqref{fp4} and \eqref{f1.27}, we obtain  
\allowdisplaybreaks
\begin{align*}
m_N & = \frac{1}{\alpha \CL }((-|\xi'|^2-B^2 +2AB)i\xi'\cdot\hat\bh'
-A(B^2-|\xi'|^2)\hat h_N), \\
n_j & = \frac{i\xi_j}{\alpha \CL}(-2Bi\xi'\cdot\hat\bh' + (B^2+|\xi'|^2)\hat h_N), \\
n_N & = \frac{A}{\alpha \CL}(2Bi\xi'\cdot\hat\bh' - (B^2+|\xi'|^2)\hat h_N), \\
m_j & = \frac{1}{\alpha B}\hat h_j + \frac{i\xi_j}{\alpha \CL B}
((-3B^2-|\xi'|^2+4AB)i\xi'\cdot\hat\bh'- B(2AB-B^2-|\xi'|^2)\hat h_N).
\end{align*}
By \eqref{Stokes-kernel} we have solution formulas:
\allowdisplaybreaks
\begin{align*}
\CF[u_j](\xi', x_N) & = \frac{1}{\alpha B}\hat h_j e^{-Bx_N} \\
&\qquad + e^{-Bx_N}\frac{i\xi_j}{\alpha \CL B}((-3B^2-|\xi'|^2+4AB)i\xi'\cdot\hat\bh'
- B(2AB-B^2-|\xi'|^2)\hat h_N) \\
&\qquad +M(x_N)\frac{i\xi_j(B-A)}{\alpha \CL}(-2Bi\xi'\cdot\hat\bh' + (B^2+|\xi'|^2)\hat h_N),\\
\CF[u_N](\xi', x_N) & = e^{-Bx_N}\frac{1}{\alpha \CL }((-|\xi'|^2-B^2 +2AB)i\xi'\cdot\hat\bh'
-A(B^2-|\xi'|^2)\hat h_N) \\
&\qquad + M(x_N)\frac{A(B-A)}{\alpha \CL}(2Bi\xi'\cdot\hat\bh' - (B^2+|\xi'|^2)\hat h_N). T
\end{align*}
Using the Volevich idea and the identity: 
$$\pd_NM(x_N) = -e^{-Bx_N} - AM(x_N),$$
we have 
\begin{align}
&u_j(x) \nonumber \\ & = \int^\infty_0\CF^{-1}_{\xi'}\Bigl[Be^{-B(x_N+y_N)}\Bigl\{
\frac{\lambda^{1/2}}{\alpha^2B^3 }\CF'[\lambda^{1/2}h_j](\xi', y_N)
-\sum_{\ell=1}^{N-1}\frac{i\xi_\ell}{\alpha B^3}\CF'[\pd_\ell  h_j](\xi', y_N) \nonumber  \\
&\hskip8.7cm-\frac{1}{\alpha B^2}\CF'[\pd_Nh_j](\xi', y_N)\Bigr\}\Bigr](x')\,dy_N \nonumber \\
& \quad +\int^\infty_0\CF^{-1}_{\xi'}\Bigl[Be^{-B(x_N+y_N)}\Bigl\{\frac{1}{\alpha \CL B}
(-3B^2-|\xi'|^2+4AB)i\xi'\cdot\CF'[\pd_j\bh'](\xi', y_N) \nonumber \\
&\hskip6.2cm
-B(2AB-B^2-|\xi'|^2)\CF'[\pd_jh_N](\xi', y_N)) \nonumber \\
& \hskip2.5cm 
-\frac{i\xi_j}{\alpha \CL B^2}((-3B^2-|\xi'|^2+4AB)i\xi'\cdot\CF'[\pd_N\bh'](\xi', y_N)
\nonumber \\
&\hskip6.3cm-B(2AB-B^2-|\xi|^2)\CF'[\pd_Nh_N](\xi', y_N))\Bigr\}\Bigr](x')\,dy_N
\nonumber \\
&\quad +\int^\infty_0\CF^{-1}_{\xi'}\Bigl[(e^{-B(x_N+y_N)}+AM(x_N+y_N))
\frac{B-A}{\alpha\CL}(-2Bi\xi'\cdot\CF'[\pd_j\bh'](\xi', y_N) \nonumber \\
&\hskip7.7cm+(B^2+|\xi'|^2)\CF'[\pd_jh_N](\xi', y_N)) \nonumber \\
&\qquad 
-M(x_N+y_N)\frac{i\xi_j(B-A)}{\alpha\CL}(-2Bi\xi'\cdot\CF'[\pd_N\bh'](\xi', y_N)
+(B^2+|\xi'|^2)\CF'[\pd_Nh_N](\xi', y_N))\Bigr](x')\,dy_N; \label{solf.j} \\
&u_N(x) \nonumber \\
& = \int^\infty_0 \CF^{-1}_{\xi'}\Bigl[
Be^{-B(x_N+y_N)}\frac{1}{\alpha \CL }((-|\xi'|^2-B^2 +2AB)\CF'[\dv'\bh'](\xi', y_N)
\nonumber \\
&\hskip1cm
-\frac{A\lambda^{1/2}}{\alpha B^2}(B^2-|\xi'|^2)\CF[\lambda^{1/2}h_N](\xi', y_N)
+\sum_{\ell=1}^{N-1}\frac{Ai\xi_\ell}{\alpha B^2}
(B^2-|\xi'|^2)\CF[\pd_\ell h_N](\xi', y_N)) \nonumber \\
&\hskip1cm
-Be^{-B(x_N+y_N)}\frac{1}{\alpha B\CL }((-|\xi'|^2-B^2 +2AB)
i\xi'\cdot\CF'[\pd_N\bh'](\xi', y_N) \nonumber \\
& \hskip1cm 
-A(B^2-|\xi'|^2)\CF'[\pd_N h_N](\xi', y_N)) \nonumber \\
&\hskip1cm +(e^{-B(x_N+y_N)} + AM(x_N+y_N))
\frac{A(B-A)}{\alpha \CL}\bigl\{2B\CF'[\dv'\bh'](\xi', y_N)\nonumber \\
&\hskip1cm
- (B^2+|\xi'|^2)(\frac{\lambda^{1/2}}{\alpha B^2}\CF'[\lambda^{1/2}h_N](\xi', y_N)
-\sum_{\ell=1}^{N-1}\frac{i\xi_\ell}{B^2}\CF'[\pd_\ell h_N](\xi', y_N))\bigr\} \nonumber \\
&\hskip1cm
-M(x_N+y_N)\frac{A(B-A)}{\alpha \CL}(2Bi\xi'\cdot\CF'[\pd_N\bh'](\xi', y_N) -
 (B^2+|\xi'|^2)\CF'[\pd_Nh_N](\xi', y_N))\Bigr](x')\,dy_N. \label{solf.N}
\end{align}
\begin{dfn} A class of $\bN_\ell$ is 
the set of all symbols $m(\lambda, \xi') \in \bM_\ell$ such that 
$\pd_\lambda m(\lambda, \xi') \in \bM_{\ell-2}$. 
\end{dfn}
\begin{lem}\label{lem:9}
Let $\CL=4AB|\xi'|^2- (B^2+|\xi'|^2)^2$ and  $M = \CL/\lambda$.  Then, 
there holds $M^{-1} \in \bN_{-2}$.
\end{lem}
\begin{proof}
Let $r>0$ be a sufficiently small positive number determined later. 
First, we consider the case where
 $|\lambda|/|\xi'|^2 \leq r$ and $(\lambda, \xi') 
\in \Sigma_{\epsilon} \times (\BR^{N-1}\setminus\{0\})$. 
The Taylor expansion of $(1+t)^{1/2}$ is that
$$
(1+t)^{1/2} = 1 + \frac{t}{2} -\frac{t^2}{4}\int^1_0(1-\theta)(1+\theta t)^{-3/2}\,d\theta.
$$
Thus, setting $\gamma_A = (\alpha+\beta)^{-1}$ and $\gamma_B = \alpha^{-1}$, 
for $E \in \{A, B\}$, we have
\begin{equation}\label{c.1}
E = |\xi'|\sqrt{1+\gamma_E \lambda|\xi'|^{-2}}
= |\xi'|\Bigl\{1 + \frac{\gamma_E}{2}\lambda|\xi'|^{-2} 
+\gamma_E^2z_E(\lambda|\xi'|^{-2})^2\Bigr\}
\end{equation}
where we have set 
$$
z_E = -\frac14\int^1_0(1-\theta)(1+\theta \gamma_E\lambda|\xi'|^{-2})^{-3/2}\,d\theta.
$$
Thus, 
\begin{equation}\label{ab.1}
AB = |\xi'|^2\Bigl\{1+\frac12(\gamma_A + \gamma_B)\lambda|\xi'|^{-2}+
y_{AB}\Bigr\},
\end{equation}
where we have set 
\begin{align}y_{AB} &=
\frac{\gamma_A\gamma_B}{4}\lambda^2|\xi'|^{-4} 
+\gamma_A^2 z_A\lambda^2|\xi'|^{-4} 
+\gamma_B^2 z_B\lambda^2|\xi'|^{-4} 
+\frac{\gamma_B\gamma_A^2}{2} z_A\lambda^3|\xi'|^{-6} 
+\frac{\gamma_A\gamma_B^2}{2} z_B\lambda^3|\xi'|^{-6} \nonumber\\
& \qquad + \gamma_A^2\gamma_B^2 z_Az_B\lambda^4|\xi'|^{-8}.
\label{yAB}
\end{align}
First, we shall show that 
\begin{align}
|D^{\kappa'}_{\xi'}z_E| &\leq C_{\kappa'}|\xi'|^{-|\kappa'|}, \label{estop.2*} \\
\label{estop.3*}
|D_{\xi'}^{\kappa'}\pd_\lambda z_E| &\leq C_{\kappa'}|\xi'|^{-2-|\kappa'|}, \\
|D_{\xi'}^{\kappa'}\lambda^{-1}y_{AB}| &\leq C|\xi'|^{-2-|\kappa'|}, \label{estop.4*}\\
|D_{\xi'}^{\kappa'}\pd_\lambda(\lambda^{-1}y_{AB})| &\leq C|\xi'|^{-4-|\kappa'|}
\label{estop.5*}
\end{align}
where $E \in \{A, B\}$.  To prove \eqref{estop.2*}, we observe that 
for any $\kappa'\in \BN_0^{N-1}$, $\ell  \in \BN$,  and $E \in \{A, B\}$, we have 
\begin{equation}\label{estlop.2}
|D_{\xi'}^{\kappa'}(1 + \theta \gamma_E\lambda|\xi|^{-2})^{-\ell/2}|
  \leq C_{\kappa'}|\xi'|^{-|\kappa'|}.
\end{equation}
In fact, by the Bell formula and \eqref{top.1}, we have
\begin{align*}
&|D_{\xi'}^{\kappa'}(1 + \theta\gamma_E\lambda|\xi'|^{-2})^{-\ell/2}| \\
 &\quad \leq C_{\kappa'}\sum_{m=1}^{|\kappa'|}
|1+\theta \gamma_E\lambda|\xi'|^{-2}|^{-\ell/2 - m}\sum_{\kappa'_1+\cdots + \kappa'_m = \kappa'}
|D_{\xi'}^{\kappa'_1}\theta \gamma_E \lambda |\xi'|^{-2}| \cdots |D_{\xi'}^{\kappa'_m}\theta 
\gamma_E \lambda|\xi'|^{-2}| \\
& \quad \leq C_{\kappa'}\sum_{m=1}^{|\kappa'|}
|1+\theta \gamma_E\lambda|\xi'|^{-2}|^{-\ell/2 - m} |\theta \gamma_E \lambda |\xi'|^{-2}|^{m}
|\xi'|^{-|\kappa'|} \\
& \quad \leq C_{\kappa'}|\xi'|^{-|\kappa'|}.
\end{align*}
Thus, we have \eqref{estlop.2}, and so \eqref{estop.2*}. Since 
$$\pd_\lambda z_E = \frac{3}{8}\int^1_0(1-\theta)(1+\theta\gamma_E\lambda|\xi'|^{-2})^{-5/2}
\theta\,d\theta \gamma_E|\xi'|^{-2},$$
by Leibniz's rule and \eqref{estlop.2}, we have \eqref{estop.3*}. 
From \eqref{yAB}, it follows that 
\begin{align*}
\lambda^{-1}y_{AB} &=
\frac{\gamma_A\gamma_B}{4}\lambda|\xi'|^{-4} 
+\gamma_A^2 z_A\lambda|\xi'|^{-4} 
+\gamma_B^2 z_B\lambda|\xi'|^{-4} 
+\frac{\gamma_B\gamma_A^2}{2} z_A\lambda^2|\xi'|^{-6} 
+\frac{\gamma_A\gamma_B^2}{2} z_B\lambda^2|\xi'|^{-6} \nonumber\\
& \qquad + \gamma_A^2\gamma_B^2 z_Az_B\lambda^3|\xi'|^{-8}.
\end{align*}
Thus, using \eqref{estop.2*} and Leibniz's rule, we have
\begin{align*}
|D_{\xi'}^{\kappa'}\lambda^{-1}y_{AB}|&\leq C_{\kappa'}(|\lambda||\xi'|^{-4-|\kappa'|}
+ |\lambda|^2|\xi'|^{-6-|\kappa'|}+ |\lambda|^3|\xi'|^{-8-|\kappa'|})\\
&\leq C_{\kappa'}(r+r^2+r^3)|\xi'|^{-2-|\kappa'|},
\end{align*}
which shows \eqref{estop.4*}. Since
\begin{align*}
\pd_\lambda(\lambda^{-1}y_{AB})&=
\frac{\gamma_A\gamma_B}{4}|\xi'|^{-4} 
+\gamma_A^2 ((\pd_\lambda z_A)\lambda + z_A) |\xi'|^{-4} 
+\gamma_B^2 ((\pd_\lambda z_B)\lambda + z_B)|\xi'|^{-4} \\
& \qquad 
+\frac{\gamma_B\gamma_A^2}{2} ((\pd_\lambda z_A)\lambda^2+ 2z_A\lambda)|\xi'|^{-6} 
+\frac{\gamma_A\gamma_B^2}{2} ((\pd_\lambda z_B)\lambda^2+2z_B\lambda)|\xi'|^{-6} \\
& \qquad + \gamma_A^2\gamma_B^2 ((\pd_\lambda z_A)z_B\lambda^3
+z_A(\pd_\lambda z_B)\lambda^3 + 3z_Az_B\lambda^2)|\xi'|^{-8}.
\end{align*}
Thus, by Leibniz's rule, \eqref{estop.2*} and \eqref{estop.3*}, we have
\begin{align*}
|D^{\kappa'}_{\xi'}(\lambda^{-1}y_{AB})| & \leq C_{\kappa'}(|\xi'|^{-4-|\kappa'|}
+ |\lambda||\xi'|^{-6-|\kappa'|} + |\lambda|^2|\xi'|^{-8-|\kappa'|}
+ |\lambda|^3|\xi'|^{-10-|\kappa'|}) 
\leq C_{\kappa'}|\xi'|^{-4-|\kappa'|}.
\end{align*}
This shows \eqref{estop.5*}. \par
Since $(B^2+|\xi'|^2)^2 = (\alpha^{-1}\lambda + 2|\xi'|^2)^2 
= \alpha^{-2}\lambda^2 + 4\gamma_B\lambda|\xi'|^2 + 4|\xi'|^4$, we have
$$\CL = 4AB|\xi'|^2 - (B^2+|\xi'|^2)^2 
= 2c_0\lambda|\xi'|^2 -\alpha^{-2}\lambda^2+ 4|\xi'|^4y_{AB}$$
with $c_0 = \gamma_A-\gamma_B = (\alpha+\beta)^{-1}-\alpha^{-1}$. 
Thus, we have
\begin{equation}\label{formM.1}
M = \frac{\CL}{\lambda} = 2c_0|\xi'|^2-\alpha^{-2}\lambda + 4|\xi'|^4\lambda^{-1}y_{AB}.
\end{equation}
Thus, by \eqref{estop.4*} and Leibniz's rule, we have
\begin{equation}\label{yform.1}
|D_{\xi'}^{\kappa'} M| \leq C_{\kappa'}|\xi'|^{2-|\kappa'|}.
\end{equation}
In particular,
$|M| \geq 2|c_0||\xi'|^2- C(|\lambda| + |\lambda|^2|\xi'|^{-2} + |\lambda|^3|\xi'|^{-4})
\geq 2|c_0| - C|\xi'|^2(r + r^2 + r^3)$.  Thus, choosing $r>0$ so small that 
$C'(r+r^2+r^3) \leq |c_0|$, we have
\begin{equation}\label{lop.1}
|M|\geq |c_0||\xi'|^2.
\end{equation}
By \eqref{yform.1}, \eqref{lop.1} and Bell's formula, we have
\begin{equation}\label{first.1}\begin{aligned}
|D_{\xi'}^{\kappa'}M^{-1}| &\leq C_{\kappa'}\sum_{\ell=1}^{|\kappa'|}
|M|^{-\ell-1}\sum_{\kappa_1'+\cdots+\kappa_\ell'=\kappa'} |D_{\xi'}^{\kappa_1'}M|
\cdots |D_{\xi'}^{\kappa_\ell'}M| \leq C_{\kappa'}|\xi'|^{-2-|\kappa'|}.
\end{aligned}\end{equation}
\par
We now consider $\pd_\lambda M^{-1} = -M^{-2}\pd_\lambda M$. From \eqref{formM.1}
it follows that 
$$\pd_\lambda M = -\alpha^{-2} + 4|\xi'|^4\pd_\lambda(\lambda^{-1}y_{AB}).$$
Thus, by \eqref{estop.5*}, we have
\begin{equation}\label{first.2}
|D^{\kappa'}_{\xi'}(\pd_\lambda M)| \leq C|\xi'|^{-|\kappa'|}.
\end{equation}
Since $\pd_\lambda M^{-1}=-M^{-2}\pd_\lambda M$, by Leibniz's rule, \eqref{first.1},
and \eqref{first.2}, we have
\begin{equation}\label{first.2.2}
|D_{\xi'}^{\kappa'}\pd_\lambda M^{-1}| \leq C_{\kappa'}|\xi|^{-4-|\kappa'|}.
\end{equation}
Noting that $|\xi'| \geq (1/2)(|\xi'| + r^{-1}|\lambda|^{1/2})$, by \eqref{first.1}
and \eqref{first.2.2}, we have
\begin{equation*}
\begin{aligned}
|D_{\xi'}^{\kappa'}M^{-1}| & \leq C(|\lambda|^{1/2}+|\xi'|)^{-2-|\kappa'|},\\
|D_{\xi'}^{\kappa'}(\pd_\lambda M^{-1})| & \leq C(|\lambda|^{1/2}+|\xi'|)^{-4-|\kappa'|}.
\end{aligned}\end{equation*}
\par
Next, we consider the case where $|\xi'|^2/|\lambda| \leq r$  
and $(\lambda, \xi') \in \Sigma_{\epsilon} \times (\BR^{N-1}\setminus\{0\})$. 
In this case, we have
\begin{align*}
A&=((\alpha+\beta)^{-1}\lambda)^{1/2}(1+|\xi'|^2(\alpha+\beta)\lambda^{-1})^{1/2},\\
B & = (\alpha^{-1}\lambda)^{1/2}(1+|\xi'|^2\alpha\lambda^{-1})^{1/2},
\end{align*}
and so 
$$\CL= 4(\alpha+\beta)^{-1/2}\alpha^{-1/2}\lambda|\xi'|^2(1 +
O(|\xi'|^2/\lambda) )-\alpha^{-2}\lambda^2 -4\alpha^{-1}\lambda|\xi'|^2
-4|\xi'|^4
$$
This implies that 
$$|\CL| \geq \alpha^{-2}|\lambda|^2 - Cr|\lambda|^2$$
for some constant $C$.  Choosing $r>0$ so small that 
$Cr \leq \alpha^{-2}/2$, we have 
$$
|\CL| \geq (\alpha^{-2}/2)|\lambda|^{2}.
$$
Thus, we have
\begin{equation*}
|M| \geq (\alpha^{-2}/2)|\lambda| \geq c_1(|\lambda|^{1/2}+|\xi'|)^2
\end{equation*}
with some positive number $c_1$, where we have used the fact that
$|\lambda|^{1/2}\geq (1/2)(|\lambda|^{1/2} + r^{-1/2}|\xi'|)$.  
We know that
$$|D^{\kappa'}|\xi'|^2| \leq C_{\kappa'}(|\lambda|^{1/2}+|\xi'|)^{2-|\kappa'|}.$$
Moreover, for $E \in \{A, B\}$ we know that 
\begin{align*}
|D^{\kappa'}_{\xi'}E| &\leq C(|\lambda|^{1/2}+|\xi'|)^{1-|\kappa'|}, \\
|D^{\kappa'}_{\xi'}(\pd_\lambda E)| &\leq C(|\lambda|^{1/2}+|\xi'|)^{-1-|\kappa'|}.
\end{align*}
Thus, we have 
$$|D^{\kappa'}_{\xi'}\CL| \leq C(|\lambda|^{1/2}+|\xi'|)^{4-|\kappa'|}.$$
Recalling that $M=\CL/\lambda$, we have
$$ |D^{\kappa'}_{\xi'}M| \leq C(|\lambda|^{1/2}+|\xi'|)^{2-|\kappa'|}. $$
By Bell's formular, 
\begin{align*}
|D_{\xi'}^{\kappa'}M^{-1}|
&\leq C_{\kappa'} \sum_{\ell=1}^{|\kappa'|}
(|\lambda|^{1/2}+|\xi'|)^{-2(\ell+1)}(|\lambda|^{1/2}+|\xi'|)^{2\ell-|\kappa'|}. 
\end{align*}
Namely, we have 
\begin{equation}\label{first.3}
|D_{\xi'}^{\kappa'}M^{-1}| 
\leq C_{\kappa'}(|\lambda|^{1/2}+|\xi'|)^{-2-|\kappa'|}.
\end{equation}
We observe that $\pd_\lambda M^{-1} = -M^{-2}(\pd_\lambda M)$ and 
$\pd_\lambda M = \lambda^{-1}\pd_\lambda \CL - \lambda^{-2}\CL$.  Since
$\pd_\lambda \CL = 4((\pd_\lambda A)B+A(\pd_\lambda B))|\xi'|^2 - 2(\alpha^{-1}\lambda
+2|\xi'|^2)\alpha^{-1}$, we have
$$|D_{\xi'}^{\kappa'}\pd_\lambda \CL| \leq C_{\kappa'}(|\lambda|^{1/2}+|\xi|)^{2-|\kappa'|},
$$
and so
$$|D_{\xi'}^{\kappa'}(\pd_\lambda M)| \leq C_{\kappa'}(|\lambda|^{1/2}+|\xi'|)^{-|\kappa'|}.$$
Therefore, we have
\begin{equation}\label{first.4}
|D^{\kappa'}_{\xi'}(\pd_\lambda M^{-1})| \leq C_{\kappa'}(|\lambda|^{1/2}+|\xi'|)^{-4-|\kappa'|}.
\end{equation}
\par
Finally, we consider the case where 
\begin{equation}\label{lop:5}
r|\lambda| \leq |\xi'|^2 \leq r^{-1}|\lambda|, 
\quad (\lambda, \xi') \in \Sigma_{\epsilon, \lambda_0}\times \BR^{N-1}.
\end{equation}
First, we shall show that 
\begin{equation}\label{nonzero}
\CL\not=0 \quad\text{
for any $(\lambda, \xi') \in \Sigma_\epsilon\times(\BR^{N-1}\setminus\{0\}).$}
\end{equation}
We know that the uniqueness of solutions of the forms $\hat u_J = a_J e^{-Bx_N} + 
b_J(e^{-Bx_N} -e^{-Ax_N})$ to equations \eqref{fp1} guarantees $\CL\not=0$.  
Thus, we shall show the uniqueness.
 Let us write $\hat{u}_J = v_J$ for the notational simplicity and 
we assume that $v_J$ satisfy equations:
\begin{equation*}
\begin{aligned}
&(\lambda +\alpha |i\xi'|^2)v_j-\alpha \pd_N^2v_j
-\beta i\xi_j(i\xi'\cdot \bv' + \pd_Nv_N) =0 \qquad(x_N>0), \\
&(\lambda +\alpha |i\xi'|^2)v_N-\alpha \pd_N^2v_N
-\beta \pd_N(i\xi'\cdot \bv' + \pd_Nv_N) =0 \quad(x_N>0), \\
&(\pd_Nv_j+ i\xi_j v_N)|_{x_N=0} = 0, \\
&(\alpha (\pd_Nv_N - i\xi'\cdot\bv') + \beta(i\xi'\cdot\bv'+\pd_Nv_N))|_{x_N=0} =0.
\end{aligned}\end{equation*}
Here, $\bv' = (v_1, \ldots, v_{N-1})$. Notice that 
\begin{align}
0 &= \lambda v_j +\alpha |\xi'|^2v_j -\alpha \pd_N^2v_j
-\beta i\xi_j(i\xi'\cdot\bv'+\pd_Nv_N)
\nonumber \\ 
&=\lambda v_j - \alpha \sum_{k=1}^{N-1}i\xi_k(i\xi_kv_j + i\xi_jv_k) - \alpha
\pd_N(\pd_Nv_j + i\xi_jv_N) -(\beta-\alpha)i\xi_j(i\xi'\cdot\bv' + \pd_Nv_N),  \label{part.j}\\
0 &= \lambda v_N +\alpha |\xi'|^2v_N -\alpha \pd_N^2v_N
-\beta \pd_N(i\xi'\cdot\bv'+\pd_Nv_N) \nonumber 
\\
&=\lambda v_N - \alpha \sum_{k=1}^{N-1}i\xi_k(i\xi_k v_N + \pd_N v_k)-
2\alpha \pd_N^2v_N -(\beta-\alpha)\pd_N(i\xi'\cdot\bv' + \pd_Nv_N).
\label{part.N}
\end{align}
Since $v_j$ decay exponentially as $x_N\to\infty$, multiplying \eqref{part.j}
with $\overline{v_j}$ and \eqref{part.N} with $\overline{v_N}$, and using integration by parts
\begin{equation}\label{energy:1}\begin{aligned}
0=&\lambda\|\bv\|^2 + \alpha\sum_{j,k=1}^{N-1} \|i\xi_kv_j\|^2 
+ \alpha \sum_{j=1}^{N-1}\|\pd_Nv_j\|^2 + \alpha \|i\xi'\cdot\bv'\|^2 \\
&+ \alpha 
\sum_{j=1}^{N-1}(i\xi_jv_N, \pd_Nv_j) +\alpha\sum_{j=1}^{N-1}\|i\xi_j v_N\|^2
+ \alpha\sum_{j=1}^{N-1}(\pd_Nv_j, i\xi_j v_N) + 2\alpha\|\pd_Nv_N\|^2 \\
&+(\beta-\alpha)(\|i\xi'\cdot\bv'\|^2 + (\pd_Nv_N, i\xi'\cdot\bv')
+ (i\xi'\cdot\bv', \pd_Nv_N) + \|\pd_Nv_N\|^2).
\end{aligned}\end{equation}
Here $(f, g) = \int_{\BR_+}f(x_N)\overline{g(x_N)}\,dx_N$ and $\|f\|^2 = (f, f)$.   Taking
the real part and the imaginary part of \eqref{energy:1}  yields
\begin{align}
&({\rm Im}\,\lambda)\|\bv\|^2=0, \label{energy.2}\\
&({\rm Re}\,\lambda)\|\bv\|^2 + \alpha \sum_{j,k=1}^{N-1}\|i\xi_kv_j\|^2 + 
\alpha \sum_{j=1}^{N-1}(\|i\xi_jv_N\|^2 + \|\pd_Nv_j\|^2 
+ (i\xi_jv_N, \pd_Nv_j) + (\pd_Nv_j, i\xi_jv_N))\nonumber \\
&+\alpha\|i\xi'\cdot\bv'\|^2 + 2\alpha\|\pd_Nv_N\|^2
+ \beta(\|i\xi'\cdot \bv'\|^2 + \|\pd_Nv_N\|^2 
+ (\pd_Nv_N, i\xi'\cdot\bv') + (i\xi'\cdot\bv', \pd_Nv_N))\nonumber \\
&-\alpha(\|i\xi'\cdot\bv'\|^2 + (\pd_Nv_N, i\xi'\cdot\bv')
+(i\xi'\cdot\bv', \pd_Nv_N) +\|\pd_Nv_N\|^2)=0. 
\label{energy.3}
\end{align}
If ${\rm Im}\, \lambda \not=0$, by \eqref{energy.2}, we have $\bv=0$.  
If ${\rm Im}\,\lambda = 0$, then ${\rm Re}\,\lambda >0$.  Thus, by \eqref{energy.3}, 
we have
\begin{align*}
0 &\geq {\rm Re}\,\lambda \|\bv\|^2 + \alpha \sum_{j,k=1}^{N-1}\|i\xi_kv_j\|^2
+ \alpha\sum_{j=1}^{N-1}(\|\pd_Nv_j\|^2 -2\|i\xi_j v_N\|\|\pd_N v_j\| 
+\|i\xi_j v_N\|^2)\\
&+\alpha\|i\xi'\cdot\bv'\|^2 + 2\alpha\|\pd_Nv_N\|^2
-\alpha\|i\xi'\cdot\bv'\|^2 - \alpha\|\pd_Nv_N\|^2
-2\alpha\|\pd_Nv_N\|\|i\xi'\cdot\bv'\|\\
&+ \beta(\|i\xi'\cdot \bv'\|^2 +\|\pd_Nv_N\|^2 
-2 \|\pd_Nv_N\|\|i\xi'\cdot\bv'\| ).
\end{align*}
Here, we have 
\begin{align*}
&\alpha(\|\pd_Nv_j\|^2-2\|i\xi_j v_N\|\|\pd_N v_j\| +\|i\xi_j v_N\|^2)
\geq \alpha(\|\pd_Nv_j\|-\|i\xi_jv_N\|)^2 \geq 0, \\
 &\beta(\|i\xi'\cdot \bv'\|^2 + \|\pd_Nv_N\|^2 
-2\|\pd_Nv_N\|\|i\xi'\cdot\bv'\|
\geq \beta(\|i\xi'\cdot \bv'\|-\|\pd_Nv_N\|)^2 \geq 0.
\end{align*}
Moreover, noting that $\sum_{j,k=1}^{N-1}\|i\xi_kv_j\|^2 
= |\xi'|^2\|\bv'\|^2$ and $\|i\xi'\cdot\bv'\| \leq |\xi'|\|\bv'\|$, we have
\begin{align*}
&\alpha \sum_{j,k=1}^{N-1}\|i\xi_kv_j\|^2
+\alpha\|i\xi'\cdot\bv'\|^2 + 2\alpha\|\pd_Nv_N\|^2
-\alpha\|i\xi'\cdot\bv'\|^2 - \alpha\|\pd_Nv_N\|^2
-2\alpha\|\pd_Nv_N\|\|i\xi'\cdot\bv'\|\\
&\geq \alpha(|\xi'|^2\|\bv'\|^2 + \|\pd_Nv_N\|^2 - 2\|\pd_Nv_N\| |\xi'|\|\bv'\|)
\geq \alpha(|\xi'|\|\bv'\|-\|\pd_Nv_N\|)^2 \geq 0.
\end{align*}
Thus, we have $0 \geq {\rm Re}\,\lambda \|\bv\|^2$, which imples that $v=0$.
In particular, we have proved \eqref{nonzero}.

\par

Let 
\begin{gather*}
\tilde\lambda = \frac{\lambda}{(|\lambda|^{1/2}+ |\xi'|)^2}, \quad
\tilde \xi' = \frac{\xi'}{|\lambda|^{1/2}+ |\xi'|}, \\
\tilde A= ((\alpha + \beta)^{-1}\tilde\lambda 
+|\tilde\xi'|^2)^{1/2}, \quad 
\tilde B = (\alpha^{-1}\tilde\lambda + |\tilde\xi'|^2)^{1/2}, \\
\tilde\CL =4\tilde A\tilde B |\tilde\xi'|^2
-(\tilde B^2+ |\tilde\xi'|^2)^2.
\end{gather*}
We have $\CL= (|\lambda|^{1/2} + |\xi'|)^4 \tilde \CL$.  Notice that
$\tilde\lambda \in \Sigma_\epsilon$ provided that $\lambda \in \Sigma_\epsilon$.
Moreover, from \eqref{lop:5} it follows that 
\begin{equation}\label{lop:7}
\frac{r^{1/2}}{1+r^{1/2}} \leq |\tilde \xi| \leq \frac{1}{1+r^{1/2}}, \quad 
\frac{r}{(1+r^{1/2})^2} \leq |\tilde \lambda| \leq \frac{1}{(1+r^{1/2})^2}.
\end{equation}
Thus, we set
$$U = \{(\tilde \lambda, \tilde\xi') \in \Sigma_\epsilon \times (\BR^{N-1}\setminus\{0\}) 
\mid \text{\eqref{lop:7} holds}\}.
$$
Since $U$ is a compact set and $\tilde L \not=0$
for $(\tilde\lambda, \tilde\xi') \in U$, 
there exists a constant $c_2>0$ such that
$$\inf_{(\tilde\lambda, \tilde\xi') \in U} 
|\tilde \CL| = c_2> 0.
$$
Thus, we have 
\begin{equation*}
|\CL| \geq c_2(|\lambda|^{1/2}+|\xi'|)^4
\end{equation*}
for $(\lambda, \xi') \in \Sigma_\epsilon \times(\BR^{N-1}\setminus\{0\})$ satisfying
the condition \eqref{nonzero}. 
Employing the same argument as in the proof of \eqref{first.3} and \eqref{first.4}, 
we have
\begin{equation*}
\begin{aligned}
|D^{\kappa'}_{\xi'}M^{-1}| &\leq C_{\kappa'}(|\lambda|^{1/2} + |\xi'|)^{-2-|\kappa'|}, \\
|D^{\kappa'}_{\xi'}(\pd_\lambda M^{-1})| &\leq C_{\kappa'}(|\lambda|^{1/2} + |\xi'|)^{-4-|\kappa'|}.
\end{aligned}\end{equation*}
Thus, we have proved $M^{-1} \in \bM_{-2}$ 
and $\pd_\lambda M^{-1} \in \bM_{-4}$. Namely, $M^{-1} \in \bN_{-2}$.
\end{proof} 
\begin{lem}\label{lem:10} Let $m_1=-3B^2-|\xi'|^2+4AB$, $m_2=2AB-B^2-|\xi'|^2$ and 
$m_3=B-A$ be symbols appearing in the solution formula.  Then, 
$m_i/\lambda \in \bN_0$ for $i=1, 2$ and $m_3/\lambda \in \bN_{-1}$. 
\end{lem}
\begin{proof} To prove the lemma, we use the symbols given in the proof of Lemma \ref{lem:9}.
First, we consider the case where $|\lambda|/|\xi'|^2\leq r$ and $(\lambda, \xi')
\in \Sigma_\epsilon\times(\BR^{N-1}\setminus\{0\})$ with some small $r>0$.
Using \eqref{ab.1}, 
we write $m_1/\lambda= (2\gamma_A-\gamma_B) + 4|\xi'|^2\lambda^{-1}y_{AB}$. 
By \eqref{estop.4*}, \eqref{estop.5*} and Leibniz's rule, we have 
\begin{align*}
|D_{\xi'}^{\kappa'}(m_1/\lambda)| &\leq C_{\kappa'}|\xi'|^{-|\kappa'|}, \\
|D_{\xi'}^{\kappa'}\pd_\lambda (m_1/\lambda)| &\leq C_{\kappa'}|\xi'|^{-2-|\kappa'|}.
\end{align*}
Analogously, using \eqref{ab.1} we write 
$m_2/\lambda = \gamma_A + 2|\xi'|^2\lambda^{-1}y_{AB}$. 
By \eqref{estop.4*}, \eqref{estop.5*} and Leibniz's rule, we have
\begin{align*}
|D_{\xi'}^{\kappa'}(m_2/\lambda)| &\leq C_{\kappa'}|\xi'|^{-|\kappa'|}, \\
|D_{\xi'}^{\kappa'}\pd_\lambda (m_2/\lambda)| &\leq C_{\kappa'}|\xi'|^{-2-|\kappa'|}.
\end{align*}
Using \eqref{c.1}, we have $m_3/\lambda = (\gamma_B-\gamma_A)|\xi'|^{-1}/2 
+(\gamma_B^2z_B-\gamma_A^2z_A)\lambda|\xi'|^{-3}$.
Since 
\begin{gather*}
|D^{\kappa'}_{\xi'}(\lambda|\xi'|^{-3})| \leq C_{\kappa'}|\lambda||\xi'|^{-3-|\kappa'|} 
\leq C_{\kappa'}r|\xi'|^{-1-|\kappa'|}, \\
|D^{\kappa'}_{\xi'}\pd_\lambda(\lambda|\xi'|^{-3})| \leq C_{\kappa'}|\xi'|^{-3-|\kappa'|}.
\end{gather*}
Using this, \eqref{estop.2*}, \eqref{estop.3*} and Leibniz's rule, we have
\begin{align*}
|D_{\xi'}^{\kappa'}(m_3/\lambda)| &\leq C_{\kappa'}|\xi'|^{-1-|\kappa'|}, \\
|D_{\xi'}^{\kappa'}\pd_\lambda (m_3/\lambda)| &\leq C_{\kappa'}|\xi'|^{-3-|\kappa'|}.
\end{align*}
\par
We now consider the case where $|\lambda|/|\xi'|^2 \geq r$ and $(\lambda, \xi')
\in \Sigma_\epsilon\times(\BR^{N-1}\setminus\{0\})$.  In view of  \eqref{top.1}, 
we see easily that 
\begin{align*}
|D^{\kappa'}_{\xi'} E|\leq C_{\kappa'}(|\lambda^{1/2}+|\xi'|)^{1-|\kappa'|}, \\
|D^{\kappa'}_{\xi'} \pd_\lambda E|\leq C_{\kappa'}(|\lambda^{1/2}+|\xi'|)^{-1-|\kappa'|} 
\end{align*}
for $E \in \{A, B\}$. Noticing that $|\lambda|^{1/2} \geq (1/2)(|\lambda|^{1/2}
+ r^{1/2}|\xi'|)$, we see easily that 
\begin{align*}
|D^{\kappa'}_{\xi'}(m_i/\lambda)| \leq C_{\kappa'}(|\lambda|^{1/2}+|\xi'|)^{-|\kappa'|}, \\
|D^{\kappa'}_{\xi'}\pd_\lambda(m_i/\lambda)| \leq C_{\kappa'}(|\lambda|^{1/2}+|\xi'|)^{-2-|\kappa'|}
\end{align*}
for $i=1,2$.  And 
\begin{align*}
|D^{\kappa'}_{\xi'}(m_3/\lambda)| \leq C_{\kappa'}(|\lambda|^{1/2}+|\xi'|)^{-1-|\kappa'|}, \\
|D^{\kappa'}_{\xi'}\pd_\lambda(m_3/\lambda)| \leq C_{\kappa'}(|\lambda|^{1/2}+|\xi'|)^{-3-|\kappa'|}. 
\end{align*}
Combining these results, we have proved that
$m_i/\lambda \in \bN_0$ ($i=1,2$) and $m_3/\lambda \in \bN_{-1}$
This  completes the proof.
\end{proof}

\subsection{Estimates of solution formulas in the half-space}

Let $\tilde H = (H_2, H_3)$ be an $N+N^2$ vertical vector such that 
$H_2$ and $H_3$ correspond to $\lambda^{1/2}\bh= (\lambda^{1/2}h_1,
\ldots, \lambda^{1/2}h_N)^\top$ and 
$\nabla \bh = (\pd_ih_j \mid i, j=1, \ldots, N)^\top$, respectively. 
Here, $\bk^\top$ denotes the transposed vector of $\bk=(k_1,\ldots, k_N)$. 
In the solution formulas \eqref{solf.j} and \eqref{solf.N}, writing 
\begin{alignat*}2
\frac{-3B^2-|\xi'|^2+4AB}{\CL} &= \frac{-3B^2-|\xi'|^2+4AB}{\lambda M},
& \quad
\frac{2AB-B^2-|\xi'|^2}{\CL} &= \frac{2AB-B^2-|\xi'|^2}{\lambda M}, \\
\frac{B^2-|\xi'|^2}{\CL} &= \frac{1}{\alpha M},
&\quad  \frac{B-A}{\CL} &= \frac{B-A}{\lambda M},
\end{alignat*}
by Lemmas \ref{lem:9} and \ref{lem:10}, we see that the first three terms  belong 
to $\bN_{-2}$ and the last term belongs to $\bN_{-3}$. 
 Thus, all the symbols appearing in the solutions formulas
\eqref{solf.j} and \eqref{solf.N} belong to $\bN_{-2}$.  Using this fact, we see that
there exist $N\times N$ matrices $\CM_1(\lambda, \xi')$ and 
$\CM_2(\lambda, \xi')$ of $\bN_{-2}$ symbols such that setting
\begin{align*}
\CT(\lambda)\tilde H &= \int^\infty_0 \CF^{-1}_{\xi'}[Be^{-B(x_N+y_N)}
\CM_1(\lambda, \xi')\CF'[\tilde H](\xi', y_N)](x')\,dy_N \\
& \quad + \int^\infty_0 \CF^{-1}_{\xi'}[B^2M(x_N+y_N) \CM_2(\lambda, \xi')
\CF'[\tilde H](\xi', y_N)](x')\,dy_N
\end{align*}
we have 
 $\bu = \CT(\lambda)(\lambda^{1/2}\bh, \nabla\bh)$ is a solution of 
equations \eqref{free.1}.  Here, $\CT(\lambda)\tilde H$ and $\tilde H$ are
both vertical vectors.  Moreover,  $\CT(\lambda)$ has the following 
properties. 
\begin{thm}\label{thm:3.1} Let $1 < q < \infty$, $-1+1/q <s < 1/q$, 
$0 < \epsilon < \pi/2$, $\lambda_0 > 0$ and  let $\sigma$ be a small positive 
number such that $-1+1/q < s-\sigma < s < s+\sigma < 1/q$. 
Then, $\CT(\lambda) \in {\rm Hol}\,(\Sigma_{\epsilon, \lambda_0},
\CL(B^s_{q,1}(\HS)^{N+N^2}, B^{s+2}_{q,1}(\HS)^N))$ has 
$(s, \sigma, q)$ properties. 
\end{thm}
\begin{proof} To prove the theorem, we shall use Theorem \ref{thm:5.2}
 in Subsection \ref{subsec:2.5}. 
Notice that 
$\|f\|_{W^1_q(\HS)} = \|\bar\nabla f\|_{L_q(\HS)}$ and $\|f\|_{W^2_q(\HS)} = 
\|\bar\nabla^2f\|_{L_q(\HS)}$.  In what follows, 
we may assume that $\tilde H \in C^\infty_0(\HS)^{N+N^2}$.  In fact, 
$C^\infty_0(\HS)$ is dense in $B^s_{q,1}(\HS)$ for $1 < q < \infty$ and 
$-1+1/q < s < 1/q$ (cf. Proposition 2.24, Lemma 2.32, and
Corollaries 2.26 and 2.34 in \cite{Gaudin}). 
Using the formulas:
$$\pd_N^\ell M(x_N ) = (-1)^\ell(A^\ell M(x_N) +
 \frac{A^\ell-B^\ell}{A-B}e^{-Bx_N})\quad(\ell \geq 1),\\
$$
and setting 
\begin{align*}
\CM_1^{(0)}(\lambda) &= \CM_1(\lambda,\xi'), \quad 
\CM_1^{(\ell)}(\lambda)= (-B)^\ell \CM_1(\lambda, \xi') 
+ (-1)^\ell \frac{A^\ell-B^{\ell}}{A-B}B \CM_2(\lambda,\xi') \quad(\ell \geq 1), 
\\
\CM_2^{(0)}(\lambda) &= \CM_2(\lambda, \xi'), \quad
\CM_2^{(\ell)}(\lambda) = (-1)^\ell A^\ell \CM_2(\lambda,\xi') \quad(\ell \geq 1)
\end{align*}
for the notational simplicity, 
we write
\begin{align}
\pd_N^\ell \CT(\lambda)\tilde H = 
\int^\infty_0\CF^{-1}_{\xi'}&\Bigl[\CM_{1}^{(\ell)}(\lambda)
\CF'[\tilde H](\xi', y_N)Be^{-B(x_N+y_N)} 
\nonumber \\
&+ 
\CM_{2}^{(\ell)}(\lambda)\CF'[\tilde H](\xi', y_N)B^2M(x_N+y_N)\Bigr](x')\,dy_N.
\label{dif.1.1}
\end{align}
Using these symbols, we can write
\begin{align*}
\lambda^{k}\pd_{x'}^{\kappa'}\pd_N^\ell \CT(\lambda)\tilde H = 
\int^\infty_0\CF^{-1}_{\xi'}\Bigl[&\lambda^{k}(i\xi')^{\kappa'}\CM_{1}^{(\ell)}(\lambda)
\CF'[\tilde H](\xi', y_N)Be^{-B(x_N+y_N)} \\
+ &\lambda^{k}(i\xi')^{\kappa'}\CM_{2}^{(\ell)}(\lambda)\CF'[\tilde H](\xi', y_N)
B^2M(x_N+y_N)\Bigr](x')\,dy_N.
\end{align*}
If $2k+|\kappa'| + \ell\leq 2$, then $\lambda^{k}(i\xi')^{\kappa'}\CM_{1}^{(\ell)}(\lambda)
\in \bM_0$ and $\lambda^{k}(i\xi')^{\kappa'}\CM_{2}^{(\ell)}(\lambda) \in \bM_0$.  Thus, 
by Proposition \ref{prop:2} we have
\begin{equation}\label{5.3}
\|(\lambda, \lambda^{1/2}\bar\nabla, \bar\nabla^2)\CT(\lambda)\tilde H \|_{L_q(\HS)} 
\leq C\|\tilde H\|_{L_q(\HS)}.
\end{equation}
To obtain the $W^1_q(\HS)$ estimate, noting that $\tilde H \in C^\infty_0(\HS)^{N+N^2}$,
using the formulas:
$$
\pd_N (-B)^{-1} e^{-B(x_N+y_N)} =e^{-B(x_N+y_N)}, \quad 
\pd_N ((AB)^{-1}e^{-B(x_N+y_N)} - A^{-1}M(x_N+y_N)) = M(x_N+y_N) 
$$
and setting 
$$\tilde\CM^{(\ell-1)}_1(\lambda) 
= (B^{-1}\CM^{(\ell)}_1(\lambda)- A^{-1}\CM^{(\ell)}_2(\lambda)), \quad 
\tilde \CM^{(\ell-1)}_2(\lambda) = A^{-1}\CM^{(\ell)}_2(\lambda), 
$$
by integration
by parts, we have 
\begin{align}
\pd_N^\ell \CT(\lambda)\tilde H = 
\int^\infty_0\CF^{-1}_{\xi'}\Bigl[
&\tilde\CM^{(\ell-1)}_1(\lambda)\CF'[\pd_N\tilde H](\xi', y_N)Be^{-B(x_N+y_N)} \nonumber \\
&+ \tilde\CM^{(\ell-1)}_2(\lambda)\CF'[\pd_N\tilde H](\xi', y_N)B^2M(x_N+y_N)\Bigr](x')\,dy_N.
\label{dif.3}
\end{align}
Thus, we have
\begin{align*}
\lambda^{k}\pd_{x'}^{\kappa'}\pd_N^\ell \CT(\lambda)\tilde H = 
\int^\infty_0\CF^{-1}_{\xi'}\Bigl[
&\lambda^{k}(i\xi')^{\kappa'} \tilde\CM^{(\ell-1)}_1(\lambda)
\CF'[\pd_N\tilde H](\xi', y_N)Be^{-B(x_N+y_N)} \\
&+\lambda^{k}(i\xi')^{\kappa'}\tilde\CM_2^{(\ell-1)}(\lambda)
\CF'[\pd_N\tilde H](\xi', y_N)B^2M(x_N+y_N)\Bigr](x')\,dy_N.
\end{align*}
If $2k+|\kappa'| + \ell \leq 3$, both  
$\lambda^{k}(i\xi')^{\kappa'}\tilde\CM^{(\ell-1)}_1(\lambda)$ and 
$\lambda^{k}(i\xi')^{\kappa'}\tilde\CM^{(\ell-1)}_2(\lambda)$ belong to 
$\bM_0$,  and so by Proposition \ref{prop:2}, we have
\begin{equation}\label{5.4}\begin{aligned}
\|(\lambda, \lambda^{1/2}\bar\nabla, \bar\nabla^2) \CT(\lambda)\tilde H
\|_{W^1_q(\HS)} 
&\leq C\|\tilde H\|_{W^1_q(\HS)}, \\
\|(\lambda^{1/2}\bar\nabla, \bar\nabla^2)\CT(\lambda)\tilde H\|_{L_q(\HS)} 
&\leq C|\lambda|^{-1/2}\|\tilde H\|_{W^1_q(\HS)}.
\end{aligned}\end{equation}
\par
We next consider $\CT(\lambda)^*$, which is defined by exchanging $\CF'$ and $\CF^{-1}_{\xi'}$
in the formula of $\CT(\lambda)$.  Namely, 
\begin{align*}
\CT(\lambda)^*\tilde H = 
\int^\infty_0\CF'\Bigl[&\CM^{(0)}_{1}(\lambda)\CF^{-1}_{\xi'}[\tilde H](\xi', y_N)
Be^{-B(x_N+y_N)} \\
&+ 
\CM^{(0)}_{2}(\lambda)\CF^{-1}_{\xi'}[\tilde H](\xi', y_N)B^2M(x_N+y_N)\Bigr](x')\,dy_N.
\end{align*}
Obvisouly, for $\tilde H, \tilde G \in C^\infty_0(\HS)^{N+N^2}$, we see that 
$$((\lambda, \lambda^{1/2}\bar\nabla, \bar\nabla^2) \CT(\lambda)\tilde H, 
\tilde G)
= (\tilde H, (\lambda, \lambda^{1/2}\bar\nabla, \bar\nabla^2)
\CT^*(\lambda)\tilde G).
$$
In particular, $(\lambda, \lambda^{1/2}\bar\nabla, \bar\nabla^2)\CT(\lambda)^*
= ((\lambda, \lambda^{1/2}\bar\nabla, \bar\nabla^2 \CT(\lambda))^*$. 
\par
Employing the same argument as in the proof of 
 \eqref{5.3} and \eqref{5.4}, we have
\begin{equation}\label{5.5}\begin{aligned}
\|(\lambda, \lambda^{1/2}\bar\nabla, \bar\nabla^2) \CT(\lambda)^*
\tilde H\|_{L_{q'}(\HS)} 
& \leq C\|\tilde H\|_{L_{q'}(\HS)}, \\
\|(\lambda, \lambda^{1/2}\bar\nabla, 
\bar\nabla^2)\CT(\lambda)^*\tilde H
\|_{W^1_{q'}(\HS)} & \leq C\|\tilde H\|_{W^1_{q'}(\HS)}, \\
\|(\lambda^{1/2}\bar\nabla, \bar\nabla^2)
 \CT(\lambda)^*\tilde H\|_{L_{q'}(\HS)} &
\leq C|\lambda|^{-1/2}\|\tilde H\|_{W^1_{q'}(\HS)}.
\end{aligned}\end{equation}
In view of \eqref{5.3}, \eqref{5.4} and \eqref{5.5},  
the assertion \thetag1 of Theorem \ref{thm:5.2} 
implies that  for any $\lambda \in \Sigma_{\epsilon, \lambda_0}$, 
there hold
\begin{align*}
\|(\lambda, \lambda^{1/2}\bar\nabla, \bar\nabla^2)\CT(\lambda)
\tilde H\|_{B^\nu_{q,1}(\HS)}
&\leq C\|\tilde H\|_{B^\nu_{q,1}(\HS)} \quad (\nu\in\{s-\sigma, s, s+\sigma \}), \\
\|(\lambda, \lambda^{1/2}\bar\nabla, \bar\nabla^2)
\CT(\lambda)\tilde H\|_{B^s_{q,1}(\HS)}
&\leq C|\lambda|^{-\frac{\sigma}{2}}\|\tilde H\|_{B^{s+\sigma}_{q,1}(\HS)}.
\end{align*}
\par
We now consider $\pd_\lambda \CT(\lambda)\tilde H$.  
From \eqref{dif.1.1} it follows that 
\begin{align*}
\pd_N^\ell\pd_\lambda  \CT(\lambda)\tilde H 
=\pd_\lambda  \pd_N^\ell \CT(\lambda)\tilde H 
&= \int^\infty_0\CF^{-1}_{\xi'}\Bigl[(\pd_\lambda 
\CM^{(\ell)}_1(\lambda))\CF'[\tilde H](\xi', y_N)Be^{-B(x_N+y_N)} 
\nonumber \\
& \qquad   + B^{-2}\CM^{(\ell)}_1(\lambda)\CF'[\tilde H](\xi', y_N)
B^2\pd_\lambda(Be^{-B(x_N+y_N)} )\nonumber \\
& \qquad +(\pd_\lambda\CM^{(\ell)}_2(\lambda))\CF'[\tilde H](\xi', y_N)B^2M(x_N+y_N)
\nonumber \\
& \qquad +B^{-2}\CM_{2}^{(\ell)}(\lambda)\CF'[\tilde H](\xi', y_N)
B^2\pd_\lambda(B^2M(x_N+y_N))\Bigr](x')\,dy_N.
\end{align*}
Thus, we have
\begin{align*}
\lambda^k\pd_{x'}^{\kappa'}\pd_N^\ell\pd_\lambda \CT(\lambda)\tilde H  
& = \int^\infty_0\CF^{-1}_{\xi'}\Bigl[
\lambda^k(i\xi')^{\kappa'}(\pd_\lambda \CM^{(\ell)}_1(\lambda))
\CF'[\tilde H](\xi', y_N)Be^{-B(x_N+y_N)} 
\nonumber \\
& \qquad + \lambda^k(i\xi')^{\kappa'}B^{-2}\CM^{(\ell)}_1(\lambda)
\CF'[\tilde H](\xi', y_N)
B^2\pd_\lambda(Be^{-B(x_N+y_N)} ) \nonumber \\
& \qquad +\lambda^k(i\xi')^{\kappa'}(\pd_\lambda\CM^{(\ell)}_2(\lambda))
\CF'[\tilde H](\xi', y_N)B^2M(x_N+y_N)
\nonumber \\
& \qquad +\lambda^k(i\xi')^{\kappa'} B^{-2}\CM_{2}^{(\ell)}(\lambda)
\CF'[\tilde H](\xi', y_N)
B^2 \pd_\lambda(B^2M(x_N+y_N))\Bigr](x')\,dy_N. 
\end{align*}
If $2k+|\kappa'| + \ell \leq 4$,  the following symbols:
\begin{alignat*}2
\lambda^k(i\xi')^{\kappa'}(\pd_\lambda \CM^{(\ell)}_1(\lambda)),& &\quad 
&\lambda^k(i\xi')^{\kappa'}B^{-2}\CM^{(\ell)}_1(\lambda), \\
\lambda^k(i\xi')^{\kappa'}(\pd_\lambda\CM^{(\ell)}_2(\lambda)),  
&& \quad 
&\lambda^k(i\xi')^{\kappa'} B^{-2}\CM_{2}^{(\ell)}(\lambda)
\end{alignat*}
belong to $\bM_0$, and so by Proposition \ref{prop:2} we have
\begin{equation}\label{5.6}
\|(\lambda, \lambda^{1/2}\bar\nabla, \bar\nabla^2)
\pd_\lambda \CT(\lambda)\tilde H\|_{L_q(\HS)}
\leq C|\lambda|^{-1}\|\tilde H\|_{L_q(\HS)}. 
\end{equation}
 \par 
Moreover, from \eqref{dif.3} it follows that 
\begin{align*}
\pd_N^\ell \pd_\lambda  \CT(\lambda)\tilde H = 
\int^\infty_0\CF^{-1}_{\xi'}\Bigl[
&(\pd_\lambda \tilde\CM^{(\ell-1)}_1(\lambda))\CF'[\pd_N\tilde H](\xi', y_N)
Be^{-B(x_N+y_N)} \nonumber \\
&+B^{-2}\tilde\CM^{(\ell-1)}_1(\lambda)\CF'[\pd_N\tilde H](\xi', y_N)
B^2\pd_\lambda(Be^{-B(x_N+y_N)}) \nonumber \\
&+(\pd_\lambda \tilde\CM^{(\ell-1)}_2(\lambda))
\CF'[\pd_N\tilde H](\xi', y_N)B^2M(x_N+y_N)\Bigr](x')\,dy_N
\\
&+B^{-2}\tilde\CM^{(\ell-1)}_2(\lambda)\CF'[\pd_N\tilde H](\xi', y_N)
B^2\pd_\lambda(B^2M(x_N+y_N))\Bigr](x')\,dy_N.
\end{align*}
Thus, we have
\begin{align*}
\lambda^k\pd_{x'}^{\kappa'}\pd_N^\ell \pd_\lambda  \CT(\lambda)\tilde H = 
&\int^\infty_0\CF^{-1}_{\xi'}\Bigl[
\lambda^k(i\xi')^{\kappa'}(\pd_\lambda \tilde\CM^{(\ell-1)}_1(\lambda))
\CF'[\pd_N\tilde H](\xi', y_N)
Be^{-B(x_N+y_N)} \nonumber \\
&+\lambda^k(i\xi')^{\kappa'}B^{-2}\tilde\CM^{(\ell-1)}_1(\lambda)
\CF'[\pd_N\tilde H](\xi', y_N)
B^2\pd_\lambda(Be^{-B(x_N+y_N)}) \nonumber \\
&+\lambda^k(i\xi')^{\kappa'}
(\pd_\lambda \tilde\CM^{(\ell-1)}_2(\lambda))\CF'[\pd_N\tilde H]
(\xi', y_N)B^2M(x_N+y_N)\Bigr](x')\,dy_N
\nonumber \\
&+\lambda^k(i\xi')^{\kappa'}B^{-2}\tilde\CM^{(\ell-1)}_2(\lambda)
\CF'[\pd_N\tilde H](\xi', y_N)
B^2\pd_\lambda(B^2M(x_N+y_N))\Bigr](x')\,dy_N. 
\end{align*}
If $2k+|\kappa'| + \ell \leq 5$,  the following symbols:
\begin{alignat*}2
\lambda^k(i\xi')^{\kappa'}(\pd_\lambda \tilde\CM^{(\ell-1)}_1(\lambda)), &
&\quad&\lambda^k(i\xi')^{\kappa'}B^{-2}\tilde\CM^{(\ell-1)}_1(\lambda), \\
\lambda^k(i\xi')^{\kappa'}
(\pd_\lambda \tilde\CM^{(\ell-1)}_2(\lambda)),&&\quad&
\lambda^k(i\xi')^{\kappa'}B^{-2}\tilde\CM^{(\ell-1)}_2(\lambda))
\end{alignat*}
 all belong to $\bM_0$, and  so by Proposition \ref{prop:2} we have
\begin{equation}\label{5.7}\begin{aligned}
\|(\lambda, \lambda^{1/2}\bar\nabla, \bar\nabla^2)\pd_\lambda 
\CT(\lambda)\tilde H \|_{W^1_q(\HS)}
&\leq C|\lambda|^{-1}\|\tilde H\|_{W^1_q(\HS)}, \\
\|(\lambda, \lambda^{1/2}\bar\nabla, \bar\nabla^2)\pd_\lambda 
\CT(\lambda)\tilde H \|_{W^1_q(\HS)}
&\leq C|\lambda|^{-1/2}\|\tilde H\|_{L_q(\HS)}.
\end{aligned}\end{equation}
\par
Since $\pd_\lambda \CT_1(\lambda)^*$ is obtained by exchanging $\CF^{-1}_{\xi'}$
and $\CF'$, that is 
\begin{align*}
\pd_\lambda\CT(\lambda)^*\tilde H &= 
\int^\infty_0\CF'\Bigl[(\pd_\lambda\CM^{(0)}_{1}(\lambda))
\CF^{-1}_{\xi'}[\tilde H](\xi', y_N)Be^{-B(x_N+y_N)} \\
&\qquad + B^{-2}\CM^{(0)}_{1}(\lambda)\CF^{-1}_{\xi'}[\tilde H]
(\xi', y_N)B^2\pd_\lambda(Be^{-B(x_N+y_N)}) \\
&\qquad + 
(\pd_\lambda\CM^{(0)}_{2}(\lambda))\CF^{-1}_{\xi'}[\tilde H](\xi', y_N)B^2M(x_N+y_N) \\
&\qquad + 
B^{-2}\CM^{(0)}_{2}(\lambda)\CF^{-1}_{\xi'}[\tilde H](\xi', y_N)B^2
\pd_\lambda(B^2M(x_N+y_N))\Bigr](x')\,dy_N. 
\end{align*}
Employing the same argument, we see that 
\begin{equation}\label{5.8}\begin{aligned}
\|(\lambda, \lambda^{1/2}\bar\nabla, \bar\nabla^2)\pd_\lambda \CT(\lambda)^*\tilde H\|_{L_{q'}(\HS)}
\leq C|\lambda|^{-1}\|\tilde H\|_{L_{q'}(\HS)}, \\
\|(\lambda, \lambda^{1/2}\bar\nabla, \bar\nabla^2)\pd_\lambda \CT(\lambda)^*\tilde H\|_{W^1_{q'}(\HS)}
\leq C|\lambda|^{-1}\|\tilde H\|_{W^1_{q'}(\HS)},  \\
\|(\lambda, \lambda^{1/2}\bar\nabla, \bar\nabla^2)\pd_\lambda \CT(\lambda)^*\tilde H\|_{W^1_{q'}(\HS)}
\leq C|\lambda|^{-1/2}\|\tilde H\|_{L_{q'}(\HS)}.
\end{aligned}\end{equation}
Since $\tilde H \in C^\infty_0(\HS)^{N+N^2}$,  we see that
 $((\lambda, \lambda^{1/2}\bar\nabla. \bar\nabla^2)(\pd_\lambda\CT(\lambda))^{*}
=(\lambda, \lambda^{1/2}\bar\nabla, \bar\nabla^2)(\pd_\lambda\CT(\lambda))^*$. 
 Thus, in view of \eqref{5.6}, \eqref{5.7} and \eqref{5.8}, by Theorem \ref{thm:5.2}, 
we see that for any $\lambda \in \Sigma_{\epsilon, \lambda_0}$
there holds 
\begin{align*}
\|(\lambda, \lambda^{1/2}\bar\nabla, \bar\nabla^2)\pd_\lambda \CT(\lambda)
\tilde H\|_{B^\nu_{q,1}(\HS)} & \leq C|\lambda|^{-1}\|\tilde H\|_{B^\nu_{q,1}(\HS)}
\quad (\nu \in \{ s-\sigma, s, s+\sigma\} ), \\
\|(\lambda, \lambda^{1/2}\bar\nabla, \bar\nabla^2)\pd_\lambda \CT(\lambda)
\tilde H\|_{B^s_{q,1}(\HS)} & \leq C|\lambda|^{-(1-\frac{\sigma}{2})}
\|\tilde H\|_{B^{s-\sigma}_{q,1}(\HS)}.
\end{align*}
\par
Finally, we shall see that for any $\lambda \in \Sigma_{\epsilon, \lambda_0}$
there holds 
\begin{equation}\label{5.1.3}
\|(1, \lambda^{-1/2}\bar\nabla)\CT(\lambda)\tilde H\|_{B^s_{q,1}(\HS)}
\leq C|\lambda|^{-(1-\frac{\sigma}{2})}\|\tilde H\|_{B^{s-\sigma}_{q,1}(\HS)}.
\end{equation}
In fact, by \eqref{5.3}, we have
\begin{align*}
\|(1, \lambda^{-1/2}\bar\nabla)\CT(\lambda)\tilde H\|_{L_q(\HS)} 
& \leq C|\lambda|^{-1}\|\tilde H\|_{L_q(\HS)}, \\
\|(1, \lambda^{-1/2}\bar\nabla)\CT(\lambda)\tilde H\|_{W^1_q(\HS)} 
& \leq C|\lambda|^{-1/2}\|\tilde H\|_{L_q(\HS)}.
\end{align*}
And, by \eqref{5.4}, we have
$$\|(1, \lambda^{-1/2}\bar\nabla)\CT(\lambda)\tilde H\|_{W^1_q(\HS)}
 \leq C|\lambda|^{-1}\|\tilde H\|_{W^1_q(\HS)}.
$$
Moreover, by \eqref{5.5}, we have
\begin{align*}
\|(1, \lambda^{-1/2}\bar\nabla)\CT(\lambda)^*\tilde H\|_{L_{q'}(\HS)}
& \leq C|\lambda|^{-1}\|\tilde H\|_{L_{q'}(\HS)}, \\
\|(1, \lambda^{-1/2}\bar\nabla)\CT(\lambda)^*\tilde H\|_{W^1_{q'}(\HS)}& \leq C|\lambda|^{-1/2}
\|\tilde H\|_{L_{q'}(\HS)}, \\
\|(1, \lambda^{-1/2}\bar\nabla)\CT(\lambda)^*\tilde H\|_{W^1_{q'}(\HS)}
& \leq C|\lambda|^{-1}\|\tilde H\|_{W^1_{q'}(\HS)}.
\end{align*}
Therefore, by Theorem \ref{thm:5.2}, we see that 
\eqref{5.1.3} holds. 
This completes the proof of Theorem \ref{thm:3.1}. 
\end{proof}

Finally, we consider the full equations:
\begin{equation}\label{full:Lame}
\begin{aligned}
\lambda \bu - \DV(\alpha\BD(\bu) + (\beta-\alpha)\dv\bu\BI)
&= \bg&\quad&\text{in $\HS$}, \\
(\alpha \BD(\bu) + (\beta-\alpha)\dv\bu)\bn\BI &= \bh
&\quad&\text{on $\pd\HS$}.
\end{aligned}\end{equation}
Combining Theorems \ref{thm:L.1} and \ref{thm:3.1}, we have
the following corollary immediately.
\begin{cor}\label{col:half} Let $1 < q < \infty$, $0 < \epsilon < \pi/2$ and 
$\lambda_0>0$.  Let $\sigma>0$ be a small constant such that
$-1+1/q < s-\sigma < s+\sigma < 1/q$.  Then, there exists an operator
$\CS_{\HS}(\lambda)$ with
$$\CS_{\HS}(\lambda) \in {\rm Hol}\,(\Sigma_{\epsilon, \lambda_0}, 
\CL(B^s_{q,1}(\HS)^{m(N)}, B^{s+2}_{q,1}(\HS)^N))$$
 haviing 
$(s, \sigma, q)$ properties such that 
for any $(\bg, \bh) \in \CH^s_{q,1}(\HS)$, $\bu=\CS_{\HS}(\lambda)\CO_\lambda
(\bg, \bh)$ is a unique solution of equations \eqref{full:Lame}.
Here and in the following, an operation $\CO_\lambda$ is defined by
$$\CO_\lambda(\bg, \bh) = (\bg, \lambda^{1/2}\bh, \nabla \bh).
$$
\end{cor}


\subsection{Spectral analysis of  Lam\'e equations 
 in a bent half space} \label{subsec.4.4}

Let $x_0 \in \pd\Omega$. As was seen in  \cite[Appendix]{ES13}
or in \cite[Subsec. 3.2.1]{S20}, there exist a constant $d>0$, a 
diffeomorphism of  $C^3$ class $\Phi: \BR^N \to \BR^N$, $x \mapsto y=\Phi(x)$ 
and its inverse map $\Phi^{-1}: \BR^N \to \BR^N$, $y \mapsto x=\Phi^{-1}(y)$  
such that $B_d(x_0) \cap \Omega \subset \Phi(\Omega)$, 
and $B_d(x_0) \cap \pd\Omega \subset \Phi(\pd\Omega)$ and 
\begin{equation*}
\nabla \Phi = \CA + \CB(x), \quad \nabla \Phi^{-1}(y) = \CA_- + \CB_-(y)
\end{equation*}
where  $\CA$ and $\CA_-$ are $N\times N$ orthogonal matrices of constant coefficients
such that $\CA\CA_- = \CA_-\CA = I$ and $\CB(x)$ and $\CB_-(y)$ are $N\times N$
matrices of $C^2$ functions. 
Here and in the following, 
we write $B_d(x_0) = \{y \in \BR^N \mid |y-x_0| < d\}$ and 
$B_d = \{x \in \BR^N \mid |x| < 1\}$. 
Moreover, we may assume that for any small constant  $M_1 > 0$
we can choose  $0 < d < 1$ such that 
\begin{equation}\label{coef.11}
\|(\CB, \CB_-)\|_{L_\infty(\BR^N)}  \leq M_1. 
\end{equation}
Furthermore, we may assume that there exist constants 
$D$ and $M_2$ such that  $D$ is independent of 
$M_1$, but $M_2$ depends on $M_1$, and 
\begin{equation}\label{coef.12}\begin{aligned}
\|\nabla(\CB, \CB_-)\|_{L_\infty(\BR^N)}&\leq D, \\
\|\nabla^2(\CB, \CB_-)\|_{L_\infty(\BR^N)}& \leq M_2.
\end{aligned}\end{equation}
We may assume that $M_1 < 1 \leq D \leq M_2$. 
\par
Let 
\begin{equation*}
\Omega_+ = \Phi(\HS), \quad \Gamma_+=\Phi(\BR^N_0). 
\end{equation*}
Let $\bn_+$ denote the unit outer normal to $\Gamma_+$.
In this section, first we consider 
the Lam\'e equations 
\begin{equation}\label{resol.1}\begin{aligned}
\lambda \bu - \DV(\alpha\BD(\bu) + (\beta-\alpha)\dv \bu\BI)  &= \tilde\bg&\quad
&\text{in $\Omega_+$}, \\ 
(\alpha\BD(\bu) + (\beta-\alpha)\dv\bu\BI)\bn_+&= \tilde\bh &\quad
&\text{on $\pd\Omega_+$}.
\end{aligned}\end{equation}
%
%

We shall show the following theorem.
\begin{thm}\label{thm:bent} Let $x_0 \in \pd\Omega$. 
Let $\Phi$ and $\Phi^{-1}$ be a $C^3$ diffeomorphism on $\BR^N$ and 
its inverse given above. 
Let $1 < q < \infty$,  $-1+1/q < s < 1/q$, and $0< \epsilon<\pi/2$.
 Let $\sigma$ be a small 
positive number such that $-1+1/q < s-\sigma < s+\sigma < 1/q$.  
Let $\nu \in \{s-\sigma, s, s+\sigma\}$. 
Then,
there exist a small constant $d>0$, a large constant $\lambda_1 \geq \max(1, \lambda_0)$ 
and an operator $\CS_{\Omega_+}(\lambda)
 \in {\rm Hol}\,(\Sigma_{\epsilon, \lambda_1},
\CL(B^\nu_{q,1}(\Omega_+)^{m(N)}, B^{\nu+2}_{q,1}(\Omega_+)^N))$ 
having  $(s,\sigma,q)$ properties 
such that for any $(\tilde \bg, \tilde\bh) \in \CH^{\nu}_{q,1}(\Omega_+)$,
$\bu = \CS_{\Omega_+}(\lambda)\CO_\lambda(\bg, \bh)$ 
is a unique solution of equations \eqref{resol.1}. 
\end{thm}
\begin{proof}
First, we shall reduce problem \eqref{resol.1}
to that in the half-space $\HS$.
Let $a_{jk}$ and $b_{jk}(x)$ be the $(j,k)$th components of $\CA_-$ and $\CB_-(\Phi(x))$, 
and then we have
\begin{equation}\label{change.1}
\frac{\pd}{\pd y_j} = \sum_{k=1}^N(a_{kj} + b_{kj}(x))\frac{\pd}{\pd x_k}\quad(j=1, \ldots, N).
\end{equation}
Since $y= \Phi(x)$ denotes the point of $\Gamma_+$ for $x \in \BR^N_0$, 
$dy = \nabla\Phi(x) dx$.  Thus, $\bn_+ = \nabla \Phi(x) \bn_0/|\nabla \Phi(x)\bn_0|$.
Recalling that $\bn_0=(0, \ldots, 0, -1)^t$ and denoting the $(j,k)$th component of
$\nabla \Phi=\CA+\CB(x)$ by $a_{kj} + \tilde b _{kj}(x)$, we have
\begin{equation*}
\bn_+ =-(a_{N1}+\tilde b_{N1}(x), \ldots, a_{NN} +\tilde b_{NN}(x))/|\bn_+|
\end{equation*}
with $|\bn_+| = (\sum_{j=1}^N(a_{Nj}+\tilde b_{Nj})^2)^{1/2}$.
Notice that 
\begin{equation}\label{change.2}
\sum_{j=1}^N a_{jk}a_{j\ell} =\sum_{j=1}^N  a_{kj}a_{\ell j} = \delta_{k\ell}.
\end{equation}
In particular, we see that
\begin{equation}\label{est.d}
|d^2-1| \leq \sum_{j=1}^N( 2|\tilde b_{Nj}||a_{Nj}| + |\tilde b_{Nj}|^2) \leq CM_1.
\end{equation}
Let $\tilde \bv(x) = \bv(y)$.  We write $\tilde \bv(x) = (\tilde v_1(x), \ldots, \tilde v_N(x))^\top$ 
and $\bv(y)=(v_1(y), \ldots, v_N(y))^\top$, where $A^\top$ denotes the transposed $A$ 
for any vector or matrix $A$.  
By \eqref{change.1} we have
\begin{equation}\label{dv.1}
\dv_y \bv(y) = \sum_{\ell, m=1}^N (a_{m\ell} +b_{m\ell}(x))\frac{\pd \tilde v_\ell}{\pd x_m}.
\end{equation}
Moreover, we set $\tilde v_\ell = \sum_{k=1}^N a_{k\ell}w_k$, 
and $\bw = (w_1, \ldots, w_N)^\top$. 
From \eqref{change.1}, \eqref{change.2} and \eqref{dv.1}, we have
\allowdisplaybreaks
\begin{align*}
\tilde g_j& = \lambda v_j - \sum_{k=1}^N\frac{\pd}{\pd y_k}(\alpha D_{jk}(\bv)
+(\beta-\alpha)\delta_{jk}\dv\bv) \\
& = \lambda \sum_{n=1}^Na_{nj}w_n 
- \alpha \sum_{\ell, n=1}^Na_{\ell j}\frac{\pd}{\pd x_n}D_{\ell n}(\bw)
-(\beta-\alpha)\sum_{n=1}^Na_{nj}\frac{\pd}{\pd x_n}\dv\bw
-R_j^1\bw.
\end{align*}
Here, we have set
\begin{align*}
R^1_j\bw &= \alpha \sum_{p, k=1}^Na_{pk}\frac{\pd}{\pd x_p} 
\Bigl(\sum_{\ell, n=1}^N(b_{\ell j}a_{nk} + b_{\ell k}a_{nj})\frac{\pd w_n}{\pd x_\ell} \Bigr)\\
& \qquad + \alpha \sum_{p, k=1}^Nb_{pk}\frac{\pd}{\pd x_p}
\Bigl( \sum_{\ell, n=1}^N(a_{\ell j}a_{nk}+a_{\ell k}a_{nj}
+b_{\ell j}a_{nk} + b_{\ell k}a_{nj})\frac{\pd w_n}{\pd x_\ell} \Bigr)\\
& \qquad +(\beta-\alpha)\Bigl(\sum_{p=1}^N a_{pj}\frac{\pd}{\pd x_p}
(\sum_{\ell, m, n=1}^N
a_{\ell m}b_{nm}\frac{\pd w_\ell}{\pd x_n})
+ \sum_{p=1}^Nb_{pj}\frac{\pd}{\pd x_p}(\dv\bw + \sum_{\ell, m, n=1}^N
a_{\ell m}b_{nm}\frac{\pd w_\ell}{\pd x_n})\Bigr).
\end{align*}
Using \eqref{change.2}, we have
\begin{equation}\label{reduceeq.1}
\sum_{j=1}^N a_{sj}\tilde g_j = \lambda w_s -\alpha \sum_{n=1}^N \frac{\pd}{\pd x_n}
D_{sn}(\bw) - (\beta-\alpha)\sum_{n=1}^N\delta_{sn}\frac{\pd}{\pd x_n}\dv\bw
-\sum_{j=1}^Na_{sj}R_j^1\bw.
\end{equation}
\par
We next consider the boundary conditions. 
\begin{align*}
\tilde h_j & = \sum_{k=1}^N (\alpha D_{jk}(\bv) + (\beta-\alpha)\delta_{jk}\dv\bv)n_k
= -(\sum_{\ell=1}^Na_{\ell j}\alpha D_{N\ell}(\bw) + (\beta-\alpha)a_{Nj}
\dv\bw + R^2_j \bw)
\end{align*}
where we have set
\begin{align*}
R^2_j\bw & = 
(d^{-1}-1)(\alpha \sum_{k,\ell, n=1}^N(a_{\ell j}a_{nk}+a_{\ell k}a_{nj})
\frac{\pd w_n}{\pd x_\ell}a_{Nk}+ (\beta-\alpha)\delta_{jk}a_{Nk}\dv \bw ) \\
& \qquad +d^{-1}\alpha \sum_{k, \ell, n=1}^N(a_{\ell j}a_{nk} + a_{\ell k}a_{nj})
\frac{\pd w_n}{\pd x_\ell} \tilde b_{Nk} 
+d^{-1}\alpha  \sum_{k,\ell, n=1}^N(b_{\ell j}a_{nk} + b_{\ell k}a_{nj})
\frac{\pd w_n}{\pd x_\ell} (a_{Nk}+\tilde b_{Nk}) \\
& \qquad +d^{-1}(\beta-\alpha)(\dv\bw\tilde b_{Nj} 
+ \sum_{\ell,m,n=1}^N a_{n\ell}b_{m\ell}\frac{\pd w_n}{\pd x_\ell}(a_{Nj}+\tilde b_{Nj})). 
\end{align*}
Using \eqref{change.2}, we have
\begin{equation}\label{reduceeq.2}
\sum_{j=1}^Na_{sj} \tilde{h}_j = -(\alpha D_{Ns}(\bw) + (\beta-\alpha)\delta_{Ns}\dv\bw
+\sum_{j=1}^Na_{sj}R^2_j \bw).
\end{equation}
Set 
\begin{equation*}
\begin{aligned}
\bg &=(\sum_{j=1}^Na_{1j}\tilde g_j(\Phi(x)), \ldots, \sum_{j=1}^Na_{Nj}\tilde g_j(\Phi(x))), \quad 
\bh =(\sum_{j=1}^Na_{1j}\tilde h_j(\Phi(x)), \ldots, \sum_{j=1}^Na_{Nj}\tilde h_j(\Phi(x))), 
\end{aligned}\end{equation*}
and define remainder terms $\CR^k$ by 
$$\CR^k\bw  = \Bigl(\sum_{j=1}^N a_{1j}R^k_j\bw, \ldots, \sum_{j=1}^Na_{Nj}R^k_j\bw \Bigr)  
\qquad (k=1,2). $$
Then, from \eqref{reduceeq.1} and \eqref{reduceeq.2} we have
\begin{equation}\label{reduce.3}\begin{aligned}
\lambda \bw - \DV(\alpha \BD(\bw) + (\beta-\alpha)\dv\bw\BI)
- \CR^1\bw = \bg&&
\quad&\text{in $\HS$}, \\
(\alpha\BD(\bw) + (\beta-\alpha)\dv\bw \BI)\bn_0 -\CR^2\bw
 = \bh&&\quad&\text{on $\pd\HS$}.
\end{aligned}\end{equation}
\par
Let $\lambda_0 > 0$ and 
let $\CS_{\HS}(\lambda)$ be the solution operator given in Corollary \ref{col:half}.
Recall that $\bw = \CS_{\HS}(\lambda)\CO_\lambda(\bg, \bh)$
is a unique solution of equations \eqref{free.1}. Inserting this formula into equations
\eqref{reduce.3}, we have
\begin{equation}\label{reduce.4}\begin{aligned}
\lambda \bw - \DV(\alpha \BD(\bw) + (\beta-\alpha)\dv\bw\BI)
- \CR^1\bw &=  \bg- \CR^1\bw&
\quad&\text{in $\HS$}, \\
(\alpha\BD(\bw) + (\beta-\alpha)\dv\bw \BI)\bn_0 -\CR^2\bw
 &= \bh-\CR^2\bw&\quad&\text{on $\pd\HS$}.
\end{aligned}\end{equation}
Define an operator $\CU$ by 
\begin{equation*}
\CU H= (\CR^1\CS_{\HS}(\lambda)H, \CR^2\CS_{\HS}(\lambda)H).
\end{equation*}
for $H \in B^{\nu}_{q,1}(\HS)^{m(N)}$. 
Obvisouly, 
\begin{equation}\label{cont.0}
(\bg-\CR^1\bw, \bh-\CR^2\bw) = (\bI- \CU\CO_\lambda)(\bg, \bh).
\end{equation}

To estimate  $\CU$, we prepare some lemmas concerning the estimate of
products in the Besov space.
\begin{lem}\label{lem:APH} Let $1 < q < \infty$ 
and $-1+1/q < \nu < 1/q$.
Let $p_2$ be an exponent such that $N < p_2 < \min(q, q')N$.  Then, we have
\begin{equation}\label{p.prod1}
\|uv\|_{B^\nu_{q,1}(\HS)} \leq C\|u\|_{B^\nu_{q,1}(\HS)}
\|v\|_{B^{N/p_2}_{p_2,\infty}(\HS)  \cap L_\infty(\HS)}.
\end{equation}
\end{lem}
\begin{remark} Since $B^{N/p_2}_{p_2,1}$ is continously imbedded into 
$B^{N/p_2}_{p_2,\infty}(\HS) \cap L_\infty(\HS)$, and therefore, 
from \eqref{p.prod1} it follows that 
\begin{equation}\label{p.prod1*}
\|uv\|_{B^\nu_{q,1}(\HS)} \leq C\|u\|_{B^\nu_{q,1}(\HS)}
\|v\|_{B^{N/p_2}_{p_2,1}(\HS)}.
\end{equation}
\end{remark}
\begin{proof}
By using an extension map from $\HS$ into $\BR^N$, 
it is sufficient to prove the lemma
in the case where the domain is $\BR^N$ instead of $\HS$.  
Below, we omit $\BR^N$. 
We shall use the Abidi-Paicu theory \cite[Cor.2.5]{AP07} 
and the Haspot theory \cite[Prop. 2.3]{H11}.
According to the Abidi-Paicu-Haspot theory, we have
$$\|uv\|_{B^{s_1+s_2-N(\frac{1}{p_1}+\frac{1}{p_2}-\frac{1}{q})}_{q,1} }
\leq C\|u\|_{B^{s_1}_{q,1}}
\|v\|_{B^{s_2}_{p_2, r}\cap L_\infty}$$
provided that $1/q \leq 1/p_1+1/\lambda_1 \leq 1$, 
$1/q \leq 1/p_2 + 1/\lambda_2 \leq 1$, 
$1/q \leq 1/p_1 + 1/p_2$, $p_1 \leq \lambda_2$, $p_2 \leq \lambda_1$, 
$s_1+s_2+N\inf(0, 1-1/p_1-1/p_2)>0$, $s_1 + N/\lambda_2 < N/p_1$
 and $s_2 + N/\lambda_1 \leq N/p_2$. 
We choose $p_1=q$, $s_1=\nu$ and $s_2 = N(\frac{1}{p_1}+\frac{1}{p_2}
-\frac{1}{q})= N/p_2$. 
In particular, $s_1+s_2-N(\frac{1}{p_1}+\frac{1}{p_2}-\frac{1}{q})=\nu$. 
Let $\lambda_1=\infty$, and then
$1/q \leq 1/q + 0 \leq 1$, $p_2 \leq \lambda_1$.  We choose $\lambda_2$ 
in such a way that
$1/\lambda_2 = 1/q -1/p_2$ when $1/q \geq 1/p_2$ and $\lambda_2=\infty$ when
$1/q < 1/p_2$. In this case, $s_1 + N/\lambda_2 = \nu < 1/q < N/q$ when $1/q < 1/p_2$.
When $1/q \geq 1/p_2$,
$s_1 + N/\lambda_2 = \nu + N(1/q-1/p_2) < N/q$ , namely we choose $p_2$ such that  
$\nu- N/p_2 <0$. Since $\nu < 1/q$, we choose $p_2$ such that $1/q \leq N/p_2$, 
that is $p_2 \leq qN$.  Thus, so far we choose $p_2$ in such a way that $N < p_2 < qN$.
Since $\lambda_1=\infty$, the condition $p_2 \leq \lambda_1$ is 
satisfied. When $1/q \geq 1/p_2$, 
$\lambda_2^{-1} = 1/q -1/p_2 < 1/q$, and so $q < \lambda_2$.
When $1/q < 1/p_2$, $\lambda_2=\infty$, and so $q \leq \lambda_2$.
When $1-1/q-1/p_2 \geq 0$, that is $p_2 \geq q'$, 
$s_1+s_2+N\inf(0, 1/p_1-1/p_2) = \nu+N/p_2 > 0$.
Since $\nu > -1+1/q = -1/q'$, we have $-N/p_2 < -1/q'$ provided that $p_2 \leq Nq'$. 
When $1-1/q-1/p_2 < 0$, that is  $p_2 < q'$, 
$s_1+s_2+N\inf(0, 1-1/p_1-1/p_2) = \nu+N/p_2+N/q'-N/p_2 = \nu+N/q' > 0$
because $s > -1/q'$.  Summing up, if $N < p_2 < \min(q,q')N$, 
then the Abidi-Paicu-Haspot 
conditions are all satisfied.  
Thus, we have \eqref{p.prod1}.  This completes the proof of 
Lemma \ref{lem:APH}.
\end{proof}
\begin{lem}\label{lem:dif1} Let $1 < q < \infty$  and
 $-1+1/q < \nu < 1/q$.
Then, for $f \in B^\nu_{q,1}(\HS)$ and $g \in W^1_\infty(\HS)$, there holds 
\begin{equation}\label{dif.1}
\|fg\|_{B^\nu_{q,1}(\HS)} \leq C_s\|f\|_{B^\nu_{q,1}(\HS)}\|g
\|_{L_\infty(\HS)}^{1-|\nu|}\|g\|_{W^1_\infty(\HS)}^{|\nu|}.
\end{equation}
provided that $s\not=0$ and 
\begin{equation*}
\|fg\|_{B^0_{q,1}(\HS)} \leq C_\epsilon\|f\|_{B^0_{q,1}(\HS)}\|g\|_{L_\infty(\HS)}^{1-\epsilon}
\|g\|_{W^1_\infty(\HS)}^{\epsilon}.
\end{equation*}
with any small $\epsilon > 0$. Here, 
 $C_s$ and $C_\epsilon$ denote constants being independent of $f$ and $g$.
\end{lem}
\begin{proof} First, we consider the case where $0 < \nu < 1/q$.  
Since $C^\infty_0(\HS)$ is 
dense in $B^\nu_{q,1}(\HS)$, we may assume that $f \in C^\infty_0(\HS)$.  We know that 
\begin{equation*}
(L_q(\HS), W^1_q(\HS))_{\nu, r} = B^\nu_{q,1}(\HS).
\end{equation*}
Here, $(\cdot, \cdot)_{\nu, r}$ denotes the real interpolation functor. 
We see easily that 
$$\|fg\|_{L_q(\HS)} \leq \|f\|_{L_q(\HS)} \|g\|_{L_\infty(\HS)}, \quad
\|fg\|_{W^1_q(\HS)} \leq \|f\|_{W^1_q(\HS)} \|g\|_{W^1_\infty(\HS)}.$$	
Since $(\cdot, \cdot)_{\nu,r}$ is an exact interpolation functor of exponent $\nu$ (cf. \cite[p.41,
in the proof of Theorem 3.1.2]{BL}), we have
$$\|fg\|_{B^\nu_{q,1}(\BR^N)} \leq C\|g\|_{L_\infty(\HS)}^{1-\nu}
\|g\|_{W^1_\infty(\HS)}^\nu\|f\|_{B^\nu_{q,1}(\HS)}.
$$
This shows \eqref{dif.1} for $0 < \nu < 1/q$. \par
Next we consider the case where $-1+1/q < \nu < 0$.  For any $\varphi \in C^\infty_0(\BR^N)$, 
we have
\begin{align*}
|(fg, \varphi)_{\BR^N}| &= |(f, g\varphi)_{\BR^N}| \leq \|f\|_{B^\nu_{q,1}(\BR^N)}
\|g\varphi\|_{B^{-\nu}_{q',r'}(\BR^N)} 
\leq C\|g\|_{L_\infty(\HS)}^{1-|\nu|}
\|g\|_{W^1_\infty(\HS)}^{|\nu|} \|f\|_{B^\nu_{q,1}(\HS)}
\|\varphi\|_{B^{-\nu}_{q', r'}(\BR^N)}.
\end{align*}
Since $C^\infty_0(\BR^N)$ is dense in $B^{-\nu}_{q', r'}(\HS)$, we have
$$\|fg\|_{B^\nu_{q,1}(\HS)} \leq C\|g\|_{L_\infty(\HS)}^{1-|\nu|}
\|g\|_{W^1_\infty(\HS)}^{|\nu|} \|f\|_{B^s_{q,1}(\HS)}.$$
Since $B^0_{q,1}(\BR^N) = (B^\epsilon_{q,1}(\BR^N), B^{-\epsilon}_{q,1}(\BR^N))_{1/2, r}$ 
for any $\epsilon > 0$,  we have
$$\|fg\|_{B^0_{q,1}(\HS)} \leq C_\epsilon 
\|g\|_{L_\infty(\HS)}^{1-\epsilon}
\|g\|_{W^1_\infty(\HS)}^{\epsilon} \|f\|_{B^0_{q,1}(\HS)}.$$
This completes the proof of Lemma \ref{lem:dif1}. \end{proof}

{\bf Continuation of the proof of Theorem \ref{thm:bent}.}  
Since $\CS_{\HS}$ has $(s, \sigma, q)$ properties, 
by \eqref{coef.11}, \eqref{coef.12}, \eqref{est.d} and Lemma \ref{lem:dif1}, 
we see that for some $\omega \in (0, 1)$ 
\begin{equation*}
\begin{aligned}
\|\CO_\lambda \CU H\|_{B^\nu_{q,1}(\HS)}
& \leq C((M_1)^{1-\omega}D^\omega+D^{1-\omega}M_2^\omega
|\lambda|^{-1/2})\|H\|_{B^\nu_{q,1}(\HS)}.
\end{aligned}\end{equation*}
Choosing $M_1$ so small that $(CM_1)^{1-\omega}D^\omega \leq 1/4$ and 
choosing $\lambda_1 \geq \lambda_0$
so large that $CD^{1-\omega}M_2^{\omega}\lambda_1^{-1/2} \leq 1/4$, we have
\begin{equation}\label{small.3.1}
\|\CO_\lambda\CU H\|_{B^\nu_{q,1}(\HS)}
\leq (1/2)\|H\|_{B^\nu_{q,1}(\HS)}.
\end{equation}
Here, we choose $d>0$ so small that $(CM_1)^{1-\omega}D^\omega \leq 1/4$
and fix such $d$.  After this procedure, $M_2$ is fixed, and so  we can choose $\lambda_1$
so large according to these fixed $d$ and  $M_2$. \par 

Let us define $\CR_\infty(\lambda)$ by setting 
$$\CR_\infty(\lambda) = \sum_{\ell=0}^\infty(\CU\CO_\lambda)^\ell. $$
By \eqref{small.3.1}, 
we see that  for $(\bg, \bh) \in \CH^\nu_{q,1}(\HS)$
and $\lambda \in \Sigma_{\epsilon, \lambda_1}$, 
\begin{align*}
\|(\CU\CO_\lambda)^{\ell}(\bg, \bh)\|_{\CH^\nu_{q,1}(\HS)}
&\leq \max(1, \lambda_1^{-1/2})\|\CO_\lambda
(\CU\CO_\lambda)^{\ell}(\bg, \bh)\|_{B^\nu_{q,1}(\HS)} \\
&= \max(1, \lambda_1^{-1/2})\|(\CO_\lambda \CU)^\ell 
\CO_\lambda(\bg, \bh)\|_{B^\nu_{q,1}(\HS)}\\
& \leq \max(1, \lambda_1^{-1/2})(1/2)^\ell \|\CO_\lambda(\bg, \bh)\|_{B^\nu_{q,1}(\HS)}
\\
& \leq \max(1, \lambda_1^{-1/2})\max(1, |\lambda|^{1/2})(1/2)^\ell 
\|(\bg, \bh)\|_{\CH^\nu_{q,1}(\HS)}.
\end{align*}
Thus, $\CR_\infty(\lambda) \in \CL(\CH^\nu_{q,1}(\HS))$.
If we define $\bv \in B^{\nu+2}_{q,1}(\HS)$ by 
\begin{align*}
\bv = \CS_{\HS}(\lambda)\CO_\lambda\CR_\infty(\lambda)(\bg, \bh)
\end{align*}
then, in veiw of \eqref{cont.0}, 
\begin{align*}
\CR_\infty(\lambda)(\bg, \bh) -\CU\CO_\lambda \CR_\infty(\lambda)(\bg, \bh)
= (\bg, \bh).
\end{align*}
Thus, from \eqref{reduce.4}, $\bv$ satisfies equations: 
\begin{equation}\label{reduce.5}\begin{aligned}
\lambda \bv - \DV(\alpha \BD(\bv) + (\beta-\alpha)\dv\bv\BI)
- \CR^1\bv &=  \bg&
\quad&\text{in $\HS$}, \\
(\alpha\BD(\bv) + (\beta-\alpha)\dv\bv \BI)\bn_0 -\CR^2\bv
 &= \bh&\quad&\text{on $\pd\HS$}.
\end{aligned}\end{equation}
Moreover, we observe that 
$$\CO_\lambda \CR_\infty(\lambda) 
= \sum_{\ell=0}^\infty \CO_\lambda (\CU\CO_\lambda)^\ell
= (\sum_{\ell=0}^\infty (\CO_\lambda \CU))^\ell )\CO_\lambda.
$$
Let $\CQ_\infty(\lambda)$ be an operator defined by 
$$\CQ_\infty(\lambda)H = \sum_{\ell=0}^\infty (\CO_\lambda \CU)^\ell H, $$
and then by \eqref{small.3.1} we see that $\CQ_\infty(\lambda)
\in {\rm Hol}\,(\Sigma_{\epsilon, \lambda_1}, \CL(B^\nu_{q,1}(\HS)^{m(N)}))$
and its operator norm does not exceed $2$ for any $\lambda \in \Sigma_{\epsilon, \lambda_1}$.
We define an operator $\CT_{\Omega_+}(\lambda)$ by 
$$\CT_{\Omega_+}(\lambda)H= \CS_{\HS}(\lambda)\CQ_\infty(\lambda)H.$$
Since $\CT_{\Omega_+}
(\lambda)\CO_\lambda = \CT(\lambda)\CO_\lambda\CR_\infty(\lambda)$, 
$\bv = \CT_{\Omega_+}(\lambda)\CO_\lambda(\bg, \bh)$ is a solution of equations \eqref{reduce.5}.
Moreover, $\CS_{\HS}(\lambda)$ has $(s, \sigma, q)$ properties, 
we see that 
\begin{equation*}
\begin{aligned}
\|(\lambda, \lambda^{1/2}\bar\nabla, \bar\nabla^2)\CT_{\Omega_+}(\lambda)H
\|_{B^\nu_{q,1}(\HS)} & \leq C\|H\|_{B^\nu_{q,1}(\HS)}
&\quad&\text{for $H \in B^{\nu}_{q,1}(\HS)^{m(N)}$}, \\
\|(\lambda^{1/2}\bar\nabla, \bar\nabla^2)\CT_{\Omega_+}(\lambda)H
\|_{B^s_{q,1}(\HS)} & \leq C|\lambda|^{-\frac{\sigma}{2}}\|H\|_{B^{s+\sigma}_{q,1}(\HS)}
&\quad&\text{for $H \in B^{s+\sigma}_{q,1}(\HS)^{m(N)}$}, \\
\|(1, \lambda^{-1/2}\bar\nabla)\CT_{\Omega_+}H\|_{B^s_{q,1}(\HS)}
& \leq C|\lambda|^{-(1-\frac{\sigma}{2})}\|H\|_{B^{s-\sigma}_{q,1}(\HS)}
&\quad&\text{for $H \in B^{\nu}_{q,1}(\HS)^{m(N)}$}.
\end{aligned}\end{equation*}
\par
To estimate $\pd_\lambda \CT_{\Omega_+}(\lambda)$, we write
$$\pd_\lambda \CT_{\Omega_+}(\lambda)H 
= (\pd_\lambda \CS_{\HS}(\lambda))\CQ_\infty H + 
\CS_{\HS}(\lambda)(\pd_\lambda \CQ_\infty)H.
$$
Since $\CS_{\HS}(\lambda)$ has $(s, \sigma, q)$ properties, we have
$$\|(\lambda, \lambda\bar\nabla, \bar\nabla^2)
(\pd_\lambda \CS_{\HS}(\lambda))\CQ_\infty H\|_{B^s_{q,1}(\HS)} 
\leq C|\lambda|^{-1}\|\CQ_\infty H\|_{B^s_{q,1}(\HS)} 
\leq C|\lambda|^{-1}\|H\|_{B^s_{q,1}(\HS)}.$$
Moreover, 
$$\|(\lambda, \lambda\bar\nabla, \bar\nabla^2)
(\pd_\lambda \CS_{\HS}(\lambda))\CQ_\infty H\|_{B^s_{q,1}(\HS)} 
\leq C|\lambda|^{-(1-\frac{\sigma}{2})}\|\CQ_\infty H\|_{B^{s-\sigma}_{q,1}(\HS)}.
$$
Since $-1+1/q < s-\sigma < 1/q$, we can show that 
$$\|(\lambda, \lambda^{1/2}\bar\nabla, \bar\nabla^2)\CS_{\HS}(\lambda)H
\|_{B^{s-\sigma}_{q,1}(\HS)} \leq C\|H\|_{B^{s-\sigma}_{q,1}(\HS)}
$$
and so,  we may have 
\begin{equation}\label{est.cu}
\|\CO_\lambda \CU H\|_{H^{s-\sigma}_{q,1}(\HS)} 
\leq (1/2)\|H\|_{B^{s-\sigma}_{q,1}(\HS)}
\end{equation}
for $\lambda \in \Sigma_{\epsilon, \lambda_1}$.  Thus, we have
$$\|\CQ_\infty H\|_{B^{s-\sigma}_{q,1}(\HS)} 
\leq C\|H\|_{B^{s-\sigma}_{q,1}(\HS)},$$
which yields
$$\|(\lambda, \lambda^{1/2}\bar\nabla, \bar\nabla^2)
(\pd_\lambda \CS_{\HS}(\lambda))\CQ_\infty H\|_{B^s_{q,1}(\HS)} 
\leq C|\lambda|^{-(1-\frac{\sigma}{2})}\|H\|_{B^{s-\sigma}_{q,1}(\HS)}.
$$
To estimate $\pd_\lambda \CQ_\infty$, we write 
$$\pd_\lambda (\CO_\lambda \CU)^\ell 
= \sum (\CO_\lambda \CU)\cdots (\pd_\lambda \CO_\lambda \CU)
\cdots (\CO_\lambda \CU).$$
Since 
$$\pd_\lambda (\CO_\lambda \CU)
= (1/2)\lambda^{-1/2}R^2 \CR_{\HS}(\lambda)H 
+ \CO_\lambda(R^1\pd_\lambda \CS_{\HS}(\lambda)H, 
R^2\pd_\lambda \CS_{\HS}(\lambda)H), 
$$
we have
$$
\|\pd_\lambda (\CO_\lambda \CU)H\|_{B^s_{q,1}(\HS)}
\leq C|\lambda|^{-1}\|H\|_{B^s_{q,1}(\HS)}
$$
and so 
$$\|\pd_\lambda (\CO_\lambda \CU)^\ell\|_{B^s_{q,1}(\HS)}
\leq C\ell (1/2)^{\ell-1}|\lambda|^{-1}\|H\|_{B^s_{q,1}(\HS)},
$$
which yields 
$$\|\pd_\lambda \CQ_\infty H\|_{B^s_{q,1}(\HS)} 
\leq C|\lambda|^{-1}\sum_{\ell=1}^\infty \ell(1/2)^{\ell-1}
\leq 4C|\lambda|^{-1}\|H\|_{B^s_{q,1}(\HS)}.
$$
Thus, we have 
$$\|(\lambda, \lambda\bar\nabla, \bar\nabla^2)\CS_{\HS}(\lambda)
\pd_\lambda \CQ_\infty H\|_{B^s_{q,1}(\HS)} 
\leq C|\lambda|^{-1}\|H\|_{B^s_{q,1}(\HS)}.
$$
Moreover, 
\begin{align*}
\|\pd_\lambda(\CO_\lambda \CU H)\|_{B^s_{q,1}(\HS)} 
&\leq C(|\lambda|^{-1/2}\|\bar\nabla \CS_{\HS}(\lambda)H\|_{B^s_{q,1}(\HS)}
+ \|(\lambda, \lambda^{1/2}\bar\nabla, \bar\nabla^2)\pd_\lambda \CS_{\HS}
(\lambda)H\|_{B^s_{q,1}(\HS)}) \\
&\leq C|\lambda|^{-(1-\frac{\sigma}{2})}\|H\|_{B^{s-\sigma}_{q,1}(\HS)}.
\end{align*}
Thus, by \eqref{small.3.1} and \eqref{est.cu}, we have
\begin{align*}
\|\pd_\lambda \CQ_\infty H\|_{B^s_{q,1}(\HS)}
&\leq C\sum_{\ell=1}^\infty \ell(1/2)^{\ell-1}|\lambda|^{-(1-\frac{\sigma}{2})}
\|H\|_{B^{s-\sigma}_{q,1}(\HS)} \\
&\leq 4C|\lambda|^{-(1-\frac{\sigma}{2})}\|H\|_{B^{s-\sigma}_{q,1}(\HS)},
\end{align*}
which yields
$$\|(\lambda, \lambda^{1/2}\bar\nabla, \bar\nabla^2)
\CS_{\HS}(\pd_\lambda \CQ_\infty(\lambda)H\|_{B^s_{q,1}(\HS)}
\leq C|\lambda|^{-(1-\frac{\sigma}{2})}\|H\|_{B^{s-\sigma}_{q,1}(\HS)}.
$$
Summing up, we have 
\begin{equation}\label{deriv.4.4.1}\begin{aligned}
\|(\lambda, \lambda^{1/2}\bar\nabla, \bar\nabla^2)
\pd_\lambda \CT_{\Omega_+}(\lambda)H\|_{B^s_{q,1}(\HS)}
&\leq C|\lambda|^{-1}\|H\|_{B^s_{q,1}(\HS)}, \\
\|(\lambda, \lambda^{1/2}\bar\nabla, \bar\nabla^2)
\pd_\lambda \CT_{\Omega_+}(\lambda)H\|_{B^s_{q,1}(\HS)}
&\leq C|\lambda|^{-(1-\frac{\sigma}{2})}\|H\|_{B^{s-\sigma}_{q,1}(\HS)}.
\end{aligned}\end{equation}
Therefore, we see that $\CT_{\Omega_+}(\lambda)$
has $(s, \sigma, q)$ properties. \par 

Finally, according to \eqref{reduce.4}, 
for any $H =(H_1, H_2, H_3) \in B^\nu_{q,1}(\Omega_+)^{m(N)}$, 
we define $\CU_{\Omega_+}(\lambda)$ by setting
$$\CU_{\Omega_+}(\lambda)H =
\CA^\top (\CT_{\Omega_+}(\lambda)(\CA H_1, \CA H_2, 
(\nabla \Phi)^\top\CA H_3)\circ \Phi))\circ\Phi^{-1}.$$
Obvisouly, for any $(\tilde\bu, \tilde\bh) \in \CH^\nu_{q,1}(\Omega_+)$, 
$\bu = \CU_{\Omega_+}(\lambda)\CO_\lambda(\bg, \bh)$
is a solution of equations \eqref{resol.1}.  
The uniqueness follows from the existence of solutions to the dual problem. 
Since $\CT_{\Omega_+}(\lambda)$ has $(s, \sigma, q)$ properties, 
so $\tilde \CT_{\Omega_+}(\lambda)$ does. 
This completes the proof of Theorem \ref{thm:bent}. 
\end{proof}

\subsection{Spectral  analysis of generalized
 Lam\'e equations in $\Omega$, A proof of Theorem \ref{thm:Lame}}

In this subsection, we consider the equations \eqref{Eq:Lame} and 
we prove Theorem \ref{thm:Lame}. 
We only consider the case where $\Omega$  
is an exterior domain and $\tilde\eta_0(x)\not\equiv0$.
Other cases can be treated in the same manner. \par
First, we consider the problem in $(B_R)^c$ with large $R>0$. 
Let $\psi$ and $\tilde\psi$ be two $C^\infty(\BR^N)$ functions such 
that $\psi(x)$ equals to $1$ for $|x| >3$ and $0$ for $|x| <2$ and $\tilde\psi(x)$ equals to 
$1$ for $|x| > 2$ and $0$ for $|x| < 1$.  Set $\psi_R(x) = \psi(x/R)$ and $\tilde\psi_R(x)=
\tilde\psi(x/R)$.  Notice that $\psi_R(x)\tilde\psi_R(x) = \psi_R(x)$. 
Let $\lambda_0>0$ and  $\CS(\lambda) \in {\rm Hol}\, (\Sigma_{\epsilon, \lambda_0},
\CL(B^\nu_{q,1}(\BR^N), B^{\nu+2}_{q,1}(\BR^N)^N))$ 
be the operator given in Theorem \ref{thm:L.1} and then
$\bv_R=\CS(\rho_*\lambda)\tilde\psi_R\bg$ satisfies equations
\begin{equation}\label{pert.8.1.1}\begin{aligned}
\rho_*\lambda \bv_R - \DV(\alpha\BD(\bv_R) + (\beta-\alpha)
\dv \bv_R\BI) = \psi_R\bg&&\quad&\text{in $\BR^N$}.
\end{aligned}\end{equation}
for $\lambda \in \Sigma_{\epsilon, \rho_1^{-1}\lambda_0}$. 
Here, notice that $\rho_1 < \rho_* < \rho_2$. Let
$$A_R =\rho_* + \tilde\psi_R(x)(\eta_0(x)-\rho_*) = \rho_* + \tilde\psi_R(x)\tilde\eta_0(x). $$
We have
\begin{equation*}
A_R\lambda \bv_R - \DV(\alpha\BD(\bv_R) + (\beta-\alpha)\dv\bv_R\BI)
= \tilde\psi_R \bg - \CR_R(\lambda)\tilde\psi_R\bg\quad\text{in $\BR^N$},
\end{equation*}
where we have set
$$\CR_R(\lambda)\bff 
 = -\tilde\psi_R(x)\tilde\eta_0(x)\lambda \CS(\rho_*\lambda)\bff.
$$
By Lemma \ref{lem:prod} and Theorem \ref{thm:L.1}, we have
$$\|\CR_R(\lambda)\bff\|_{B^\nu_{q,1}(\BR^N)} 
\leq C\|\tilde\psi_R\tilde\eta_0\|_{B^{N/q}_{q,1}(\BR^N)}\|\bff\|_{B^{\nu}_{q,1}(\BR^N)}. 
$$
By Lemma 12 in \cite{KS23},  for any $\delta > 0$ there exists an $R_0>1$ such that 
$\|\tilde\psi_R \tilde\eta_0\|_{B^{N/q}_{q,1}(\Omega)} < \delta
$
for any $R > R_0$, and so   we have
$\|\CR_R(\lambda)\bff\|_{B^\nu_{q,1}(\Omega)} 
\leq C\delta\|\CO_\lambda(\bg, \bh)\|_{B^\nu_{q,1}(\Omega)}
$
for $\nu \in \{s-\sigma, s, s+\sigma\}$ and $R > R_0$. 
We choose $\delta > 0$ in such a way that 
$C\delta \leq 1/2$, we have 
$$\|\CR_R(\lambda)\bff \|_{B^\nu_{q,1}(\Omega)}
\leq (1/2)\|\bff\|_{B^\nu_{q,1}(\Omega)}$$
for $\nu \in \{s-\sigma, s, s+\sigma\}$ and $R > R_0$. 
In the following this $R$ is fixed.   Thus, we can define 
$$\CR_{R, \infty}(\lambda)=
(I-\CR_R(\lambda))^{-1} = \sum_{\ell=0}^\infty \CR_R(\lambda)^\ell.$$  
Let $\CS_{R}(\lambda)H= \CS(\rho_*\lambda)\CR_{R, \infty}(\lambda)\tilde\psi_RH_1$
for $\lambda \in \Sigma_{\epsilon, \rho_1^{-1}\lambda_1}$ and $H \in B^\nu_{q,1}(\HS)^{m(N)}$,
and then 
from \eqref{pert.8.1.1}, $\bw_R = \CS_{R}(\lambda)\CO_\lambda(\bg, \bh)$ satisfies equations 
\begin{equation}\label{per.8.1.2}
A_R\lambda \bw_R - \DV(\alpha\BD(\bw_R) + (\beta-\alpha)\dv\bw_R \BI) = 
\tilde\psi_R\bg \quad\text{in $\BR^N$}
\end{equation}
for any $\lambda \in \Sigma_{\epsilon, \rho_1^{-1}\lambda_1}$ and $(\bg, \bh) 
\in \CH^\nu_{q,1}(\Omega)$. 
Moreover, by Theorem \ref{thm:L.1}, 
\begin{equation}\label{est.8.1.1}\begin{aligned}
\|(\lambda, \lambda^{1/2}\bar\nabla, \bar\nabla^2)\CS_{R}(\lambda)H \|_{B^\nu_{q,1}(\Omega)}
&\leq C\|H\|_{B^\nu_{q,1}(\Omega)}&\quad&\text{for $H \in B^\nu_{q,1}(\Omega)^{m(N)}$}, \\
\|(\lambda^{1/2}\bar\nabla, \bar\nabla^2)\CS_{R}(\lambda)H \|_{B^s_{q,1}(\Omega)}
&\leq C|\lambda|^{-\frac{\sigma}{2}}\|H\|_{B^{s+\sigma}_{q,1}(\Omega)}
&\quad&\text{for $H \in B^{s+\sigma}_{q,1}(\Omega)^{m(N)}$}, \\
\|(1, \lambda^{-1/2}\bar\nabla) \CS_{R}(\lambda)H\|_{B^s_{q,1}(\Omega)}
& \leq C|\lambda|^{-(1-\frac{\sigma}{2})}\|H\|_{B^{s-\sigma}_{q,1}(\Omega)}
&\quad&\text{for $H \in B^s_{q,1}(\Omega)^{m(N)}$}
\end{aligned}\end{equation}
for $\lambda \in \Sigma_{\epsilon, \rho_1^{-1}\lambda_0}$. 
Moreover, employing the similar argument to the poof of \eqref{deriv.4.4.1}, we see
that 
\begin{equation}\label{est.8.1.1*}\begin{aligned}
\|(\lambda, \lambda^{1/2}\bar\nabla, \bar\nabla^2)
\pd_\lambda \CS_{R}(\lambda)H \|_{B^\nu_{q,1}(\Omega)}
&\leq C|\lambda|^{-1}
\|H\|_{B^\nu_{q,1}(\Omega)}&\quad&\text{for $H \in B^\nu_{q,1}(\Omega)^{m(N)}$}, \\
\|(\lambda, \lambda^{1/2}\bar\nabla, \bar\nabla^2)\pd_\lambda 
\CS_{R}(\lambda)H \|_{B^s_{q,1}(\Omega)}
&\leq C|\lambda|^{-(1-\frac{\sigma}{2})}\|H\|_{B^{s-\sigma}_{q,1}(\Omega)}
&\quad&\text{for $H \in B^{s}_{q,1}(\Omega)^{m(N)}$}
\end{aligned}\end{equation}
for $\lambda \in \Sigma_{\epsilon, \rho_1^{-1}\lambda_1}$.  
Let $\bu_R = \psi_R\CS_R(\lambda)\tilde\psi_R\bg$ and set 
\begin{align*}
U_R(\lambda)H &= \psi_R\DV(\alpha\BD(\CS_R(\lambda)\tilde\psi_RH_1) 
+ (\beta-\alpha)\dv(\CS_R(\lambda)\tilde\psi_RH_1)\BI) \\
& \qquad - \DV(\alpha\BD(\psi_R\CS_R(\lambda)\tilde\psi_RH_1) 
+ (\beta-\alpha)\dv(\psi_R\CS_R(\lambda)\tilde\psi_RH_1)\BI)
\end{align*}
for $\lambda \in \Sigma_{\epsilon, \rho_1^{-1}\lambda_1}$
and  $H \in B^\nu_{q,1}(\Omega)^{m(N)}$. 
From \eqref{est.8.1.1} and Lemma \ref{lem:dif1}, it follows that 
\begin{equation}\label{residue.8.1}
\|U_R(\lambda)H \|_{B^\nu_{q,1}(\Omega)} 
\leq C_R|\lambda|^{-1/2}\|\tilde \psi_RH_1\|_{B^\nu_{q,1}(\Omega)}
\end{equation}
for $\lambda \in \Sigma_{\epsilon, \rho_1^{-1}}$ and $H \in B^\nu_{q,1}(\Omega)^{m(N)}$.
Moreover, by \eqref{per.8.1.2} $\bu_R$ satisfies equations
\begin{equation*}
\eta_0 \lambda \bu_R - \DV(\alpha\BD(\bu_R) + (\beta-\alpha)\dv\bu_R \BI) = 
\psi_R\bg - U_R(\lambda)\CO_s(\bg, \bh) \quad\text{in $\BR^N$} 
\end{equation*}
for any $\lambda \in \Sigma_{\epsilon, \rho_1^{-1}\lambda_1}$
and $(\bg, \bh) \in \CH^\nu_{q,1}(\Omega)$. \par
Let $x_0 \in \Omega$ and let $d_{x_0}$ be a positive number such that 
$B_{4d_{x_0}}(x_0) \subset \Omega$.  Let $\varphi$ and $\tilde\varphi$ be 
two $C^\infty_0(\BR^N)$ functions such that 
$\varphi$ equals $1$ for $|x| < 1$ and $0$ for $|x| > 2$ and $\tilde\varphi(x)$
equals $1$ for $|x| <2$ and $0$ for $|x| >3$.  Set $\varphi_{x_0, d_{x_0}}(x)
=\varphi((x-x_0)/d_{x_0})$ and $\tilde\varphi_{x_0, d_0}(x)
= \tilde\varphi((x-x_0)/d_{x_0})$.  Notice that 
$\varphi_{x_0, d_{x_0}}\tilde\varphi_{x_0, d_0}= \varphi_{x_0, d_{x_0}}$. 
Let $\CS(\lambda)$ be the operator given in Theorem \ref{thm:L.1} again and set 
$\bv_{x_0}=\CS(\eta_0(x_0)\lambda)\tilde\varphi_{x_0, d_{x_0}}\bg$
for $\lambda \in \Sigma_{\epsilon, \rho_1^{-1}\lambda_0}$.  
Here, notice that $\rho_1 < \eta_0(x_0) < \rho_2$. Then, $\bv_{x_0}$
satisfies equations: 
\begin{equation}\label{pert.8.2.1}\begin{aligned}
\eta_0(x_0) \lambda \bv_{x_0} - \DV(\alpha\BD(\bv_{x_0}) + (\beta-\alpha)
\dv \bv_{x_0}\BI) = \tilde\varphi_{x_0, d_{x_0}}\bg&&\quad&\text{in $\BR^N$}.
\end{aligned}\end{equation}
 Let
$$A_{x_0} =\eta_0(x_0) + \tilde\varphi_{x_0, d_{x_0}}(x)(\eta_0(x)-\eta_0(x_0)). $$
We have
\begin{equation*}
A_{x_0}\lambda \bv_{x_0} - \DV(\alpha\BD(\bv_{x_0}) + (\beta-\alpha)\dv\bv_{x_0}\BI)
= \tilde\varphi_{x_0, d_{x_0}} \bg - \CR_{x_0}(\lambda)\tilde\varphi_{x_0, d_{x_0}}\bg
\quad\text{in $\BR^N$}
\end{equation*}
where we have set
$$\CR_{x_0}(\lambda)\bff 
 = -\tilde\varphi_{x_0, d_{x_0}}(x)(\eta_0(x)- \eta_0(x_0))
\lambda \CS(\lambda)\bff.
$$
By Lemma \ref{lem:prod} and Theorem \ref{thm:L.1}, we have
$$\|\CR_{x_0}(\lambda)\bff\|_{B^\nu_{q,1}(\BR^N)} 
\leq C\|\tilde\varphi_{x_0, d_{x_0}}(\eta_0(\cdot)- \eta_0(x_0))
\|_{B^{N/q}_{q,1}(\BR^N)}\|\bff\|_{B^{\nu}_{q,1}(\BR^N)}.
$$
By Appendix in \cite{DT22},  for any $\delta > 0$ there exists a $d_0$ 
uniformly with respect to $x_0$ such that 
$$\|\tilde\varphi_{x_0, d_{x_0}}(\eta_0(\cdot)-\eta_0(x_0))\|_{B^{N/q}_{q,1}(\BR^N)} < \delta
$$
provided $0 < d_{x_0} \leq d_0$,  and so   we have
$\|\CR_{x_0}(\lambda)\bff\|_{B^\nu_{q,1}(\Omega)} 
\leq C\delta\|\bff\|_{B^\nu_{q,1}(\Omega)}
$
for $\nu \in \{s-\sigma, s, s+\sigma\}$. We choose $\delta > 0$ in such a way that 
$C\delta \leq 1/2$, we have 
$$\|\CR_{x_0}(\lambda)\bff \|_{B^\nu_{q,1}(\Omega)}
\leq (1/2)\|\bff\|_{B^\nu_{q,1}(\Omega)}$$
for $\nu \in \{s-\sigma, s, s+\sigma\}$ and $0<d_{x_0} < d_0$. 
In the following $d_{x_0}$ is fixed. Thus, we can define $\CR_{x_0, \infty}(\lambda)=
(I-\CR_{x_0}(\lambda))^{-1} = \sum_{\ell=0}^\infty \CR_{x_0}(\lambda)^\ell$.  
Let $\CS_{x_0}(\lambda)H = \CS(\lambda)\CR_{x_0, \infty}(\lambda)\tilde\varphi_{x_0, d_{x_0}}H_1$, 
then from \eqref{pert.8.2.1}, $\bw_{x_0} = \CS_{x_0}(\lambda)\CO_\lambda(\bg, \bh)$ 
satisfies equations 
\begin{equation}\label{per.8.2.2}
A_{x_0}\lambda \bw_{x_0} - \DV(\alpha\BD(\bw_{x_0}) + (\beta-\alpha)\dv\bw_{x_0}\BI) = 
\tilde \varphi_{x_0, d_{x_0}}\bg \quad\text{in $\BR^N$}. 
\end{equation}
Moreover, by Theorem \ref{thm:L.1}, 
\begin{equation}\label{est.8.2.1}\begin{aligned}
\|(\lambda, \lambda^{1/2}\bar\nabla, \bar\nabla^2)
\CS_{x_0}(\lambda)H\|_{B^\nu_{q,1}(\Omega)}
&\leq C\|H\|_{B^\nu_{q,1}(\Omega)}&\quad&
\text{for $H \in B^\nu_{q,1}(\Omega)^{m(N)}$}, \\
\|(\lambda^{1/2}\bar\nabla, \bar\nabla^2)\CS_{x_0}(\lambda)H 
\|_{B^\nu_{q,1}(\Omega)}
&\leq C|\lambda|^{-\frac{\sigma}{2}}\|H\|_{B^{s+\sigma}_{q,1}(\Omega)}
&\quad&\text{for $H \in B^{s+\sigma}_{q,1}(\Omega)^{m(N)}$}, \\
\|(1, \lambda^{-1/2}\bar\nabla) \CS_{x_0}(\lambda)H\|_{B^s_{q,1}(\Omega)}
& \leq C|\lambda|^{-(1-\frac{\sigma}{2})}\|H\|_{B^{s-\sigma}_{q,1}(\Omega)}
&\quad&\text{for $H \in B^s_{q,1}(\Omega)^{m(N)}$}
\end{aligned}\end{equation}
 for $\lambda \in \Sigma_{\epsilon, \rho_1^{-1}\lambda_1}$.  
Moreover, employing the similar argument to the proof of \eqref{deriv.4.4.1}, we see
that 
\begin{equation}\label{est.8.2.1*}\begin{aligned}
\|(\lambda, \lambda^{1/2}\bar\nabla, \bar\nabla^2)\pd_\lambda
\CS_{x_0}(\lambda)H\|_{B^\nu_{q,1}(\Omega)}
&\leq C|\lambda|^{-1}\|H\|_{B^\nu_{q,1}(\Omega)}&\quad&
\text{for $H \in B^\nu_{q,1}(\Omega)^{m(N)}$}, \\
\|(\lambda, \lambda^{1/2}\bar\nabla, \bar\nabla^2)\pd_\lambda\CS_{x_0}(\lambda)H 
\|_{B^s_{q,1}(\Omega)}
&\leq C|\lambda|^{-(1-\frac{\sigma}{2})}\|H\|_{B^{s-\sigma}_{q,1}(\Omega)}
&\quad&\text{for $H \in B^{s-\sigma}_{q,1}(\Omega)^{m(N)}$}
\end{aligned}\end{equation}
for $\lambda \in \Sigma_{\epsilon, \rho_1^{-1}\lambda_1}$.   
Let $\bu_{x_0} = \varphi_{x_0, d_{x_0}}\bv_{x_0}$ and set 
\begin{align*}
U_{x_0}(\lambda)H &= \varphi_{x_0, d_{x_0}}\DV(\alpha\BD(\CS_{x_0}(\lambda)H) 
+ (\beta-\alpha)\dv(\CS_{x_0}(\lambda)H)\BI) \\
&\qquad 
- \DV(\alpha\BD(\varphi_{x_0, d_{x_0}}\CS_{x_0}(\lambda)H) + (\beta-\alpha)
\dv(\varphi_{x_0, d_{x_0}}\CS_{x_0}(\lambda)H)\BI).
\end{align*}
From \eqref{est.8.2.1} and Lemma \ref{lem:dif1}, it follows that 
\begin{equation}\label{residue.8.2}
\|U_{x_0}(\lambda)H\|_{B^\nu_{q,1}(\BR^N)} 
\leq C_{x_0}|\lambda|^{-1/2}\|H\|_{B^\nu_{q,1}(\Omega)}
\end{equation}
for any $\lambda \in \Sigma_{\epsilon, \rho_1^{-1}\lambda_1}$ and $H
\in B^\nu_{q,1}(\Omega)$. 
Moreover, by \eqref{per.8.2.2} $\bu_{x_0}$ satisfies equations
\begin{equation*}
\eta_0 \lambda \bu_{x_0} - \DV(\alpha\BD(\bu_{x_0}) 
+ (\beta-\alpha)\dv\bu_{x_0} \BI) = 
\varphi_{x_0}\bg - U_{x_0}(\lambda)
\CO_\lambda(\bg, \bh) \quad\text{in $\BR^N$}.
\end{equation*}
\par
For $x_1 \in \pd\Omega$, let $d_{x_1}>0$ be a small number and let $\Omega_+$
a bent half-space given in Subsection \ref{subsec.4.4} such that 
$B_{4d_{x_1}}(x_1) \cap \Omega \subset \Omega_+$.  Let $\lambda_1 \geq 1$ and 
$\CS_{\Omega_+} \in {\rm Hol}\,(\Sigma_{\epsilon, \lambda_1}, \CL(B^\nu_{q,1}(\Omega_+)^{m(N)},
B^{\nu+2}_{q,1}(\Omega_+)^N)) $
be the operator given in Theorem \ref{thm:bent}.
Note that $|\eta_0(x_1)\lambda| \geq \lambda_1$ for $\lambda \in \Sigma_{\epsilon, 
\rho_1^{-1}\lambda_1}$. Set 
$\bv_{x_1}=\CS_{\Omega_+}
(\eta_0(x_1)\lambda)\CO_\lambda\tilde\varphi_{x_0, d_{x_0}}(\bg, \bh)$, and then $\bv_{x_1}$
satisfies equations: 
\begin{equation*}
\begin{aligned}
\eta_0(x_1) \lambda \bv_{x_1} - \DV(\alpha\BD(\bv_{x_1}) + (\beta-\alpha)
\dv \bv_{x_1}\BI) = \tilde\varphi_{x_1, d_{x_1}}\bg&&\quad&\text{in $\Omega$}, \\
(\alpha\BD(\bv_{x_1}+(\beta-\alpha)\dv\bv_{x_1}\BI)\bn = \tilde\varphi_{x_1, d_{x_1}}\bh
&&\quad&\text{on $\pd\Omega$}.
\end{aligned}\end{equation*}
 Let
$$A_{x_1} =\eta_0(x_1) + \tilde\varphi_{x_1, d_{x_1}}(x)(\eta_0(x)-\eta_0(x_1)). $$
We have
\begin{equation*}
\begin{aligned}
A_{x_1}\lambda \bv_{x_1} - \DV(\alpha\BD(\bv_{x_1}) + (\beta-\alpha)\dv\bv_{x_1}\BI)
&= \tilde\varphi_{x_1, d_{x_1}} \bg - \CR_{x_1}^1(\lambda)\CO_\lambda(\bg, \bh)
&\quad&\text{in $\Omega$}, \\
(\alpha\BD(\bv_{x_1}+(\beta-\alpha)\dv\bv_{x_1}\BI)\bn
& = \tilde\varphi_{x_1, d_{x_1}}\bh &\quad&\text{on $\pd\Omega$}
\end{aligned}\end{equation*}
for any $\lambda \in \Sigma_{\epsilon, \rho_1^{-1}\lambda_1}$ and 
$H\in B^\nu_{q,r}(\Omega)$, where we have set
$$\CR_{x_1}^1(\lambda)H
 = -\tilde\varphi_{x_1, d_{x_1}}(x)(\eta_0(x)- \eta_0(x_1))
\lambda \CS_{\Omega_+}(\eta_0(x_1)\lambda)\tilde\varphi_{x_1, d_{x_1}}H.
$$
By Lemma \ref{lem:prod} and Theorem \ref{thm:bent}, we have
$$\|\CR^1_{x_1}(\lambda)H\|_{B^\nu_{q,1}(\BR^N)} 
\leq C\|\tilde\varphi_{x_1, d_{x_1}}(\eta_0(\cdot)- \eta_0(x_1))
\|_{B^{N/q}_{q,1}(\BR^N)}\|H\|_{B^{\nu}_{q,1}(\BR^N)}.
$$
By Appendix in \cite{DT22},  for any $\delta > 0$ there exists a $d_0$ 
uniformly with respect to $x_0$ such that 
$$\|\tilde\varphi_{x_1, d_{x_1}}(\eta_0(\cdot)-\eta_0(x_1))\|_{B^{N/q}_{q,1}(\BR^N)} < \delta
$$
provided $0 < d_{x_1} \leq d_0$,  and so   we have
$\|\CR^1_{x_1}(\lambda)H\|_{B^\nu_{q,1}(\Omega)} 
\leq C\delta\|H\|_{B^\nu_{q,1}(\Omega)}
$
for $\nu \in \{s-\sigma, s, s+\sigma\}$. We choose $\delta > 0$ in such a way that 
$C\delta \leq 1/2$, we have 
$$\|\CR^1_{x_1}(\lambda)H \|_{B^\nu_{q,1}(\Omega_+)}
\leq (1/2)\|H\|_{B^\nu_{q,1}(\Omega_+)}$$
for $\nu \in \{s-\sigma, s, s+\sigma\}$ for $0<d_{x_1} < d_0$. 
In the following $d_{x_1}$ is fixed.  Let $\CU_{x_1}$ be an operator defined by 
$\CU_{x_1}H = (\CR^1_{x_1}(\lambda)H, 0)$ for $H=(H_1, H_2, H_3) \in \CH^\nu_{q,1}(\Omega_+)$,  
and then 
\begin{equation}\label{residue.8.3}
(\tilde\varphi_{x_1, d_{x_1}} \bg - \CR_{x_1}^1(\lambda)\tilde\varphi_{x_1, d_{x_1}}\bg, 
\tilde\varphi_{x_1, d_{x_1}}\bh)
= (\bI-\CU_{x_1}\CO_\lambda)\tilde\varphi_{x_1,d_{x_1}}(\bg, \bh).
\end{equation}
For $\lambda \in \Sigma_{\epsilon, \rho_1^{-1}\lambda_1}$, 
\begin{equation}\label{residue.8.4}
 \|(\CO_\lambda \CU_{x_1}) H\|_{B^\nu_{q,1}(\Omega_+)}
\leq \|\CR_{x_1}(\lambda)H_1 \|_{B^\nu_{q,1}(\Omega_+)}
\leq (1/2)\|H\|_{B^\nu_{q,1}(\Omega_+)}.\end{equation}
Thus,   
\begin{equation*}
\|\CO_\lambda( \CU_{x_1}\CO_\lambda)^\ell H\|_{B^\nu_{q,1}(\Omega_+)}
= \|(\CO_\lambda\CU_{x_1})^\ell \CO_\lambda H\|_{B^\nu_{q,1}(\Omega_+)}
\leq (1/2)^\ell \|\CO_\lambda  H\|_{B^\nu_{q,1}(\Omega_+)}
\end{equation*}
and so 
$\CR_{x_1, \infty}(\lambda)= \sum_{\ell=0}^\infty(\CU\CO_\lambda)^\ell
= (\bI- \CU_{x_1}\CO_\lambda)^{-1}$  can be defined.  In fact, 
\begin{align*}
\|\CO_\lambda \CR_{x_1, \infty}(\lambda)H\|_{B^\nu_{q,1}(\Omega_+)}
&\leq \sum_{\ell=0}^\infty\|(\CO_\lambda \CU_{x_1})^\ell \CO_\lambda H\|_{B^\nu_{q,1}(\Omega_+)} 
\\
&\leq \sum_{\ell=0}^\infty (1/2)^\ell \|\CO_\lambda H\|_{B^\nu_{q,1}(\Omega_+)} 
=2\|\CO_\lambda  H\|_{B^\nu_{q,1}(\Omega_+)}.
\end{align*}
If we define $\bw_{x_1}$ by 
$\bw_{x_1} = \CS_{\Omega_+}(\eta_0(x_1)\lambda)\CO_\lambda 
\CR_{x_1, \infty}(\lambda)
\tilde\varphi_{x_1, d_{x_1}}(\bg, \bh)$,
then  from \eqref{residue.8.3} we see that $\bw_{x_1}$ 
satisfies equations 
\begin{equation}\label{per.8.3.2}\begin{aligned}
A_{x_1}\lambda \bw_{x_1} - \DV(\alpha\BD(\bw_{x_1}) + (\beta-\alpha)\dv\bw_{x_1}\BI) &= 
\tilde \varphi_{x_1, d_{x_1}}\bg &\quad&\text{in $\Omega$}, \\
(\alpha \BD(\bw_{x_1}) + (\beta-\alpha)\dv\bw_{x_1}\BI)\bn &= \tilde\varphi_{x_1, d_{x_1}}\bh
&\quad&\text{on $\pd\Omega$}.
\end{aligned}\end{equation}
In view of \eqref{residue.8.4}, we can define $\tilde\CR_{x_1, \infty}(\lambda)H 
= \sum_{\ell=0}^\infty(\CO_\lambda\CU)^\ell$, and then $\CO_\lambda \tilde\CR_{x_1, \infty}
= \CR_{x_1, \infty}\CO_\lambda$ and the orerator norm of $\tilde \CR_{x_1, \infty}(\lambda)$
does not exceed $2$.  Thus, $\bw_{x_1} = \CU_{x_1}(\lambda)\CO_\lambda(\bg, \bh)$ where 
$\CU_{x_1}$ is defined by $\CU_{x_1}(\lambda)H = \CS_{\Omega_+}(\eta_0(x_1)\lambda)
\tilde \CR_{x_1, \infty}(\lambda)\varphi_{x_1, d_{x_1}}H$, and by Theorem \ref{thm:bent},
\begin{equation}\label{est.8.3.1}\begin{aligned}
\|(\lambda, \lambda^{1/2}\bar\nabla, \bar\nabla^2)\CU_{x_1}(\lambda)H 
\|_{B^\nu_{q,1}(\Omega)}
&\leq C\|H\|_{B^\nu_{q,1}(\Omega)}&\quad&\text{for $H \in B^\nu_{q,1}(\Omega)^{m(N)}$}, \\
\|(\lambda^{1/2}\bar\nabla, \bar\nabla^2)\CU_{x_1}(\lambda)H \|_{B^s_{q,1}(\Omega)}
&\leq C|\lambda|^{-\frac{\sigma}{2}}\|H\|_{B^{s+\sigma}_{q,1}(\Omega)}
&\quad&\text{for $H \in B^{s+\sigma}_{q,1}(\Omega)^{m(N)}$}, \\
\|(1, \lambda^{-1/2}\bar\nabla)\CU_{x_1}(\lambda)H \|_{B^s_{q,1}(\Omega)}
& \leq C|\lambda|^{-(1-\frac{\sigma}{2})}\|H\|_{B^{s-\sigma}_{q,1}(\Omega)}
&\quad&\text{for $H \in B^s_{q,1}(\Omega)^{m(N)}$}.
\end{aligned}\end{equation}
Moreover, employing the similar argument to the poof of \eqref{deriv.4.4.1}, we see
that 
\begin{equation}\label{est.8.3.1*}\begin{aligned}
\|(\lambda, \lambda^{1/2}\bar\nabla, \bar\nabla^2)
\pd_\lambda \CU_{x_1}(\lambda)H 
\|_{B^\nu_{q,1}(\Omega)}
&\leq C|\lambda|^{-1}\|H\|_{B^\nu_{q,1}(\Omega)}
&\quad&\text{for $H \in B^\nu_{q,1}(\Omega)^{m(N)}$}, \\
\|(\lambda, \lambda^{1/2}\bar\nabla, \bar\nabla^2)
\pd_\lambda \CU_{x_1}(\lambda)H \|_{B^s_{q,1}(\Omega)}
&\leq C|\lambda|^{-(1-\frac{\sigma}{2})}\|H\|_{B^{s-\sigma}_{q,1}(\Omega)}
&\quad&\text{for $H \in B^{s}_{q,1}(\Omega)^{m(N)}$}
\end{aligned}\end{equation}
for any $\lambda \in \Sigma_{\epsilon, \rho_1^{-1}\lambda_1}$. 
Let $\bu_{x_1} = \varphi_{x_1, d_{x_1}}\bw_{x_1}$ and set 
\begin{align*}
U_{x_1}^1(\lambda)H &= \varphi_{x_1, d_{x_1}}\DV(\alpha\BD(\CU_{x_1}(\lambda)H)
+ (\beta-\alpha)\dv(\CU_{x_1}(\lambda)H\BI)\\
&\quad 
- \DV(\alpha\BD(\varphi_{x_1, d_{x_1}}\CU_{x_1}(\lambda)H) + (\beta-\alpha)
\dv(\varphi_{x_1, d_{x_1}}\CU_{x_1}(\lambda)H)\BI), \\
U_{x_1}^2(\lambda)H &= \varphi_{x_1, d_{x_1}}(\alpha\BD(\CU_{x_1}(\lambda)H) 
+ (\beta-\alpha)\dv(\CU_{x_1}(\lambda)H)\BI)\bn \\
&\quad 
- (\alpha\BD(\varphi_{x_1, d_{x_1}}\CU_{x_1}(\lambda)H) + (\beta-\alpha)
\dv(\varphi_{x_1, d_{x_1}}\CU_{x_1}(\lambda)H)\BI)\bn.
\end{align*}
From \eqref{est.8.3.1} and Lemma \ref{lem:dif1}, it follows that 
\begin{equation}\label{residue.8.2*}
\|\CO_\lambda(U_{x_1}^1(\lambda)H, U_{x_1}^2(\lambda)H)\|_{B^\nu_{q,1}(\Omega)} 
\leq C_{x_1}|\lambda|^{-1/2}\|H\|_{B^\nu_{q,1}(\Omega)}
\end{equation}
for any $\lambda \in \Sigma_{\epsilon, \rho_1^{-1}\lambda_1}$ and 
$H \in B^\nu_{q,1}(\Omega)^{m(N)}$. 
Moreover, by \eqref{per.8.3.2} $\bu_{x_1}$ satisfies equations
\begin{equation*}
\begin{aligned}
\eta(x)\lambda \bu_{x_1} - \DV(\alpha\BD(\bu_{x_1}) + (\beta-\alpha)\dv\bu_{x_1} \BI) 
& = \varphi_{x_1, d_{x_1}}\bg -U^1_{x_1}(\bg, \bh)&\quad&\text{in $\Omega$}, \\
(\alpha \BD(\bu_{x_1}) + (\beta-\alpha)\dv\bu_{x_1}\BI)\bn &= \varphi_{x_1, d_{x_1}}\bh
-U^2_{x_1}(\bg, \bh)&\quad&\text{on $\pd\Omega$}.
\end{aligned}\end{equation*}
\par

Now, we shall show the theorem. 
let $\overline{\Omega \cap B_{2R}} = \{x \in \Omega \cup \pd\Omega \mid |x| \leq 2R\}$. 
Notice that $\Omega \cup \pd\Omega = (B_{2R})^c \cup \overline{\Omega \cap B_{2R}}.$
Since $\overline{\Omega \cap B_{2R}}$ is a compact set, 
there exist a finte set $\{x^0_j\}_{j=1}^{m_0}$ of 
points of $\Omega$ and a finite set $\{x^1_j\}_{j=1}^{m_1}$ of points of  $\pd\Omega$ such that 
$\overline{\Omega} \subset (B_{2R})^c \cup (\bigcup_{j=1}^{m_0} B_{d_{x^0_j}/2}(x^0_j)) 
\cup (\bigcup_{j=1}^{m_1}B_{d_{x^1_j}/2}(x^1_j))$. 
Let $\Phi(x) = \varphi_R(x) + (\sum_{j=1}^{m_0} \varphi_{x^0_j}(x)) + (\sum_{j=1}^{m_1}
\varphi_{x^1_j}(x))$.  Obviously, $\Phi(x) \in C^\infty(\overline{\Omega})$ 
and $\Phi(x) \geq 1$
for $x \in \overline{\Omega}$.  Thus, set 
$\omega_0(x) = \varphi_{R}(x)/\Phi(x)$, $\omega_j(x) =  \varphi_{x^0_j}(x)/\Phi(x)$ 
($j=1, \ldots, m_0$), 
and $\omega_{m_0+j}(x) = \varphi_{x^1_j}(x)/\Phi(x)$ ($j=1, \ldots, m_1$).  Then, 
$\{\omega_j\}_{j=0}^{m_0+m_1}$ is a partition of unity on $\overline{\Omega}$. 
Let us define the parametrix operator 
$\CT_\Omega(\lambda)$ and the remainder operator 
$\CU^i(\lambda)$ by 
\begin{alignat*}2
\CT_\Omega(\lambda)H &=\omega_0\CS_R(\lambda)H
+ \sum_{j=1}^{m_0} \omega_j\CS_{x^0_j}(\lambda)H
+ \sum_{j=1}^{m_1} \omega_{m_0+j}\CU_{x^1_j}(\lambda)H, & & \\
\CU^1(\lambda)H &=U_R(\lambda)H
+  \sum_{j=1}^{m_0} U^1_{x^0_j}(\lambda)H
+ \sum_{j=1}^{m_1} U^1_{x^1_{m_0+j}}(\lambda)H, &
&\CU^2(\lambda)H = \sum_{j=1}^{m_1} U^2_{x^1_{m_0 +j}}(\lambda)H.
\end{alignat*}
From \eqref{est.8.1.1}, \eqref{est.8.1.1*}, \eqref{est.8.2.1}, \eqref{est.8.2.1*},
\eqref{est.8.3.1}, \eqref{est.8.3.1*} and Lemma \ref{lem:dif1}, it follows that 
$\CT_\Omega(\lambda)$ has $(s, \sigma, q)$ properties
for $\lambda \in \Sigma_{\epsilon, \lambda_2}$. 
Moreover,  we see that $\bv=\CT_\Omega(\lambda)\CO_\lambda(\bg, \bh)$, 
$\CU^1(\lambda)\CO_\lambda(\bg, \bh)$,  
and $\CU^2(\lambda)\CO_\lambda(\bg, \bh)$ satisfy equations:
\begin{equation}\label{eq:8.1*}\begin{aligned}
\lambda \bv- \DV(\alpha\BD(\bv) + (\beta-\alpha)\dv\bv\BI) 
= \bg- \CU^1(\lambda)\CO_\lambda(\bg, \bh)
&&\quad&\text{in $\Omega$},  \\
(\alpha\BD(\bv)+(\beta-\alpha)\dv\bv\BI)\bn 
= \bh - \CU^2(\lambda)\CO_\lambda(\bg, \bh)
&&\quad&\text{on $\pd\Omega$}. 
\end{aligned}\end{equation}
Let $\lambda_2= \lambda_1\rho_1^{-1} $.
Set $\CU(\lambda)
= (\CU^1(\lambda), \CU^2(\lambda))$. By \eqref{residue.8.1},
\eqref{residue.8.2} and \eqref{residue.8.2*},  we have
$$
\|\CO_\lambda\CU(\lambda)H\|_{\CH^\nu_{q,1}(\Omega)} 
\leq C|\lambda|^{-1/2}\| H\|_{B^\nu_{q,1}(\Omega)} 
$$
for $\lambda \in \Sigma_{\epsilon, \lambda_2}$, because the summation is 
finite number.   Choosing $\lambda_3\geq\lambda_2$ so large 
that $C\lambda_3^{-1/2} \leq 1/2$, we have 
\begin{equation}\label{residue.8.4'}
\|\CO_\lambda\CU(\lambda)H\|_{B^\nu_{q,1}(\Omega)}
\leq (1/2)\|H\|_{B^\nu_{q,1}(\Omega)} 
\end{equation}
for $\lambda \in \Sigma_{\epsilon, \lambda_3}$. 
Let 
$$\CW(\lambda)
=  \sum_{\ell=0}^\infty (\CU(\lambda)\CO_\lambda)^\ell.$$
Noting that $\lambda_3^{-1/2} \leq 1$, we see that
\begin{equation*}
\begin{aligned}
\|\CW(\lambda)(\bg, \bh)\|_{\CH^\nu_{q,1}(\Omega)}
&\leq \|\CO_\lambda\CW(\lambda)(\bg, \bh)\|_{B^\nu_{q,1}(\Omega)} \\
&\leq \sum_{\ell=0}^\infty\|(\CO_\lambda\CU(\lambda))^\ell \CO_\lambda
(\bg, \bh)\|_{B^\nu_{q,1}(\Omega)} \\
& \leq \sum_{\ell=0}^\infty(1/2)^\ell\|\CO_\lambda(\bg, \bh)\|_{B^\nu_{q,1}(\Omega)}
= 2\|\CO_\lambda(\bg, \bh)\|_{B^\nu_{q,1}(\Omega)}\\
&\leq 2\max(1, |\lambda|^{1/2})\|(\bg, \bh)\|_{\CH^{\nu}_{q,1}(\Omega)}.
\end{aligned}\end{equation*}
From \eqref{eq:8.1*} we see that 
$\bu= \CT_\Omega(\lambda)\CO_\lambda\CW(\lambda)(\bg, \bh)$ satisfies 
equations:
$$\begin{aligned}
\eta_0\lambda \bu- \DV(\alpha\BD(\bu) + (\beta-\alpha)\dv\bu\BI) = \bg
&&\quad&\text{in $\Omega$}, \\
(\alpha\BD(\bu)+(\beta-\alpha)\dv\bu\BI)\bn = \bh
&&\quad&\text{on $\pd\Omega$}. 
\end{aligned}$$
Let $\tilde\CW$ be an operator defined by 
$\tilde\CW = \sum_{\ell=0}^\infty(\CO_\lambda \CU(\lambda))^\ell$,
and then by  \eqref{residue.8.4'}, we see that the operator norm of 
$\tilde \CW$ does not exceed $2$ and $\CO_\lambda\CW = \tilde\CW\CO_\lambda$.
Thus, defining $\CV(\lambda)$ by $\CV(\lambda)
= \CT_\Omega(\lambda)\tilde\CW(\lambda)$, we see that $\bu=
\CV(\lambda)\CO_\lambda(\bg, \bh)$ 
 is a  solution of equations \eqref{Eq:Lame}.  Moreover, since $\CT_\Omega(\lambda)$
has $(s, \sigma, q)$ properties, so $\CV(\lambda)$ does.  
The uniqueness of solutions follows from the existence of solutions to the dual problem. 
This completes the proof of Theorem \ref{thm:Lame}.

\section{About the general resolvent problems for the
Stokes equations}

In this section, we consider the generalized resolvent problem for 
 Stokes equations with non-homogeneous free boundary conditions, which read as
\begin{equation}\label{stokes.2}\left\{\begin{aligned}
\lambda \rho + \eta_0 \dv \bu = f&&\quad&\text{in $\Omega$},  \\
\eta_0\lambda \bu - \DV(\alpha\BD(\bu) + (\beta-\alpha)\dv\bu\BI
 - P'(\eta_0)\rho\BI) = \bg
&&\quad&\text{in $\Omega$}, \\
(\alpha\BD(\bu)+(\beta-\alpha)\dv\bu\BI - P'(\eta_0)\rho\BI)\bn = \bh&&\quad
&\text{on $\pd\Omega$}.
\end{aligned}\right.\end{equation}
Recall that 
$\eta_0(x) = \rho_* + \tilde\eta_0(x)$, where $\rho_*$ is a positive constant and 
$\tilde\eta_0 \in B^{N/q}_{q,1}(\Omega)$, that
 $P(s)$ is a $C^\infty$ function defined 
for $s \in (0, \infty)$ such that $P'(s)>0$, and that Assumption
\ref{assump.1} holds. 
For any domain $D \subset \BR^N$,  let 
$\CJ^\nu_{q,1}(D) = B^{\nu+1}_{q,1}(D) \times B^{\nu}_{q,1}(D)\times 
B^{\nu+1}_{q,1}(D)$
and  $\CD^\nu_{q,1}(D)= 
B^{\nu+1}_{q,1}(D)\times B^{\nu+2}_{q,1}(D)$.  Their norms 
$\|\cdot\|_{\CJ^\nu_{q,1}(D)}$ and 
$\|\cdot\|_{\CD^\nu_{q,1}(D)}$ are defined by 
\begin{align*}
\|(f, \bg, \bh)\|_{\CJ^\nu_{q,1}(D)} &= \|f\|_{B^{\nu+1}_{q,1}(D)}
+ \|\bg\|_{B^\nu_{q,1}(D)} + \|\bh\|_{B^{\nu+1}_{q,1}(D)}, \\
\|(\rho, \bu)\|_{\CD^\nu_{q,1}(D)} &= \|\rho\|_{B^{\nu+1}_{q,1}(D)}
+ \|\bu\|_{B^{\nu+2}_{q,1}(D)}.
\end{align*}
Moreover, for $(f, \bg, \bh) \in \CJ^\nu_{q,1}(D)$, we set
$$\|(f, \CO_\lambda(\bg, \bh))\|_{B^{\nu+1, \nu}_{q,1}(D)}
= \|f\|_{B^{\nu+1}_{q,1}(D)} + \|(\bg, \lambda^{1/2}\bh, \nabla\bh)
\|_{B^\nu_{q,1}(D)}.$$
In what follows, 
we consider the two cases where $\rho_0=\rho_*$ and $\tilde\eta_0\not\equiv0$.
We shall show the following theorem.
\begin{thm}\label{thm:main4} Assume that the following conditions \thetag1 or \thetag2 holds. 
\begin{itemize}
\item[\thetag1] If $\eta_0(x) = \rho_*$, then $1 < q < \infty$ and $-1+1/q < s < 1/q$. 
\item[\thetag2] If $\tilde\eta_0(x)\not\equiv0$, then 
$N-1 < q < 2N$ and $-1+N/q \leq  s < 1/q$. 
\end{itemize}
Then, there exists a large positive number 
$\lambda_4 \geq 1$ such that for any $\lambda \in \Sigma_{\epsilon, \lambda_4}$, 
$(f, \bg, \bh) \in \CJ^s_{q,1}(\Omega)$, 
problem \eqref{stokes.2} admits unique solutions $(\rho, \bu)
\in \CD^s_{q,1}(\Omega)$ satisfying the estimate:
\begin{align*}
\|\lambda(\rho, \bu)\|_{\CH^s_{q,1}(\Omega)}
+ \|(\rho, \bu)\|_{\CD^s_{q,1}(\Omega)}
\leq C\|(f, \bg, \bh)\|_{\CJ^s_{q,1}(\Omega)}.
\end{align*}
\par
 Moreover, there exist three operators $\CB_v(\lambda)$, 
$\CC_m(\lambda)$ 
and $\CC_v(\lambda)$ with 
\begin{align*}
&\CB_v(\lambda) \in {\rm Hol}\, (\Sigma_{\epsilon, \lambda_4}, 
\CL(B^s_{q,1}(\Omega)^{m(N)}, B^{s+2}_{q,1}(\Omega)^N)), \\
&\CC_m(\lambda) \in {\rm Hol}\, (\Sigma_{\epsilon, \lambda_2}, 
\CL(B^{s+1}_{q,1}(\Omega)\times B^{s}_{q,1}(\Omega)^{m(N)}, 
B^{s+1}_{q,1}(\Omega)), \\
&\CC_v(\lambda) \in {\rm Hol}\, (\Sigma_{\epsilon, \lambda_2}, 
\CL(B^{s+1}_{q,1}(\Omega)\times B^{s}_{q,1}(\Omega)^{m(N)}),
B^{s+2}_{q,1}(\Omega)^N)
\end{align*}
such that $\CB_v(\lambda)$ has $(s, \sigma, q)$ properties  and 
$\CC_m(\lambda)$ and $\CC_v(\lambda)$ have generalized resolvent properties
for $X = B^{s+1}_{q,1}(\Omega)\times B^s_{q,1}(\Omega)^{m(N)}$ and 
for any $\lambda \in \Sigma_{\epsilon, \lambda_4}$ in the sense of Definition
\ref{dfn.2} and solutions $\rho$ and $\bu$ are represented by 
$\rho = \CC_m(\lambda)(f, \CO_\lambda(\bg, \bh))$ and $\bu = \CB_v(\lambda)
\CO_\lambda(\bg, \bh) + \CC_v(\lambda)(f, \CO_\lambda(\bg, \bh))$. 
\end{thm}
\begin{proof}
In what follows, we shall prove the theorem 
in the case where $\tilde\eta_0\not\equiv0$ only.
In fact, in the case where $\eta_0=\rho_*$  the theorem can be proved 
by using the similar argument. 
In \eqref{stokes.2}, setting $\rho= \lambda^{-1}(f- \eta_0\dv \bu)$
and inserting this formula into the second equations, we have
\begin{equation}\label{stokes.3} \begin{aligned}
\eta_0\lambda \bu -\DV( \alpha\BD(\bu)+( \beta-\alpha)\dv\bu \BI
+ \lambda^{-1}P'(\eta_0)\eta_0\dv\bu \BI)
 &= \bg
-\lambda^{-1}\nabla(P'(\eta_0)f)
&\quad&\text{in $\Omega$}, \\
(\alpha\BD(\bu)+(\beta-\alpha)\dv\bu \BI +\lambda^{-1}P'(\eta_0)\eta_0\dv\bu\BI)\bn
&=\bh + \lambda^{-1}P'(\eta_0)f\bn&\quad&\text{on $\pd\Omega$}. 
\end{aligned}\end{equation}
For a while, setting $\bG = \bg - \lambda^{-1}\nabla(P'(\eta_0)f)$
and $\bH = \bh + \lambda^{-1}P'(\eta_0)f\bn$,  we shall solve
equations:
\begin{equation}\label{stokes.4}\begin{aligned}
\eta_0\lambda \bu -\DV( \alpha\BD(\bu)+( \beta-\alpha)\dv\bu \BI 
+ \lambda^{-1}P'(\eta_0)\eta_0\dv\bu \BI)
 &= \bG
&\quad&\text{in $\Omega$}, \\
(\alpha\BD(\bu)+(\beta-\alpha)\dv\bu \BI +\lambda^{-1}P'(\eta_0)\eta_0\dv\bu\BI)\bn
&=\bH&\quad&\text{on $\pd\Omega$}. 
\end{aligned}\end{equation}
Let $\CV_\Omega(\lambda) \in {\rm Hol}\,(\Sigma_{\epsilon, \lambda_3},
\CL(B^\nu_{q,1}(\Omega)^{m(N)}, B^{\nu+2}_{q,1}(\Omega)^N))$ 
be the solution operator of equations. Insert the formula  
$\bu = \CV_\Omega(\lambda)\CO_\lambda(\bG, \bH)$ into \eqref{stokes.4}
to obtain
\begin{equation*}
\begin{aligned}
\eta_0\lambda \bu -\DV( \alpha\BD(\bu)+( \beta-\alpha)\dv\bu \BI 
+ \lambda^{-1}P'(\eta_0)\eta_0\dv\bu \BI)
&= \bG - \lambda^{-1}\nabla(P'(\eta_0)\eta_0\dv\bu)
&\quad&\text{in $\Omega$}, \\
(\alpha\BD(\bu)+(\beta-\alpha)\dv\bu \BI 
+\lambda^{-1}P'(\eta_0)\eta_0\dv\bu\BI)\bn
&=\bH + \lambda^{-1}(P'(\eta_0)\eta_0\dv\bu\BI)\bn
&\quad&\text{on $\pd\Omega$}. 
\end{aligned}\end{equation*}
We define an operator $\CP(\lambda)$ by 
\begin{align*}
\CP(\lambda) H &= (\nabla(P'(\eta_0)\eta_0\dv\CV_\Omega(\lambda)H), 
-P'(\eta_0)\eta_0\dv(\CV_\Omega(\lambda)H)\bn).
\end{align*}
Then, we have 
\begin{equation}\label{residue.4.1}
(\bG, \bH) - (\lambda^{-1}\nabla(P'(\eta_0)\eta_0\dv\bu), 
-\lambda^{-1} (P'(\eta_0)\eta_0\dv\bu\BI)\bn)
= (\bI -\lambda^{-1}\CP\CO_\lambda)(\bG, \bH).
\end{equation}
We will show that 
\begin{equation}\label{est.p1}
\|\CO_\lambda\CP(\lambda)H\|_{B^s_{q,1}(\Omega)}
\leq C(\rho_*, \|\tilde\eta_0\|_{B^{s+1}_{q,1}(\Omega)})
\|H\|_{B^s_{q,1}(\Omega)}.
\end{equation}
Here and in what follows, $C(\rho_*, \|\tilde\eta_0\|_{B^{s+1}_{q,1}(\Omega)})$ 
denotes some constant 
depending on $\rho_*$
and $\|\eta_0\|_{B^{s+1}_{q,1}(\Omega)}$ in  the case where $\tilde\eta_0 \not\equiv 0$.
If we consider the case where $\eta_0=\rho_*$, 
$C(\rho_*, \|\tilde\eta_0\|_{B^{s+1}_{q,1}(\Omega)})$ is replaced with simply a constant 
$C(\rho_*)$.
\par 
To prove \eqref{est.p1}, we shall use Lemma \ref{lem:Hasp} and the fact that
$B^{s+1}_{q,1}$ is a Banach algebra. 
In fact, noting that $N/q \leq s+1$ by Lemma \ref{lem:prod}, we have
\begin{align}
\|uv\|_{B^{s+1}_{q,1}(\Omega)} &\leq \|(\nabla u)v\|_{B^s_{q,1}(\Omega)}
+\|u(\nabla v)\|_{B^s_{q,1}(\Omega)} +\|uv\|_{B^s_{q,1}(\Omega)} \nonumber \\
&\leq C(\|u\|_{B^{\nu+1}_{q,1}(\Omega)}\|v\|_{B^{N/q}_{q,1}(\Omega)}
+ \|u\|_{B^{N/q}_{q,1}(\Omega)}\|v\|_{B^{\nu+1}_{q,1}(\Omega)}
+ \|u\|_{B^s_{q,1}(\Omega)}\|v\|_{B^{N/q}_{q,1}(\Omega)} ) \nonumber \\
&\leq C\|u\|_{B^{s+1}_{q,1}(\Omega)}\|v\|_{B^{s+1}_{q,1}(\Omega)}.
\label{banach.1}
\end{align}
\par
To prove \eqref{est.p1},  recalling that $\eta_0 = \rho_* +\tilde\eta_0$, 
we write  $P'(\eta_0)\eta_0 = P'(\rho_*)\rho_* + \CP_1(\tilde\eta_0)$, 
where we have set
$$\CP_1(r)= P'(\rho_*)r + \int^1_0P''(\rho_* + \theta r)\,d\theta r(\rho_*+r)
$$
with $r = \tilde\eta_0$.  Note that $\CP_1(0)=0$ and 
$\rho_1-\rho_* \leq \tilde \eta_0(x) \leq \rho_2-\rho_*$ as follows 
from Assumption \ref{assump.1}.  By Lemma \ref{lem:Hasp}, we have
$$\|\CP_1(\tilde\eta_0)\|_{B^{s+1}_{q,1}(\Omega)}
\leq C\|\tilde\eta_0\|_{B^{s+1}_{q,1}(\Omega)}.
$$
Thus, 
\begin{align*}
\|\nabla(P'(\eta_0)\eta_0\dv\CV_\Omega(\lambda)H)\|_{B^s_{q,1}(\Omega)}
&\leq |P'(\rho_*)\rho_*|\|\nabla\dv\CV_\Omega(\lambda)H\|_{B^s_{q,1}(\Omega)}
+ \|\CP_1(\tilde\eta_0)
\dv\CV_\Omega(\lambda)H\|_{B^{s+1}_{q,1}(\Omega)} \nonumber \\
&\leq |P'(\rho_*)\rho_*|\|\CV_\Omega(\lambda)H\|_{B^{s+2}_{q,1}(\Omega)}
+ C\|\tilde \eta_0\|_{B^{s+1}_{q,1}(\Omega)}
\|\dv\CV_\Omega(\lambda)H\|_{B^{s+1}_{q,1}(\Omega)}.
\end{align*}
Using Theorem \ref{thm:Lame}, we have
$$\|\nabla(P'(\eta_0)\eta_0\dv\CV_\Omega(\lambda)H)\|_{B^s_{q,1}(\Omega)}
\leq C(\rho_*, \|\tilde\eta_0\|_{B^{s+1}_{q,1}(\Omega)})
\|H\|_{B^\nu_{q,1}(\Omega)}.
$$
Analogously, 
\begin{align*}
&\|(\lambda^{1/2}, \nabla)(P'(\eta_0)\eta_0\dv
\CV_\Omega(\lambda)H)\bn\|_{B^s_{q,1}(\Omega)} \\
&\quad 
\leq C(\|\nabla(P'(\eta_0)\eta_0\dv\CV_\Omega(\lambda)H)
\|_{B^s_{q,1}(\Omega)} +
|\lambda|^{1/2}\|P'(\eta_0)\eta_0\dv\CV_\Omega(\lambda)H
\|_{B^s_{q,1}(\Omega)}).
\end{align*}
Here and in the following, we may assume that $|\lambda| \geq \lambda_3 \geq 1$. 
We see that 
\begin{align*}
\|P'(\eta_0)\eta_0\dv\CV_\Omega(\lambda)H\|_{B^s_{q,1}(\Omega)}
\leq (P'(\rho_*)\rho_*
+ C\|\CP_1(\tilde\eta_0)\|_{B^{N/q}_{q,1}(\Omega)})
\|\dv\CV_\Omega(\lambda)H\|_{B^s_{q,1}(\Omega)}.
\end{align*}
Therefore, we have 
$$
\|(\lambda^{1/2}, \nabla)(P'(\eta_0)\eta_0\dv\CV_\Omega(\lambda)H
)\bn\|_{B^s_{q,1}(\Omega)}
\leq C(\rho_*, \|\tilde\eta_0\|_{B^{s+1}_{q,1}(\Omega)})
\|(\lambda^{1/2}\bar\nabla, \bar\nabla^2)
\CV_\Omega(\lambda)H\|_{B^s_{q,1}(\Omega)}.
$$
Combining these estimates with  Theorem \ref{thm:Lame}  yields \eqref{est.p1}.  \par
Choosing $\lambda_4 \geq \lambda_3$ in such a way that 
$C(\rho_*, \|\tilde\eta_0\|_{B^{s+1}_{q,1}(\Omega)})\lambda_4^{-1} \leq 1/2$ 
in \eqref{est.p1}, we have
\begin{equation}\label{residue.4.2}
\|\lambda^{-1}\CO_\lambda\CP(\lambda)H\|_{B^s_{q,1}(\Omega)} 
\leq (1/2)\|H\|_{B^s_{q,1}(\Omega)}
\end{equation}
for any $\lambda \in \Sigma_{\epsilon, \lambda_4}$.  Let us define 
$\CP_\infty(\lambda)$ by
\begin{equation*}
\CP_\infty(\lambda) = 
\sum_{\ell=0}^\infty (\lambda^{-1}\CP(\lambda)\CO_\lambda)^\ell
\end{equation*}
for $\lambda \in \Sigma_{\epsilon, \lambda_4}$.  In view of  \eqref{residue.4.2}, 
we see that $\CP_\infty(\lambda) \in \CL(\CH^\nu_{q,1}(\Omega))$. 
Let $\bv = \CV_\Omega(\lambda)\CO_\lambda
\CP_\infty(\lambda)(\bG, \bH)$, and 
then by \eqref{residue.4.1} $\bv$ is a unique solution of the following equations: 
\begin{equation*}
\begin{aligned}
\eta_0\lambda \bv -\DV( \alpha\BD(\bv)+( \beta-\alpha)\dv\bv\BI
+ \lambda^{-1}P'(\eta_0)\eta_0\dv\bv \BI)
&= \bG
&\quad&\text{in $\Omega$}, \\
(\alpha\BD(\bv)+(\beta-\alpha)\dv\bv\BI 
+\lambda^{-1}P'(\eta_0)\eta_0\dv\bv\BI)\bn
&=\bH&\quad&\text{on $\pd\Omega$}. 
\end{aligned}\end{equation*}
The uniqueness follows
from the existence of solutions to the dual problem, 
which has essentially the same forms. \par
Let $\tilde \CP_\infty(\lambda) 
= \sum_{\ell=0}^\infty(\lambda^{-1} \CO_\lambda \CP(\lambda))^\ell$, 
and then by \eqref{residue.4.2}, 
$\tilde \CP_\infty(\lambda) \in \CL(B^s_{q,1}(\Omega))$ and its operator norm
does not exceed $2$.  Moreover, we have
$\CO_\lambda \CP_\infty(\lambda) = \tilde\CP_\infty(\lambda)
\CO_\lambda$.  Let 
$\CW^2_\Omega(\lambda)H = \CV_\Omega(\lambda)\tilde\CP_\infty H$,
and then  $\bv = \CW^2_\Omega(\lambda)\CO_\lambda(\bG, \bH)$. 
Moreover, by Theorem \ref{thm:Lame} and the similar argument to
the proof of \eqref{deriv.4.4.1}, we see that 
$\CW^2_\Omega(\lambda)$ has $(s, \sigma, q)$ properties
for any $\lambda \in \Sigma_{\epsilon, \lambda_4}$. \par
Setting $\bG = \bg-\lambda^{-1}\nabla(P'(\eta_0)f)$ and 
$\bH = \bh+\lambda^{-1}P'(\eta_0)f \bn$, we define $\bu$ by 
\begin{equation}\label{sol.9.1}
\bu = \CW^2_\Omega(\lambda)\CO_\lambda
(\bg-\lambda^{-1}\nabla(P'(\eta_0)f), \bh+\lambda^{-1}P'(\eta_0)f \bn).
\end{equation}
Obvisouly, $\bu$ is a solution of equations \eqref{stokes.3} and satisfies the estimate:
$$\|\lambda \bu\|_{B^s_{q,1}(\Omega)} + \|\bu\|_{B^{s+2}_{q,1}(\Omega)}
\leq C(\rho_*, \|\tilde\eta_0\|_{B^{s+1}_{q,1}(\Omega)})
(\|f\|_{B^{s+1}_{q,1}(\Omega)} + \|\CO_\lambda(\bg, \bh)\|_{B^s_{q,1}(\Omega)}).
$$  
And therefore,
setting $\rho=\lambda^{-1}(f-\eta_0\dv\bu)$, we see that 
$\rho$ and $\bu$ are solutions of equations \eqref{stokes.2}, and satisfy the estimate:
$$\|(\lambda\rho, \rho)\|_{B^{s+1}_{q,1}(\Omega)}
+ \|(\lambda, \bar\nabla^2)\bu\|_{B^s_{q,1}(\Omega)}
\leq  C(\rho_*, \|\tilde\eta_0\|_{B^{s+1}_{q,1}(\Omega)})
(\|f\|_{B^{s+1}_{q,1}(\Omega)} + \|\CO_\lambda(\bg, \bh)\|_{B^s_{q,1}(\Omega)}).
$$ 
The uniqueness of solutions
to equations \eqref{stokes.2} follows from the uniqueness of solutions of equations
\eqref{stokes.4}. 
\par 
In view of \eqref{sol.9.1},  we define an operator $\CZ_\Omega(\lambda)$ by 
\begin{equation}\label{sol.9.1*}
\CZ_\Omega(\lambda)(f, H) = \CW^2_\Omega(\lambda)H
+ \lambda^{-1}\CW^2_\Omega(\lambda)\CO_\lambda
(-\nabla(P'(\eta_0)f), P'(\eta_0)f\bn),
\end{equation}
then, we have $\bu = \CZ_\Omega(\lambda)(f, \CO(\bg, \bh))$. Therefore, 
we define an operator $\CC_m(\lambda)$ by
$$\CC_m(\lambda)(f, H) = \lambda^{-1}f - \lambda^{-1}\eta_0(x)
\dv(\CZ_\Omega(\lambda)(f, H)).
$$
Obviously, $\rho = \CC_m(\lambda)(f, \CO_\lambda(\bg, \bh))$. 

To define $\CB_m(\lambda)$ and $\CC_m(\lambda)$, we observe that
$\sum_{\ell=1}^\infty(\lambda^{-1}\CO_\lambda\CP(\lambda))^\ell
= \lambda^{-1}\CO_\lambda\CP(\lambda)\CQ_\infty(\lambda),  
$
where we have set
$\CQ_\infty(\lambda) 
= \sum_{\ell=0.}^\infty(\lambda^{-1}\CO_\lambda\CP(\lambda))^\ell.$
By \eqref{residue.4.2}, $\CQ_\infty$ is well-defined and 
$\| \CQ_\infty(\lambda)H\|_{B^s_{q,1}(\Omega)} \leq C\|H\|_{B^s_{q,1}(\Omega)}$
for any $\lambda \in \Sigma_{\epsilon, \lambda_4}$.  Moreover, 
$\tilde \CP_\infty(\lambda) = \bI +\lambda^{-1}\CQ_\infty(\lambda)$. 
Thus, 
in view of  \eqref{sol.9.1*}, we define
$\CB_v(\lambda)$ and $\CC_v(\lambda)$ by 
$$\CB_v(\lambda)H = \CV_\Omega(\lambda)H, 
\quad 
\CC_v(\lambda)(f, H) = \lambda^{-1}\CV_\Omega(\lambda)\CQ_\infty(\lambda)H 
+ \lambda^{-1}\CW^2_\Omega(\lambda)\CO_\lambda
(-\nabla(P'(\eta_0)f), P'(\eta_0)f\bn).
$$
Using Theorem \ref{thm:Lame}, we see that 
$\CB_v(\lambda)$ has $(s, \sigma, q)$ properties.  
Using Theorem \ref{thm:Lame}, the similar argument to
the proof of \eqref{deriv.4.4.1}, and \eqref{banach.1}, we see that
$$
\|(\lambda, \lambda^{1/2}\bar\nabla, \bar\nabla^2)\pd_\lambda^\ell
\CD(\lambda)(f, H)\|_{B^s_{q,1}(\Omega)}
\leq C|\lambda|^{-1-\ell}
\|(f, H)\|_{B^s_{q,1}(\Omega)\times \CH^s_{q,1}(\Omega)}
\quad(\ell=0,1)
$$
for $\CD \in \{\CC_m, \CC_v\}$. Namely, 
$\CC_m(\lambda)$ and $\CC_v(\lambda)$ have generalized resolvent properties
with $X = B^{s+1}_{q,1}(\Omega)\times \CH^s_{q,1}(\Omega)$ in the sense of 
Definition \ref{dfn.2}.  This completes the proof of Theorem \ref{thm:main4}.
\end{proof}

\section{$L_1$ maximal regularity, A proof of Theorem \ref{thm:main.1}.}\label{sec:6}

In this section, we consider the following evolution equations \eqref{Eq:Linear} 
and we shall prove Theorem \ref{thm:main.1}. 
To prove the theorem, problem \eqref{Eq:Linear} is divided into the following
two equations:
\begin{equation}\label{st.10.2}\left\{\begin{aligned}
\pd_t\rho +\eta_0\dv\bu = f&&\quad&\text{in $\Omega\times\BR$}, \\
\eta_0\pd_t\bu-\DV(\alpha\BD(\bu)+(\beta-\alpha)\dv\bu\BI
- P'(\eta_0)\rho\BI) = \bg&&\quad&\text{in $\Omega\times\BR$}, \\
(\alpha\BD(\bu) + (\beta-\alpha)\dv\bu\BI - P'(\eta_0)\rho\BI)\bn
= \bh&&\quad&\text{on $\pd\Omega\times\BR$};
\end{aligned}\right.\end{equation}
\begin{equation}\label{st.10.3}\left\{\begin{aligned}
\pd_t\rho +\eta_0\dv\bu = 0&&\quad&\text{in $\Omega\times(0, \infty)$}, \\
\eta_0\pd_t\bu-\DV(\alpha\BD(\bu)+(\beta-\alpha)\dv\bu\BI
- P'(\eta_0)\rho\BI) = 0&&\quad&\text{in $\Omega\times(0, \infty)$}, \\
(\alpha\BD(\bu) + (\beta-\alpha)\dv\bu\BI - P'(\eta_0)\rho\BI)\bn
= 0&&\quad&\text{on $\pd\Omega\times(0, \infty)$}, \\
(\rho, \bu)|_{t=0} = (\rho_0, \bu_0)&&\quad&\text{in $\Omega$}.
\end{aligned}\right.\end{equation}
First, we consider equations \eqref{st.10.2}.  We shall prove the following
theorem.
\begin{thm}\label{thm:L1.1}
Assume that the following conditions \thetag1 or \thetag2 holds. 
\begin{itemize}
\item[\thetag1] If $\eta_0(x) = \rho_*$, then $1 < q < \infty$ and $-1+1/q < s < 1/q$. 
\item[\thetag2] If $\tilde\eta_0(x)\not\equiv0$, then 
$N-1 < q < 2N$ and $-1+N/q \leq  s < 1/q$. 
\end{itemize}
 Let $\lambda_4$ be the positive number given in
Theorem \ref{thm:main4} and $\gamma \geq \lambda_4$. 
Then, for any 
right members $f$, $\bg$ and $\bh$ with
\begin{align*}
e^{-\gamma t}f \in L_1(\BR, B^{s+1}_{q,1}(\Omega)), \quad 
e^{-\gamma t}\bg \in L_1(\BR, B^s_{q,1}(\Omega)^N), \quad 
e^{-\gamma t}\bh \in L_1(\BR, B^{s+1}_{q,1}(\Omega)^N)
\cap W^{1/2}_1(\BR, B^s_{q,1}(\Omega)^N)
\end{align*}
for some $\gamma \geq \lambda_5$, problem \eqref{st.10.2} admits 
unique solutions $\rho$ and $\bu$ with
\begin{align*}
e^{-\gamma t}\rho \in W^1_1(\BR, B^{s+1}_{q,1}(\Omega)), \quad
e^{-\gamma t}\bu \in W^1_1(\BR, B^{s}_{q,1}(\Omega)^N) \cap 
L_1(\BR, B^{s+2}_{q,1}(\Omega)^N)
\end{align*}
possessing the estimate:
\begin{align*}
&\|e^{-\gamma t}(\pd_t\rho, \rho)\|_{L_1(\BR, B^{s+1}_{q,1}(\Omega))}
+ \|e^{-\gamma t}\pd_t\bu\|_{L_1(\BR, B^s_{q,1}(\Omega))}
+ \|e^{-\gamma t}\bu\|_{L_1(\BR, B^{s+2}_{q,1}(\Omega))}\\
&\leq C\|e^{-\gamma t}(\bar\nabla f,  \bg, \Lambda^{1/2}\bh, 
\nabla \bh)\|_{L_1(\BR, B^{s}_{q,1}(\Omega))}.
\end{align*}
Here, the constant $C$ depends on $\lambda_4$ but is independent 
of $\gamma$ when 
$\gamma\geq \lambda_4$.  
\end{thm}
\begin{proof}
First, we consider equations \eqref{st.10.2}.
Applying the Laplace transform to equations \eqref{st.10.2}, we have
\begin{equation*}
\begin{aligned}
\lambda u + \eta_0\dv\hat\bv = \hat f&&\quad&\text{in $\Omega$}, \\
\eta_0\lambda\bv -\DV(\alpha\BD(\bv)+(\beta-\alpha)\dv\bv\BI-
P'(\eta_0)u\BI) = \hat \bg&&\quad&\text{in $\Omega$}, \\
(\alpha\BD(\bv) + (\beta-\alpha)\dv\bv\BI - P'(\eta_0)u\BI)\bn
=\hat\bh&&\quad&\text{on $\pd\Omega$},
\end{aligned}\end{equation*}
where $\hat f = \CL[f]$ and $\hat\bg= \CL[\bg]$.  
From Theorem \ref{thm:main4},
it follows that $u= \CC_m(\lambda)(f, \CO_\lambda(\bg, \bh))$ and 
$\bv = \CB_v(\lambda)\CO_\lambda(\bg, \bh) + 
\CC_v(\lambda)(f, \CO_\lambda(\bg, \bh))$. 
We define $\rho$ and $\bu$ by
\begin{align*}
\rho & = \CL^{-1}[u] = \CL^{-1}[\CC_m(\lambda)(f, \CO_\lambda(\bg, \bh))], \\
\bu  &= \CL^{-1}[\bv] = \CL[\CB_v(\lambda)\CO_\lambda(\bg, \bh) + 
\CC_v(\lambda)(f, \CO_\lambda(\bg, \bh))].
\end{align*}  Then, 
$\rho$ and $\bu$ are solutions of equations \eqref{st.10.2}.
Moreover,  by Propositions \ref{prop:L1} and \ref{prop:L2}, we have 
\begin{align*}
\int^\infty_{-\infty}e^{-\gamma t}\|\rho(\cdot, t)\|_{B^{s+1}_{q,1}(\Omega)}\, dt 
& \leq C\int^{\infty}_{-\infty}e^{-\gamma t}(\|f(\cdot, t)\|_{B^s_{q,1}(\Omega)}
+ \|(\bg, \Lambda^{1/2}\bh, \nabla\bh)(t)\|_{B^s_{q,1}(\Omega)})\,dt, \\
\int^\infty_{-\infty}e^{-\gamma t}\|\bu(\cdot, t)\|_{B^{s+2}_{q,1}(\Omega)}\, dt 
& \leq C\int^{\infty}_{-\infty}e^{-\gamma t}(\|f(\cdot, t)\|_{B^s_{q,1}(\Omega)}
+ \|(\bg, \Lambda^{1/2}\bh, \nabla\bh)(t)\|_{B^s_{q,1}(\Omega)})\,dt.
\end{align*}
For the estimate of the time derivatives, we use equations, and then writing 
$\pd_t\rho(t) = -\rho_0\dv\bu + f$, we have
\begin{align*}
&\int^\infty_{-\infty}e^{-\gamma t}\|\pd_t\rho(t)\|_{B^{s+1}_{q,1}(\Omega)}\,dt \\
&\leq \int^\infty_{-\infty}e^{-\gamma t}(\|f(t)\|_{B^{s+1}_{q,1}(\Omega)}
+ C(\rho_*, \|\tilde\eta_0\|_{B^{N/q+1}_{q,1}(\Omega)})\|\bu\|_{B^{s+2}_{q,1}(\Omega)})\,dt
\\
&\leq C(\rho_*, \|\tilde\eta_0\|_{B^{N/q+1}_{q,1}(\Omega)})
\int^{\infty}_{-\infty}e^{-\gamma t}(\|f(t)\|_{B^s_{q,1}(\Omega)}
+ \|(\bg, \Lambda^{1/2}\bh, \nabla\bh)(t)\|_{B^s_{q,1}(\Omega)}\,dt.
\end{align*}
To estimate $\pd_t\bu$, we have to investigate the multiplication $\eta_0(x)^{-1}$.
From Assumption \ref{assump.1} we have the following lemma: 
\begin{lem}\label{lem:operation} Let $1 < q < \infty$ and 
$-1+N/q \leq s < 1/q$.  Let $\eta_0(x) = \rho_*+ \tilde\eta_0(x)$ such that 
$\tilde\eta_0(x) \in B^{N/q}_{q,1}(\Omega)$ and $\eta_0(x)$ satisifes
Assumption \ref{assump.1}. Then, for any $u \in B^s_{q,1}(\Omega)$, there holds 
\begin{equation*}
\|u\eta_0^{-1}\|_{B^s_{q,1}(\Omega)} \leq \rho_*^{-1}\|u\|_{B^s_{q,1}(\Omega)}
+ C\|\tilde\eta_0\|_{B^{N/q}_{q,1}(\Omega)}\|u\|_{B^s_{q,1}(\Omega)}
\end{equation*}
for some constant $C>0$ depending on $\rho_*$, $\rho_1$ and $\rho_2$.
\end{lem}
\begin{proof}Note that 
$\eta_0(x)^{-1} = \rho_*^{-1}-\tilde\eta_0(x)(\rho_*\eta_0(x))^{-1}$. 
If $q_1 > N$, then 
\begin{equation}\label{op.2}
\|\tilde\eta_0(x)(\rho_*\eta_0(x))^{-1}\|_{B^{N/q_1}_{q_1, 1}(\Omega)} 
\leq C\|\tilde\eta_0\|_{B^{N/q_1}_{q_1, 1}(\Omega)}.
\end{equation}
In fact,to prove \eqref{op.2},
we use the relation $B^{N/q_1}_{q_1,\infty}(\Omega) = (L_{q_1}(\Omega), 
W^1_{q_1}(\Omega))_{N/q_1, 1}$.
Since $\rho_1 < \eta_0(x) < \rho_2$ as follows from Assumption \ref{assump.1}, 
we have
$$
\|\tilde\eta_0(x)(\rho_*\eta_0(x))^{-1}\|_{L_{q_1}(\Omega)} \leq (\rho_*\rho_1)^{-1}
\|\tilde\eta_0\|_{L_{q_1}(\Omega)}.
$$
And also, 
$$\|\nabla (\tilde\eta_0(\rho_*\eta_0)^{-1})\|_{L_{q_1}(\Omega)}
\leq \|(\nabla \tilde\eta_0) (\rho_*\eta_0)^{-1}\|_{L_{q_1}(\Omega)}
+ \rho_*^{-1}\|\tilde\eta_0(\nabla \eta_0)\eta_0^{-2}\|_{L_{q_1}(\Omega)}.
$$
Noticing that $\nabla \eta_0=\nabla\tilde\eta_0$ and that 
$|\tilde\eta_0(x)| \leq |\eta_0(x)| + \rho_* \leq \rho_2+\rho_*$, we have
$$\|\nabla (\tilde\eta_0(\rho_*\eta_0)^{-1})\|_{L_{q_1}(\Omega)}
\leq ((\rho_*\rho_1)^{-1}+ \rho_*^{-1}(\rho_2+\rho_*)
\rho_1^{-2})\|\nabla\tilde\eta_0\|_{L_{q_1}(\Omega)}.
$$
Thus, there exists a constant $C$ depending 
on $\rho_*$, $\rho_1$ and $\rho_2$ such that 
\eqref{op.2} holds. \par
First, we consider the case where $N/q < 1$. Then, using Abidi-Paicu-Haspot estimate \eqref{p.prod1*} we have
\begin{align*}
\|u\eta_0^{-1}\|_{B^s_{q,1}(\Omega)} &\leq 
\rho_*^{-1}\|u\|_{B^s_{q,1}(\Omega)} 
+ \|\tilde\eta_0(\rho_*\eta_0)^{-1}\|_{B^{N/q}_{q, 1}(\Omega)}
\|u\|_{B^s_{q,1}(\Omega)} \\
&\leq (\rho_*^{-1}+ C\|\tilde\eta_0\|_{B^{N/q}_{q,1}(\Omega)})\|u\|_{B^s_{q,1}(\Omega)}.
\end{align*}
Next, we consider the case where $N/q \geq 1$.  Since $-1 + N/q \leq s < 1/q$, 
  if we choose $q_1$ in such a way that $N < q_1 < qN$, then 
$s \in (-N/q_1, N/q_1)$ and $s \in (-N/q', N/q_1)$. Thus, since $N/q_1 < 1$, 
using Abidi-Paicu-Haspot estimate \eqref{p.prod1*}
and \eqref{op.2}
we have
\begin{align*}
\|u\eta_0^{-1}\|_{B^s_{q,1}(\Omega)} &\leq (\rho_*^{-1}\|u\|_{B^s_{q,1}(\Omega)} 
+ \|\tilde\eta_0(\rho_*\eta_0)^{-1}\|_{B^{N/q_1}_{q_1, 1}(\Omega)}
\|u\|_{B^s_{q,1}(\Omega)}) \\
& \leq (\rho_*^{-1}\|u\|_{B^s_{q,1}(\Omega)} 
+ C\|\tilde\eta_0\|_{B^{N/q_1}_{q_1, 1}(\Omega)}
+ (\rho_*\rho_1)^{-1}(\rho_*+\rho_2))\|u\|_{B^s_{q,1}(\Omega)}.
\end{align*}
Notice that $1 < q \leq  N <  q_1$. 
By the embedding theorem of the Besov spaces, we have
\begin{align*}
\|\tilde\eta_0\|_{B^{N/q_1}_{q_1, 1}(\Omega)}
\leq C\|\tilde\eta_0\|_{B^{\frac{N}{q_1}+N\left(\frac{1}{q}-\frac{1}{q_1}\right)}_{q,1}(\Omega)}
=  C\|\tilde\eta_0\|_{B^{N/q}_{q,1}(\Omega)}.
\end{align*}
This completes the proof of Lemma \ref{lem:operation}.
\end{proof}
From equations \eqref{st.10.3}, we write
$$\pd_t\bu = \eta_0^{-1}(\DV(\alpha\BD(\bu) + (\beta-\alpha)\dv\bu\BI
- P'(\eta_0)\rho\BI) + \eta_0^{-1}\bg. 
$$
Using \eqref{op.2}, we have
\begin{align*}
&\int^\infty_{-\infty}e^{-\gamma t}\|\pd_t\bu(t)\|_{B^s_{q,1}(\Omega)}\,dt \\
&\quad \leq C(\rho_*, \|\tilde\eta_0\|_{B^{N/q+1}_{q,1}(\Omega)})
\int^\infty_{-\infty} e^{-\gamma t}
(\|\bu(t)\|_{B^{s+2}_{q,1}(\Omega)} + \|\rho(t)\|_{B^{s+1}_{q,1}(\Omega)})\,dt \\
&\quad  \leq C(\rho_*, \|\tilde\eta_0\|_{B^{N/q+1}_{q,1}(\Omega)})
\int^\infty_{-\infty}e^{-\gamma t}(\|f(\cdot, t)\|_{B^{s+1}_{q,1}(\Omega)}
+ \|(\bg, \Lambda^{1/2}\bh,  \nabla\bh)(t)\|_{B^s_{q,1}(\Omega)})\,dt.
\end{align*}
This completes the proof of Theorem \ref{thm:L1.1}. 
\end{proof}
We now consider equations \eqref{st.10.3}.  We shall show the generation of $C^0$ analytic 
semigroup $\{T(t)\}_{t\geq 0}$ defined on $\CH^s_{q,1}(\Omega)$.  
To formulate equations
\eqref{st.10.3} in the semigroup setting, we define a set 
 $\CD^\nu_{q,1}(\Omega)$  by setting
\begin{align*}
\CD^\nu_{q,1} & = \{(\rho, \bu) \in B^{\nu+1}_{q,1}(\Omega)\times B^{\nu+2}_{q,1}(\Omega) 
\mid (\alpha\BD(\bu) + (\beta-\alpha)\dv\bu \BI -P'(\eta_0)\rho\BI)\bn |_{\pd\Omega} 
= 0\}.
\end{align*}
Let $\CA_m$, $\bA_v$ and $\CA$ be an operators defined by 
\begin{equation*}
\begin{aligned}
\CA_m(\rho, \bu) &=\eta_0\dv\bu, 
\quad \CA_v(\rho, \bu) = \eta_0^{-1}(-\alpha\Delta\bu - \beta\nabla\dv\bu
+ \nabla(P'(\eta_0)\rho)), \\
\CA(\rho, \bu) &= (\CA_m(\rho, \bu), \CA_v(\rho, \bu))
\end{aligned}\end{equation*}
for $(\rho, \bu) \in \CD^\nu_{q,1}$.  Then, problem \eqref{st.10.3} is written as 
\begin{equation*}
(\pd_t + \CA)(\rho, \bu) = (0, 0) \quad \text{for $t > 0$}, \quad (\rho, \bu)|_{t=0} = (\rho_0, \bu_0)
\quad\text{in $\Omega$}.
\end{equation*}
From Theorem \ref{thm:main4}, $\Sigma_{\epsilon, \lambda_4}$ is 
contained in a resolvent set of 
the operator $\CA$ and for $\lambda \in \Sigma_{\epsilon, \lambda_4}$, we have
$$\|\lambda (\lambda\bI + \CA)^{-1}(f, \bg)\|_{\CH^s_{q,1}(\Omega)}
\leq C\|(f, \bg)\|_{\CH^s_{q,1}(\Omega)}.
$$
Therefore, there exists a $C_0$ analytic semigroup $\{T(t)\}_{t\geq 0}$
associated with equations \eqref{st.10.3}.  Moreover, we have the 
$L_1$ maximal regularity of $\{T(t)\}_{t\geq 0}$ as follows. 
\begin{thm}\label{thm:t.2}
Assume that the following conditions \thetag1 or \thetag2 holds. 
\begin{itemize}
\item[\thetag1] If $\eta_0(x) = \rho_*$, then $1 < q < \infty$ and $-1+1/q < s < 1/q$. 
\item[\thetag2] If $\tilde\eta_0(x)\not\equiv0$, then 
$N-1 < q < 2N$ and $-1+N/q \leq  s < 1/q$. 
\end{itemize}
Let $T(t)(\rho_0, \bu_0) = (T_m(t)(\rho_0, \bu_0), T_v(t)(\rho_0, \bu_0))$, that is
$T_m(t)$ and $T_v(t)$ are the mass density part and the velocity part, respectively. 
Then, there holds
\begin{equation}\label{L'1.1}
\int^\infty_0e^{-\gamma t}(\|(\pd_t, \bar\nabla^2)T_v(t)(f, \bg)\|_{B^s_{q,1}(\Omega)}
+ \|(1, \pd_t)T_m(t)(f, \bg)\|_{B^{s+1}_{q,1}(\Omega)})\, dt
\leq C\|(f, \bg)\|_{{\CH^s_{q,1}(\Omega)}}.
\end{equation}
\end{thm}
\begin{proof}
In view of Theorem \ref{thm:main4}, we can write
$(\lambda \bI + \CA)^{-1}(f, \bg) = (\CC_m(\lambda)(f, \CO_\lambda(\bg, 0)),
\CB_v(\lambda)\CO_\lambda(\bg, 0) + \CC_v(\lambda)(f, \CO_\lambda(\bg, 0))$. 
As is known well, $T(t)(f, \bg) = \CL^{-1}[(\lambda\bI + \CA)^{-1}(f, \bg)]$. 
By Proposition \ref{prop:L1} and \ref{prop:L2} we see that \eqref{L'1.1}  hold.
This completes the proof of Theorem \ref{thm:t.2}. 
\end{proof}
Let $\rho$ and $\bu$ be solutions given in Theorem \ref{thm:L1.1}.
By the time trace theorem, we see that
\begin{multline*}
\|(\rho, \bu)|_{t=0}\|_{B^{s+1}_{q,1}(\Omega)\times B^s_{q,1}(\Omega)} \\
\leq C(\|e^{-\gamma t}(\pd_t\rho, \rho)\|_{L_1((0, \infty), B^{s+1}_{q,1}(\Omega))}
+ \|e^{-\gamma t}\bu\|_{L_1((0, \infty), B^{s+2}_{q,1}(\Omega))}
+ \|e^{-\gamma t}\pd_t\bu\|_{L_1((0, \infty), B^s_{q,1}(\Omega))}).
\end{multline*}
Thus, combining Theorems \ref{thm:L1.1} and \ref{thm:t.2}, we can prove
the existence part of Theorem \ref{thm:main.1}. The uniqueness of solutions to
 equations \eqref{Eq:Linear} can be proved using Dhamel's principle and 
the analytic semigroup $\{T(t)\}_{t\geq0}$. 
This completes the proof of Theorem \ref{thm:main.1}.


\end{document}